\documentclass[a4paper,10pt]{article}

\usepackage{a4wide}
\usepackage{multicol, float}
\usepackage{lineno}
\usepackage{verbatim, euscript}
\usepackage{subfigure}
\usepackage{bbm}

\usepackage{pictexwd}

\usepackage{psfrag}

\usepackage{tikz}
\usetikzlibrary{arrows.meta}
\usepackage{rotating}
\usepackage{color}
\usepackage{mathrsfs, stmaryrd}
\usepackage{pictexwd}
\usepackage{tikz}
\usepackage{etex} 
\usepackage{multido}

\usepackage{lineno,hyperref}
\usepackage{verbatim, euscript}
\usepackage{subfigure}
\usepackage{bbm,amssymb}
\usepackage{amsfonts}
\usepackage{amsmath, mathtools}
\usepackage{pictexwd}
\usepackage{amsthm}
\usepackage{graphicx}

\usepackage{tikz}
\usetikzlibrary{arrows.meta}
\usetikzlibrary{math}
\usetikzlibrary{decorations.pathreplacing}

\usepackage{rotating}
\usepackage{color}
\usepackage{mathrsfs}

\usepackage{mathtools}
\mathtoolsset{showonlyrefs}

\usetikzlibrary{arrows,automata}

\DeclarePairedDelimiter\floor{\lfloor}{\rfloor}

\numberwithin{equation}{section}

\newcommand{\eps}{\varepsilon}

\theoremstyle{plain}
\newtheorem{The}{Theorem}[section]
\newtheorem{Lem}[The]{Lemma}

\newtheorem{Def}[The]{Definition}
\newtheorem{remark}[The]{Remark}

\newtheorem{Prop}[The]{Proposition}
\newtheorem{Ass}[The]{Assumption}

\newcommand{\E}{\mathbb{E}}
\newcommand{\PP}{\mathbb{P}}
\newcommand{\1}{\mathbbm{1}} 
\newcommand{\N}{\mathbb{N}}
\newcommand{\Z}{\mathbb{Z}}
\newcommand{\R}{\mathbb{R}}

\usepackage[utf8]{inputenc}

\newcommand{\mP}{\mathbb{P}}

\newcommand{\V}{\mathbb{V}}
\newcommand{\ZZ}{\mathbf{Z}}

\begin{document}

\title{Spatial Invasion of Cooperative Parasites}
\author{
Vianney Brouard\thanks{ENS de Lyon, UMPA, CNRS UMR 5669, 46 Allée d’Italie, 69364 Lyon Cedex 07, France; E-mail: vianney.brouard@ens-lyon.fr},~
Cornelia Pokalyuk\thanks{Mathematical Institute, Goethe-University Frankfurt, 60629, Frankfurt am Main, Germany; E-mail: pokalyuk@math.uni-frankfurt.de},~
Marco Seiler\thanks{Frankfurt Institute for Advanced Studies, Ruth-Moufang-Straße 1, 60438 Frankfurt am Main, Germany; E-mail: seiler@fias.uni-frankfurt.de} \ and
Hung Tran\thanks{Frankfurt Institute for Advanced Studies, Ruth-Moufang-Straße 1, 60438 Frankfurt am Main, Germany; E-mail: htran@ae.cs.uni-frankfurt.de}}

	\maketitle

\begin{abstract}
    In this paper we study invasion probabilities and invasion times of cooperative parasites spreading in spatially structured host populations. The spatial structure of the host population is given by a random geometric graph on $[0,1]^n$, $n\in \mathbb{N}$, with a Poisson($N$)-distributed number of vertices and in which vertices are connected over an edge when they have a distance of at most $r_N\in \Theta\left(N^{\frac{\beta-1}{n}}\right)$ for some $0<\beta<1$ and $N\rightarrow \infty$. At a host infection many parasites are generated and parasites move along edges to neighbouring hosts. We assume that parasites have to cooperate to infect hosts, in the sense that at least two parasites need to attack a host simultaneously. We find lower and upper bounds on the invasion probability of the parasites in terms of survival probabilities of branching processes with cooperation. Furthermore, we characterize the asymptotic invasion time.
    
    An important ingredient of the proofs is a comparison with infection dynamics of cooperative parasites in host populations structured according to a complete graph, i.e. in well-mixed host populations. For these infection processes we can show that invasion probabilities are asymptotically equal to survival probabilities of branching processes with cooperation. 
    Furthermore, we build in the proofs on techniques developed in \cite{BrouardEtAl2022}, where an analogous invasion process has been studied for host populations structured according to a configuration model. 
    
    We substantiate our results with simulations.  
\end{abstract}


\section{Introduction}

Understanding the dynamics of infection processes is a highly relevant and active research field. In this study we are particularly interested in the spread of cooperative parasites in spatially structured host populations. Cooperative behaviour is observed in many biological systems, see \cite{Rubenstein2010}.
The main biological motivation for our model stems from observations made on phages, that is viruses infecting bacteria. Bacteria own various mechanisms to defend against phages. Defense on the basis of CRISPR-Cas system is widespread in bacteria. Certain phages, called anti-CRISPR phages, can overcome these defense mechanism by cooperation. Only when anti-CRISPR-phages infect simultaneously or subsequently a CRISPR-resistant bacterium the infection gets likely to be successful, see \cite{Landsberger, Borges}.

Besides the motivation stemming from application, models which incorporate cooperative mechanisms are also highly interesting from a mathematical point. 
For example Gonzalez Casanova, Pardo and Perez \cite{GPP21} show that for a branching process with cooperation the survival probability is positive as long as the probability to generate offspring for pairs of individuals is non-zero.  In case of survival it explodes in finite time. In the papers \cite{neuhauser1994long}, \cite{noble1992equilibrium} and more recently \cite{mach2020recursive} mean-field limits of systems with cooperative reproduction are studied. Mach et al. find in \cite{mach2020recursive} that the mean-field equation corresponding to certain interacting particles with cooperation can have more fixed points than the corresponding  mean-field equations of classical infection models such as the contact process. This can be seen as evidence that in the microscopic model there could exist more extremal invariant laws as compared to the non-cooperative infection models. Sturm and Swart  studied in \cite{sturm2015particle} such a cooperative microscopic model. To be precise they considered a nearest neighbour cooperative branching-coalescing random walk on $\Z$. In comparison with the classical branching-coalescing random walk a subcritical phase exists, where the system ends up with only one particle. Superficially, this cooperative branching-coalescing system seems to be similar to a contact process, but a closer look reveals some apparent differences. For example \cite{sturm2015particle} show that the decay rates in the subcritical regime are polynomial and not exponentially as for the contact process.

In \cite{BrouardEtAl2022} the invasion of cooperative parasites in host populations structured according to a configuration model was studied. In this paper a parameter regime was considered, in which 
parasites have many offspring and a parasite can reach many, but not all hosts. In the critical scale this resulted in an invasion process, which could be approximated initially by a Galton Watson process with roughly Poisson offspring numbers. We show in \cite{BrouardEtAl2022}, that the invasion probability is asymptotically equal to the survival probability of this approximating Galton-Watson process.

In this manuscript we consider the spread of cooperative parasites in host populations that have a finite-dimensional spatial structure. More precisely, we assume that the (immobile) hosts are distributed on an $n$-dimensional cube $[0,1]^n$ according to a Poisson point process. Parasites can move in every generation up to some fixed distance in space and attack the hosts located in this region. As in \cite{BrouardEtAl2022} we consider a parameter regime, in which parasites have many offspring and can reach many, but not all hosts, as well as hosts need to be attacked jointly by parasites for successful parasite reproduction.

However, in contrast to the case considered in \cite{BrouardEtAl2022} the invasion process is already in the initial phase badly approximated by a Galton-Watson process. The reason is that parasites generated in different hosts have in the spatial setting often a good chance to cooperate, because infected hosts are located close to each other. To arrive at lower and upper bounds on the invasion probability we compare the spread with infection dynamics caused by
 cooperative parasites spreading on complete graphs. The number of vertices of these complete graphs yield upper and lower bounds on the number of hosts, with which parasites generated on different hosts can cooperate. We show that the invasion probabilities of these infection processes on complete graphs are asymptotically equal to survival probabilities of certain branching processes with cooperation, a result that is of interest on its own. Furthermore, we show that the spatial infection processes can be coupled from above and below  with these branching processes with cooperation until either the parasite population dies or a sufficiently large amount of hosts are infected so that afterwards with high probability (i.e. asymptotically with probability 1, abbreviated as whp in the following) the parasite population will spread through the whole host population. 

Once a sufficiently large number of hosts is infected, we show that the parasite population spreads with high probability at linear, almost maximal speed. As in the considered scaling the initial phase, which is decisive for survival of the parasite population, takes place only on a negligible amount of space, invasion time is basically determined by the time frame in which the parasite population spreads linearly fast. This yields our asymptotic result on the invasion time. Here again a clear difference to the model in \cite{BrouardEtAl2022} occurs, in which the final phase of invasion is finished after a constant number of steps.

By means of simulations we study the fit of the upper and lower bound on the invasion probability and our prediction for the invasion time. Interestingly, we find that the upper bound on the invasion probability matches very well with simulated invasion probabilities.

\section{Main results}

\subsection{Model definition and main theorems}

Consider the $n$-dimensional cube $[0,1]^n$, which we will denote by $M= M_n$ in the following. Measure distances on $M$ according to the maximum metric denoted by $\rho$, i.e. for $x=(x_1, ..., x_n), y=(y_1, ..., y_n) \in [0,1]^n$ we have
$\rho(x,y) = \max_{i=1, ..., n}\{|x_i - y_i|\}$.
 Consider a homogeneous Poisson point process  with intensity $N$ on $[0,1]^n$, in particular the number of Poisson points contained in a set $S \subset [0,1]^n$ of volume $s$ is Poisson distributed with parameter $sN$. Denote the set of the Poisson points by $\mathcal{V}=\mathcal{V}^{(N)}$. Build a random geometric graph on $[0,1]^n$ by connecting all points in $ \mathcal{V}^{(N)}$  over an edge which have a distance of at most $r_N$ with respect to the metric $\rho$. Denote the set of edges by $\mathcal{E}=\mathcal{E}^{(N)}$ and the random geometric graph by $\mathcal{G}=\mathcal{G}^{(N)}=\mathcal{G}^{(r_N)}= (\mathcal{V}^{(N)}, \mathcal{E}^{(N)})$. 

On $\mathcal{G}$ we consider the following infection process. At the beginning place on each vertex a host. These hosts can get infected with parasites. Choose the vertex $x_0= x_0^{(N)}  \in \mathcal{V}$ closest to the center of the cube $[0,1]^n$. We assume that the host on this vertex gets infected in the first generation $g=0$. This means that the host dies and $v_N$ offspring parasites are generated on $x_0$. Then the infection process continues in discrete generations. At the beginning of each generation each parasite chooses uniformly at random and independently of all other parasites  an edge, that is adjacent to the vertex on which the parasite is located. Along this edge the parasite moves to the neighbouring vertex and attacks the host on this vertex, if a host is still available. After movement of parasites, offspring parasites are generated and hosts die according to the following rules. If a vertex is occupied by a host and at least two parasites attack the host, the host on the vertex gets infected, dies and $v_N$ parasites are generated. If only a single parasite attacks a host, it dies and the host stays alive. If a parasite arrives at an unoccupied vertex, it dies.

If a vertex is occupied/not occupied with a host, in the following we will call these vertices \textit{occupied/unoccupied vertices}. Sometimes we also speak of \textit{susceptible/so far uninfected vertices},
if a host on a vertex did not yet get infected. Similar, we say that a \textit{vertex} is \textit{infected}, if the host on the vertex is in the current generation infected.

Denote by $\mathcal{S}_g= \mathcal{S}_g^{(N)} \in \mathcal{V}$, $\mathcal{I}_g= \mathcal{I}_g^{(N)}$ and $\mathcal{R}_g= \mathcal{R}_g^{(N)}$, resp., the occupied and uninfected, the infected and the unoccupied, resp., vertices in generation $g$. We set $S_g:= |\mathcal{S}_g|$, $ I_g:= |\mathcal{I}_g|$ and $R_g:=|\mathcal{R}_g|$.  Furthermore  $\overline{I}^{(N)}_g= \sum_{i=0}^{g} I^{(N)}_i$ is the number hosts that got infected till generation $g$. Let $\mathbf{I}= \mathbf{I}^{(N)} = (I^{(N)}_g)_{g\geq0}$  and  $\overline{\mathbf{I}}= \overline{\mathbf{I}}^{(N)}=(\overline{I}^{(N)}_g)_{g\geq 0}$ be the corresponding processes. 

To state our main results about the infection process we need the definition of \textit{branching processes with cooperation in discrete time}.  

\begin{Def}[Branching process with cooperation in discrete time]
  Let $\mathcal{L}_o$ and $\mathcal{L}_c$ be two probability distributions on $\mathbbm{N}_0$. A discrete-time branching process with cooperation (DBPC, for short) $\mathbf{Z}=(Z_g)_{g\geq 0}$ with offspring distribution $\mathcal{L}_o$ and cooperation distribution $\mathcal{L}_c$ is recursively defined as follows. Assume $Z_0=k$  a.s. for some $k\in \mathbbm{N}$, then for any $g\geq1$ $Z_g$ is defined as
 \begin{equation*}
 Z_g := \sum_{i=1}^{Z_{g-1}} X_{g,i} + \sum_{i,j=1, i>j}^{Z_{g-1}} Y_{g,i,j},
\end{equation*}
where $(X_{g,i})_{g,i}$ and $(Y_{g,i,j})_{g,i,j}$ are sequences of independent and identically distributed random variables with $X_{g,i}\sim \mathcal{L}_o$ and 
$Y_{g,i,j} \sim \mathcal{L}_c$.
We denote by $\overline{\mathbf{Z}}= (\overline{Z}_g)_{g\geq 0}$ the total size process, i.e. 
\[\overline{Z}_g= \sum_{i=0}^g Z_i.\]
In the following we will denote the probability weights of the distributions $\mathcal{L}_o$ and  $\mathcal{L}_c$, resp.,  by  $(p_{k,o})_{k\in \mathbbm{N}_0}$ and $(p_{k,c})_{k\in \mathbbm{N}_0}$, resp.
\end{Def}

\begin{remark}
Branching processes with cooperation have been mostly studied in continuous time in more general settings, like branching process with (pairwise) interactions, see e.g. \cite{S49}, \cite{K02},\cite{K03}, \cite{GPP21}, \cite{OP20} and  \cite{CLCZ12}. In particular, in \cite{K02} formula for extinction probabilities for the case of branching processes with cooperation have been determined, see Remark \ref{Rem:FixProbabilities}
\end{remark}

A central object for our results is a DBPC with Poisson offspring and cooperation distributions or rather its survival probability. Therefore, we fix in the next definition some notation for these processes.
\begin{Def}[DBPC with Poisson offspring and cooperation distributions]\label{DefPoisDBPC}
Let $a>0$. Denote by $\mathbf{P}^{(a)}$ a DBPC with offspring distribution $\mathcal{L}_o\sim \text{Poi}\big(\frac{a^2}{2}\big)$ and the cooperation distribution $\mathcal{L}_a\sim\text{Poi}(a^2)$. 
Furthermore, we denote by $\pi(a)$ the survival probability of  $\mathbf{P}^{(a)}$.
\end{Def}

Denote by 
$$d_N:= (2r_N)^n N,$$
which is the expected number of vertices a vertex of $\mathcal{G}^{(N)}$ (with an asymptotically non-vanishing distance to the boundary of $M$) is connected to in dimension $n$.
Furthermore, 
denote by $$E_{u}^{(N)}:= \left\{\exists g\in \mathbbm{N}_0: \overline{I}^{(N)}_g \geq u \cdot \# \mathcal{V}^{(N)}  \right\}$$
the event that at least a proportion $u$ of the host population dies during the infection process.

Our main result is the following theorem.

\begin{The}\label{MainResult}
Consider the above defined sequence of infection processes $\left(\mathbf{I}^{(N)}\right)_{N \in \mathbbm{N}}$ on $[0,1]^n$ for some $n \in \mathbbm{N}$. Assume $r_N= \frac{1}{2} N^{\frac{\beta-1}{n}}$ for some $0< \beta <1$ and let $0<u\leq 1$.  Then it holds:
\begin{itemize}
\item[1)] Invasion probability
     \begin{itemize}
    \item[(i)] If $v_N \in o(\sqrt{d_N})$
 \begin{align*}
     \lim_{N\rightarrow \infty} \PP\left(E_u^{(N)}\right)=0.
     \end{align*}
\item[(ii)] If $v_N \sim a\sqrt{d_N}$ for some $0<a<\infty$
     \begin{align*}
     \pi\left(\tfrac{a}{\sqrt{2^{n}}}\right)\leq \liminf_{N\rightarrow \infty} \PP\left(E_u^{(N)}\right) \leq \limsup_{N\rightarrow \infty} \PP\left(E_u^{(N)}\right)\leq \pi(a).
     \end{align*}
     \item[(iii)] If $\sqrt{d_N}\in o(v_N)$
     \begin{align*}
     \lim_{N\rightarrow \infty} \PP\left(E_u^{(N)}\right)=1.
     \end{align*}
     \end{itemize}
\item[2)] Invasion time\\ Assume $v_N \sim a \sqrt{d_N}$ for some $0<a<\infty$.\\
Denote by $$T^{(N)}:= \inf\left\{g \in \mathbb{N} | \overline{I}_g^{(N)}=\# \mathcal{V}^{(N)}\right\}.$$
Then
\begin{align*}
\lim_{N\rightarrow \infty}\PP\left(\left\lfloor\frac{1}{2r_N} \right\rfloor \leq T^{(N)} \leq \left\lceil\frac{1}{2 r_N} \right\rceil + \mathcal{O}\left(\max \left(\log(\log(N)),\frac{\varepsilon_N}{r_N^2}\right)\right) \bigg| T^{(N)}<\infty\right)=1.
\end{align*}
with $\varepsilon_N =\left(N^{\frac{\frac{\beta}{2}-1}{n} + \delta}\right)$, for any $\delta>0$.
\end{itemize}
\end{The}

\begin{remark}\label{Rem:FixProbabilities}
\begin{itemize}
    \item[(i)] In 1) (ii) we obtain bounds for $\liminf \mP(E_u^{(N)})$ and $\limsup \mP\left(E_u^{(N)}\right)$. We believe that the limit of $\mP\left(E_u^{(N)}\right)$ exists. Simulations suggest that the upper bound provides a good approximation of the actual invasion probability, see Section \ref{Sec: Simulations}. An analysis of the initial phase of the epidemic, when infected parasites start to spread around the initially infected vertex, would be helpful to understand if the upper bound indeed gives the correct asymptotic. In case this is true, immediately also the existence of the limit  $\mP\left(E_u^{(N)}\right)$ for $N\rightarrow \infty$ would follow.
  \item[(ii)]  In Lemma \ref{lem:SurvivalDBPC} below, we prove, that for any $a>0$ the survival probability $\pi\left(a \right)$ of a DBPC $\mathbf{P}^{\left(a\right)}$ as defined in Definition \ref{DefPoisDBPC} is strictly positive. Therefore, the invasion probability is in Case 1) (ii) of Theorem \ref{MainResult} for any $a>0$ strictly positive. This contrasts the situation studied in \cite{BrouardEtAl2022} where for $\frac{a^2}{2}\leq 1$ the invasion probability is asymptotically 0 (for a host population structured according to a configuration model instead of a random geometric graph on $[0,1]^n$).
  \item[(iii)] Extinction probabilities of branching processes with cooperation in continuous time have been characterised in certain cases in \cite{K03}. In particular, for branching processes with cooperation in continuous time with offspring and cooperation events occurring at the same rate and a $\text{Poi}(\mu)$ and $\text{Poi}(\lambda)$, respectively, offspring and cooperation distribution the extinction probability for a process started in 1 solves the equation
\begin{align*}
q= \int_K \frac{u}{ \phi(u)} du.
  \end{align*}
   Here \[\phi(u)= \frac{1}{h_2(u)-u^2} \exp\left(\int_{u_0}^u\frac{h_1(v)-v}{h_2(v)-v^2} dv\right),\]
   $h_1$ is the generating function of a $\text{Poi}(\mu)$ distributed random variable and $h_2$ the generating function of the random variable $Y= X+2$ where $X\sim\text{Poi}(\lambda)$. Furthermore, $K$ is a curve in the complex u-plane  meeting the condition
  \begin{align*}
  \int_K \frac{d}{du} e^{zu}(h_2(u)-u^2)\phi(u) du=0.
  \end{align*}
and $u_0$ is the starting point of $K$. 
\item[(iv)] We assume in our model that a parasite dies, if it moves to an unoccupied vertex. This assumption is not essential, it just simplifies some proofs. Our results still hold true if one e.g. assumes that a parasite, which moves to an unoccupied vertex, stays alive and moves forward in the next generation.
    \end{itemize}
\end{remark}

\begin{remark}\label{remarkManifolds}(Host populations structured according to random geometric graphs on Riemannian manifolds)
Instead of considering the spread of the parasite population in host populations structured according to a random geometric graph on an $n$-dimensional cube it is natural to assume that the host population is located on a manifold. We can generalize our model to this setting as follows. Let $(M',g)$ be a compact, connected orientable, n-dimensional Riemannian manifold with Riemannian metric $g$. Assume without loss of generality that $vol(M')=1$, where $vol(M')$ denotes the volume of $M'$ calculated according to the volume from induced by $g$. Denote furthermore by $\rho'$ the metric on $M'$ induced by $g$. Consider a homogeneous Poisson point process  with intensity $N$ on $M'$ (for this point process the number of vertices contained in a set $S\subset M'$ with volume $vol(S)=s$ is Poisson distributed with parameter $s N$). We denote the set of the Poisson points by $\mathcal{V}'=\mathcal{V'}^{(N)}$ and build a random geometric graph on $M'$ by connecting all points in $ \mathcal{V'}^{(N)}$  over an edge which have a distance of at most $r_N$ with respect to the metric $\rho'$. Denote the set of edges by $\mathcal{E'}=\mathcal{E'}^{(N)}$ and the random geometric graph by $\mathcal{G'}=\mathcal{G'}^{(N)}=\mathcal{G'}^{(r_N)}= (\mathcal{V'}^{(N)}, \mathcal{E'}^{(N)})$. 

Given $\mathcal{G'}$ we can consider an infection process (with the components $(\mathcal{S}'_g, \mathcal{I}'_g, \mathcal{R}'_g)$) in the same way in which we defined it on the random graph on the cube.

  Denote by $$d'_N:=\frac{\pi^{n/2}}{\Gamma\left(\frac{n}{2} \right)} \left(r_N \right)^n N,$$ which is the expected number of vertices a vertex of $\mathcal{G}^{(N)}$ is connected to in dimension $n$ (if the distance of the vertex to the boundary of $M'$ is asymptotically non-vanishing, in case $M'$ has a boundary) and let $p\in M'$. Denote by 
\begin{align}
    \tau(p) := \max_{q\in M} \{\rho'(q,p)\} 
\end{align}
the maximal distance between $p$ and any other point $q\in M$. Furthermore, 
denote as before by $${E'}_{u}^{(N)}:= \left\{\exists g\in \mathbbm{N}_0: \overline{I}^{(N)}_g \geq u \cdot \# \mathcal{V'}^{(N)}  \right\}.$$
Then we believe that the following statements hold at least for $n\in \{1,2\}$.

 Assume that $r_N = \frac{\Gamma\left(\frac{n}{2} \right)}{\pi^{n/2}}  N^{\frac{\beta-1}{n}}$ for some $0<\beta<1$ and $r>0$, 
   let $0<u\leq 1.$ Assume the infection process is started in a vertex $x_0^{(N)}\in \mathcal{V'}^{(N)}$ that has asymptotically a positive distance to the boundary of $M'$ (if $M'$ has a boundary). 
\begin{itemize}
\item[1)] Invasion probability
     \begin{itemize}
    \item[(i)] If $v_N \in o(\sqrt{d'_N})$
 \begin{align*}
     \lim_{N\rightarrow \infty} \PP\left(E_u^{(N)}\right)=0.
     \end{align*}
\item[(ii)] 
    If $v_N \sim a \sqrt{d'_N}$ for some $0<a<\infty$
     \begin{align*}
     \pi\left(\tfrac{a}{\sqrt{2^n}}\right)\leq \liminf_{N\rightarrow \infty}\PP\left(E_u^{(N)}\right) \leq \limsup_{N\rightarrow \infty} \PP\left(E_u^{(N)}\right)\leq \pi(a).
     \end{align*}
     \item[(iii)] If $\sqrt{d'_N}\in o(v_n)$
     \begin{align*}
     \lim_{N\rightarrow \infty} \PP\left(E_u^{(N)}\right)=1.
     \end{align*}
     \end{itemize}
\item[2)] Invasion time: Assume $v_N \sim a \sqrt{d'_N}$ for some $0<a<\infty$.\\
Denote by ${T'}^{(N)}:= \inf\left\{g \in \mathbbm{N} | \overline{I}_g^{(N)}=\# \mathcal{V'}^{(N)}\right\}$.
Then 
\begin{align*}
&\PP\left( \left\lfloor \frac{\tau\left(x_0^{(N)}\right)}{r_N}\right\rfloor\leq {T'}^{(N)} \leq  \left\lceil \frac{\tau\left(x_0^{(N)}\right)}{r_N}\right\rceil+\mathcal{O} \left(\max \left( \log(\log(N)), \frac{\varepsilon_N}{r_N^2} \right)\right)\bigg| {T'}^{(N)}<\infty\right) \\
&\underset{N \to \infty}{\to}1,
\end{align*}
with $\varepsilon_N =\left(N^{\frac{\frac{\beta}{2}-1}{n} + \delta}\right)$, for any $\delta>0$.
\end{itemize}
The main reason that these results should hold is that the decision if eventually invasion takes place is failed in a neighbourhood of $x_0$ that has asymptotically a negligible volume, since only $N^{\varepsilon}$ many hosts need to get infected to show that whp subsequently the whole host population gets infected and $N^{\varepsilon}$ many hosts are directly connected to $x_0$ for $\varepsilon>0$ small enough. Therefore, at the beginning the invasion process is essentially the same as a corresponding process on $[0,1]^n$ with distances measured according to the Euclidean distance. Indeed, for any sequence  $h_N \rightarrow 0$ it holds $\frac{vol(B_{h_N}(x))}{vol(\tilde{B}_{h_N}(0))} = 1- \frac{S}{6(n+2)} h_N^2 + o(h_N^2)$, where $vol(B_{h_N}(x))$ denotes the volume of a (geodesic) ball of radius $h_N$ centered in $x \in M'$ and $vol(\tilde{B}_{h_N}(0))$ denotes the volume of an $n$-dimensional Euclidean ball of radius $h_N$ centered in 0 and $S$ the scalar curvature in $x$, see \cite{Chavel84}, Section XII.8. Since $M'$ is compact and scalar curvature is a continuous function on $M'$, scalar curvature of $M'$ is bounded from above and below. In particular, for $h_N= r_N$  the number of points connected to $x\in \mathcal{V'}$ is Poisson distributed with parameter $N\frac{\pi^{n/2}}{\Gamma(\frac{n}{2} +1)} r_N^n + O(N r_N^{n+2})$, since  $vol(\tilde{B}_{r_N}(0)) = \frac{\pi^{n/2}}{\Gamma(\frac{n}{2} +1)} r_N^n$.

In Theorem \ref{MainResult} we consider the maximum metric to measure distances between two points. With this metric we easily can cover $M$ with balls (that are cubes as well) to control the spread of parasites across $M$. A similar construction is also possible with Euclidean balls (at least in the case $n\in \{1,2\}$), the notation is just a bit more complicated. Therefore, considering the Euclidean metric or maximum metric should not influence the invasion probability as long as the ratio of the expected number of vertices contained in a ball and the number of offspring parasites generated at infection is asymptotically the same. The invasion time in general differs for two different metrics, because the function $\tau(p)$ depends on the metric.
\end{remark}

Next we want to give a sketch of the proof of  Theorem~\ref{MainResult}, which is formally proven in Section \ref{Sec:RGG}.  
The proof of the lower bound on the invasion probability is based on an asymptotic result on the invasion probability of an analogously defined infection process when the host population is not structured according to a random geometric graph on the cube, but according to a complete graph. This model mimics the spread of cooperative parasite in well-mixed host populations and is neither covered by the parameter regime considered in \cite{BrouardEtAl2022} nor by Theorem \ref{MainResult}. Therefore the result is of interest on its own. We state it next.

Consider a complete graph with $D_N$ vertices. 
On the complete graph we consider the same infection process as on the random geometric graph. We assume that at infection $V_N$ many parasites are generated. As in the case of the random geometric graph we count the number of infected hosts up to generation $g$, that we denote here by $\overline{J}^{(N)}_g$, and we are interested in the event $F_u^{(N)}$ that eventually a proportion $u$ of the host population gets infected, i.e.
\begin{equation*}
    F_u^{(N)} =\big\{\exists g\in \mathbbm{N}_0: \overline{J}^{(N)}_g \geq u \cdot D_N  \big\}.
\end{equation*}
We show that the invasion probability is in the critical regime asymptotically equal to the survival probability of a branching process with cooperation.

\begin{The}\label{Theorem:Invasion on complete graph}
Assume $D_N \in \Theta \left( N^{\beta}\right)$ for some $0< \beta<1$.  The following invasion regimes hold: \\
(i) Assume $V_N \in o \left( \sqrt{D_N} \right)$. Then for all $0<u \leq 1$ 
\begin{align}
    \lim_{N \to \infty} \mathbb{P} \left( F_{u}^{(N)}\right)=0.
\end{align}
(ii) Assume $V_N \sim a \sqrt{D_N}$ for some $0<a<\infty$. 
Then the invasion probability of parasites satisfies for all $0<u \leq 1$ 
\begin{align}
    \lim_{N \to \infty} \mathbb{P}\left(F_{u}^{(N)}\right)=\pi(a)>0.
\end{align}
(iii) Assume $\sqrt{D_N} \in o(V_N)$. Then 
\begin{align}
    \lim_{N \to \infty} \mathbb{P}\left(F_{u}^{(N)}\right)=1.
\end{align}
\end{The}
The proof of Theorem \ref{Theorem:Invasion on complete graph} is given in Section \ref{Section: Infection on complete graph}. Next we sketch the proofs of Theorem \ref{MainResult} and \ref{Theorem:Invasion on complete graph}.

Hereinafter we often will use the following terminology.
We call an infection a \textit{CoSame} infection (for cooperation from the same edge), if a host gets infected by two parasites (originating from the same vertex) that moved along the same edge to the vertex on which the infected host is located on, and we call an infection a \textit{CoDiff} infection (for cooperation from different edges), if a host gets infected by two parasites that moved along different edges to the vertex the infected host is located on. \\

\textit{Sketch of the proof of Theorem \ref{Theorem:Invasion on complete graph}}:
Case (ii):
To arrive at an upper bound on the invasion probability we couple whp the total number of currently infected and currently empty vertices $\overline{\mathcal{J}}$ from above with the total size of a DBPC $\overline{\mathbf{Z}}^{(N)}_u$ until $\overline{\mathbf{Z}}^{(N)}_u$ remains constant or reaches at least the level $\ell_N$ for a sequence $\ell_N$ with $\ell_N \rightarrow \infty $ sufficiently slowly, 
see Proposition \ref{Prop: Coupling with upper DBPC}.  The probability to reach with $\overline{\mathbf{Z}}^{(N)}_u$ the level $\ell_N$ is asymptotically equal to $\pi(a)$, see Proposition \ref{Proposition: asymptotic probability reaching certain level for BPI}, as the appproximating DBPC has asymptotically the survival probability $\pi(a)$. In case the level $\ell_N$ is reached we upper bound the probability  by 1, that afterwards also the remaining hosts get infected.\\
For the lower bound we couple whp $\overline{\mathcal{J}}$ from below with a DBPC $\overline{\mathbf{Z}}^{(N)}_\ell$   that has asymptotically the survival probability $\pi(\frac{a}{\sqrt{2}})$ of a DBPC $\mathbf{P}^{\left(\frac{a}{\sqrt{2}}\right)}$ until $\overline{\mathbf{Z}}^{(N)}_\ell$ remains constant or reaches the level $\ell_N$ for some sequence $\ell_N$ with $\ell_N \in \Theta( N^{\varepsilon})$ and $\varepsilon >0$ small enough, see Proposition \ref{Coupling Lower DBPC} and Proposition \ref{Asymptotic survival probability lower DBPC}. As for the lower bound the probability to reach the level $\ell_N$ is asymptotically equal to $\pi(\frac{a}{\sqrt{2}})$. 
When the level $\ell_N$ is reached we show that the total number of empty vertices grows in a finite number of generations to a level $N^{\frac{\beta}{2}+\delta}$ for some small $\delta>0$ whp, see Lemma \ref{Lem: beta over 2 + delta}.
Afterwards the remaining hosts get infected whp in a single generation. Indeed, the probability that a particular vertex gets attacked by at most one parasite can be (roughly) upper bounded by \[\left(1 - \frac{1}{D_N}\right)^{\frac{ (V_N N^{\beta/2 +\delta} )^2}{2}} +  \frac{(V_N N^{\beta/2 +\delta} )^2}{2} \left(1 -  \frac{1}{D_N}\right)^{\frac{(V_N N^{\beta/2 +\delta} )^2}{2}-1 } \frac{1}{D_N} \approx \exp\left(- N^{2\delta}/2 \right),  \] since roughly $\left(V_N N^{\beta/2+ \delta}\right)^2/2$ many pairs of parasites can be formed. Hence, the probability that at least one vertex is attacked by at most one parasites can be upper bounded by (roughly)
\[N\exp\left(- N^{2\delta}/2 \right)\rightarrow 0,\]
see Lemma
\ref{Lem: every vertices are killed} for details.\\
Case (i): We show that with asymptotically probability 1 the parasite population does not survive the first generation. \\
Case (iii): 
We show that we can whp
couple from below $\overline{\mathcal{I}}$ with the total size process of a Galton-Watson process with approximately Poi$(a^2/2)$ offspring distribution until $N^\alpha$ hosts get infected or the parasite population dies out for any $0<\alpha<\beta$ and any $a>0$. By choosing $\alpha>\beta/2$ we can show that once the level $N^\alpha$ is reached whp after one more generation the remaining hosts get infected. Since the probability to reach the level $N^\alpha$  is asymptotically equal to the survival probability $\varphi_a$ of a Galton-Watson process with Poi$(a^2/2)$ offspring distribution and $\varphi_a$  tends to 1 for $a\rightarrow \infty$  the result follows.\\
\\
We proceed with a \textit{sketch of the proof of Theorem \ref{MainResult}}:

Claim 1) (ii) and Claim 2):
For our upper bound on the survival probability we couple (as in the case of the complete graph) $\overline{\mathcal{I}}$ with a DBPC with offspring and cooperation distributions that are approximately Poisson distributions until a certain number $\ell_N$ of hosts get infected or the parasite population dies out, for a sequence $\ell_N \rightarrow \infty$ sufficiently slowly, see Proposition \ref{UpperBoundRGG}. The parameter of the approximating Poisson distribution for the offspring distribution is roughly $\frac{a^2}{2}$, since if all vertices are occupied with hosts the number of CoSame infections is on average approximately $\binom{v_N}{2} \frac{1}{d_N}$. The Poisson parameter of the cooperation distribution is roughly $a^2$, since  cooperation is maximal, if two balls centered around vertices, on which parasites have been generated in the same generation, are completely overlapping. In case of a complete overlap the number of cooperation events is on average roughly $v_N^2 \frac{1}{d_N}\sim a^2$.
Then we show that the probability to reach with the upper DBPC the level $\ell_N$ is asymptotically equal to the survival probability of the DBPC. This yields the upper bound, since again we upper bound the probability to infect the remaining hosts afterwards by 1.

For the lower bound we consider the spread of the parasites restricted to a certain \textit{complete neighbourhood} $\mathcal{C}(v_0)$ of the vertex $x_0$, that gets initially infected. The set $\mathcal{C}(x_0)$ contains all Poisson points with a distance $r_N/2$ to $x_0$. Since any two points in $\mathcal{C}(x_0)$ have a distance of at most $r_N$ any two points are connected over an edge, in other words the restriction  of $G_N$ to points in $\mathcal{C}(x_0)$ is a complete graph. Consequently, also the infection process restricted to $\mathcal{C}(x_0)$ is an infection process on a complete graph. In particular, we can apply Theorem \ref{Theorem:Invasion on complete graph} to show that the probability to infect at least $N^\eps$ vertices can be asymptotically lower bounded by the survival probability $\pi(\frac{a}{\sqrt{2^n}})$ of a DBPC $\mathbf{P}^{\left(\frac{a}{\sqrt{2^n}}\right)}$. The parameter of the offspring distribution is roughly  $\frac{a^2}{2^{n+1}}$, 
since $\binom{v_N}{2}$  pairs of parasites
can be generated per infected hosts and the probability that for pair of parasites both parasites hit the same vertex and that the vertex lies in $\mathcal{C}(x_0)$ is roughly
$\tfrac{\left(2 \tfrac{r_N}{2} \right)^n N}{\left(2r_N \right)^n N} \frac{1}{d_N}$, so the number of CoSame infections per host is roughly $\binom{v_N}{2} \tfrac{\left(2 \tfrac{r_N}{2} \right)^n N}{\left(2r_N \right)^n N} \frac{1}{d_N} \rightarrow \tfrac{a^2}{2^{n+1}}. $ 
 Similarly the parameter for the cooperation distribution is $\left(v_N\right)^2 \frac{1}{2^n} \frac{1}{d_N}\sim \frac{a^2}{2^n}$.
Once $N^\varepsilon$ many hosts are infected we show that after at most $O(\log (\log(N))$ many further generations the infection process expands from $x_0$ by a distance $r_N(1- o\left(1\right))$ per generation, see also Figure \ref{fig:rgg2d-boxcoverage}. Indeed, we show in Section \ref{section; travelling wave} that once the complete $r_N$-neighbourhood of a vertex gets infected, we can move the front from this vertex on in each generation by a distance $r_N\left(1- \frac{N^{-\frac{1-\beta/2}{n}+\delta}}{r_N}\right)$ for some $\delta>0$ small enough (which is the scale of the maximal distance possible). Consequently, after roughly at most $\frac{1}{2 r_N\left(1- \frac{N^{-\frac{1-\beta/2}{n}+\delta}}{r_N}\right)}$ many generations the complete cube is infected. On the other hand the invasion time is lower bounded by $\frac{1}{2 r_N}$, since parasites can move in any generation at most at a distance $r_N$ and the infection starts in the center of the cube. This explains our Claim 2) on the invasion time.
\\Case (i): As in the case of the complete graph we show that with asymptotically probability 1 the parasite population does not survive the first generation. \\
Case (iii): 
Again as in the case of the complete graph we show that we can whp
couple from below $\overline{\mathcal{I}}$ with the total size process of a Galton-Watson process with approximately Pois$(a^2/2)$ offspring distribution until $N^{\beta'}$ hosts get infected or the parasite population dies out for any $0<\beta'<\beta$ and any $a>0$. In addition we can show that when the level $N^{\beta'}$ is reached there exists a ball of radius $r_N$ which contains at least $N^{\beta'}/\log(N)$ infected hosts. By choosing $\beta'>\beta/2$ we can show that once the level $N^{\beta'}$ is reached whp after one more generation the remaining hosts in this ball get infected. Afterwards the infection expands by a distance $r_N(1-o(1))$ in every generation whp  (similar as in Case (ii)) until the remaining hosts are all infected.
Since the probability to reach the level $N^{\beta'}$  is asymptotically equal to the survival probability $\varphi_a$ of a Galton-Watson process with Pois$(a^2/2)$ offspring distribution and $\varphi_a$  tends to 1 for $a\rightarrow \infty$  the result follows.\\

\subsection{Simulating spatial invasion of cooperative parasites}\label{Sec: Simulations}

We supplement our findings with simulation results for moderately sized, finite $N$. We simulated invasion of parasites in host populations structured according to random geometric graphs on 
 (i) the interval $[0, 1]$ with the euclidean metric (which agrees with the maximum metric, since $n=1$), (ii) the square $[0, 1]^2$ using the maximum metric, (iii) the unit $2$-sphere $S^2$ using spherical distances (to substantiate our conjecture given in Remark \ref{remarkManifolds} at least by means of simulations).

To ease computations in the case of the $2$-sphere, we generate points on the unit $2$-sphere $S^2$, instead of the sphere with radius $\frac{1}{\sqrt{4\pi}}$ which would has as required in Remark \ref{remarkManifolds} a surface area (aka volume) of 1.
This simplification benefits both generation and evaluation of point distances in our implementation of the process and only requires appropriate rescaling.
The distance between two points $x$ and $y$ is then simply given by $\arccos(x \cdot y)$ as the radius is of length $1$.
Uniform points on $S^2$ can be generated by a two-step scheme in which first the polar angles $(\theta, \phi)$ are sampled using inverse transform sampling.
To this end, let $U_1, U_2$ be independent random variables with $U_1, U_2 \sim U(0, 1)$.
We compute $\theta = 2\pi \cdot U_1$ and $\phi = \arccos(1 - 2\cdot U_2)$ and obtain Cartesian coordinates by a standard transformation.

In general, storing and operating on an explicit representation of $\mathcal{G}$ takes space in the order of~$|\mathcal{E}^{(N)}|$ rendering parameter combinations of $N$ and $d_N = 2N^\beta$ infeasible for general-purpose compute architectures in case $Nd_N \sim 2N^{1+\beta}$ gets prohibitively large. 
Optimizations, however, are possible by implicit representations of $\mathcal{G}$ using the coordinates of~$\mathcal{V}^{(N)}$. 
Realizations of this are straight-forward for the interval $[0, 1]$ and can be adapted using Quadtrees for 2-dimensional spaces~\cite{S84}.\\
\\
\noindent
\textbf{Invasion probabilities}\\
\\
In Theorem~\autoref{MainResult} we claim that for $0 < a < \infty$ and $v_N \sim a \sqrt{d_N}$
\[\pi\left(\tfrac{a}{\sqrt{2}}\right)\leq \liminf_{N\rightarrow \infty} \PP\left(E_u^{(N)}\right) \leq  \limsup_{N\rightarrow \infty} \PP\left(E_u^{(N)}\right)\leq \pi(a). 
\]

\begin{figure}[t]
    \centering
    \includegraphics[scale=0.7]{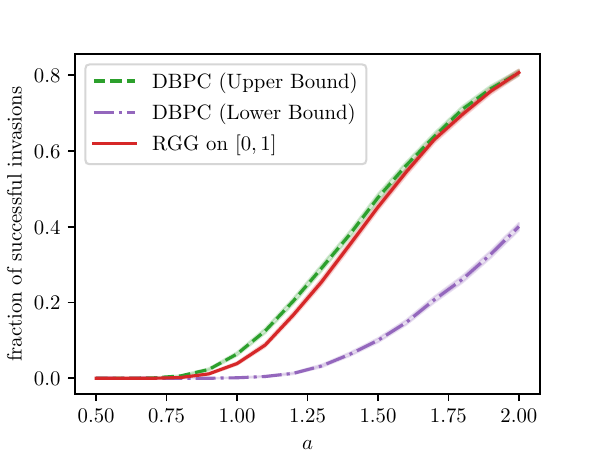}
    \hfill
    \includegraphics[scale=0.7]{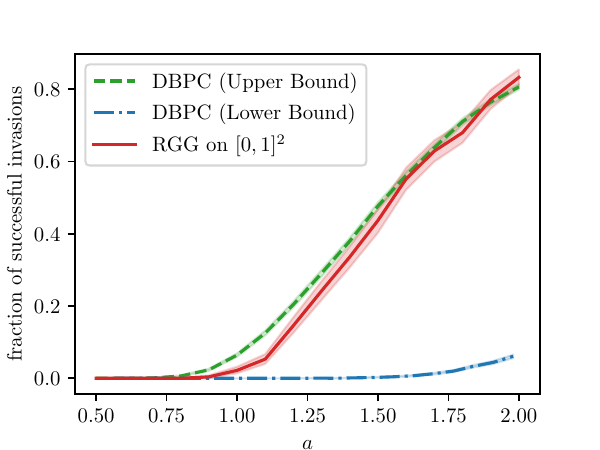}
    \caption{Simulated invasion probabilities with a host population structured by a random geometric graph (RGG, for short) on $[0, 1]$ (left) and $[0,1]^2$ (right) for $N = 10^6$ and $\beta = 0.7$ as well as simulated survival probabilities $\pi\left(\frac{a}{\sqrt{2^n}}\right)$ and $\pi\left(a\right)$.}
    \label{fig:2dim-vs-dbpc}
\end{figure}

\begin{figure}[t]
    \centering
    \includegraphics[scale=0.7]{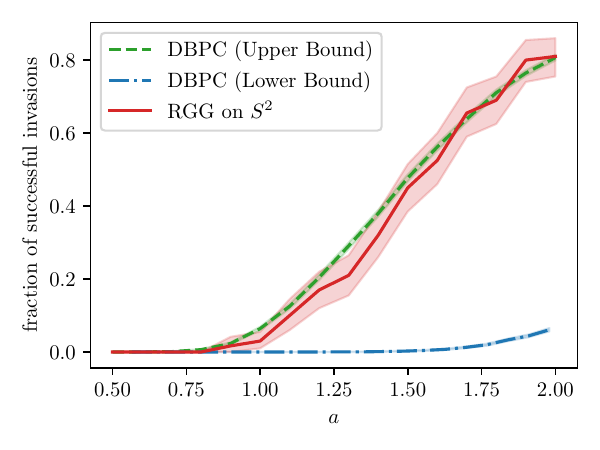}
    \hfill
    \includegraphics[scale=0.7]{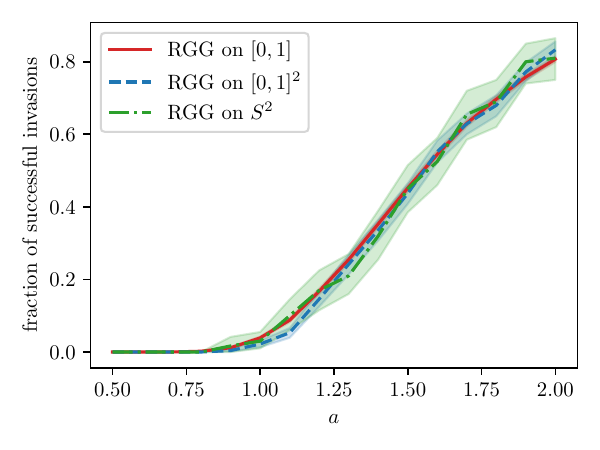}
    \caption{Fraction of successful invasions for $N = 10^6$ and $\beta = 0.7$ of infection processes spreading on host populations structured by a RGG on $S^2$ (left) and for comparison on $[0,1]$, $[0,1]^2$ and $S^2$ (right).}
    \label{fig:invasion-probability-all}
\end{figure}

In \autoref{fig:2dim-vs-dbpc} and \autoref{fig:invasion-probability-all} simulation results are depicted that show the fraction of cases, in which the host population got completely infected, for parasites spreading in host populations structured according to random geometric graphs on the interval $[0,1]$, the square $[0, 1]^2$ and the sphere $S^2$. In the simulations we assume that survival took place, if the DBPC attains size $N$. The simulated survival probabilities $\pi\left(\frac{a}{\sqrt{2^n}} \right)$ and $\pi\left(a \right)$ are based on 10000, 1000 and 200, resp. simulations for each value of $a$ for a RGG on $[0,1]$, $[0,1]^2$ and $S^2$, resp. They appear to be appropriate upper and lower bounds of the simulated invasion probabilities. The upper bound gives a particularly good approximation to the invasion probability. For the upper bound one assumes that the chance for two parasites, which have been generated on different vertices, to cooperate is roughly $\frac{1}{d_N}$, which is actually only true if the distance of the two vertices is 0. Therefore it might be surprising that the upper bound gives such a good fit. However, since parasites perform symmetric random walks a large part of parasites stays in a neighbourhood of $x_0$ and parasites that are close together have due to CoDiff infections a higher chance to generate offspring, which implies that parasites located in densely populated regions have in general more offspring parasites than parasites located in sparsely populated regions. 
 This effect remains until a significant proportion of the host population in a $r_N$-ball gets infected, but at this time point invasion is essentially already decided.
Consequently, the probability that a typical pair of parasites produces   CoDiff infections could be in the initial phase  pretty close to $\frac{1}{d_N}$.

Our asymptotic upper bound of the invasion probability does only depend on the ratio of the number $v_N$ of parasites generated on a vertex and the (asymptotic) number of direct neighbours of a typical vertex, but neither depends on the dimension nor (in the setting considered in Remark \ref{remarkManifolds}) on the curvature of the manifold. We suppose that this is also the case for the invasion probabilities. 
In \autoref{fig:invasion-probability-all} we present a direct comparison of simulated invasion probabilities for infection process on $[0,1]$, $[0, 1]^2$ and $S^2$ and see that the probabilities are very close to each other (even for finite $N$).

\begin{figure}[t]
    \centering
    \includegraphics[scale=0.7]{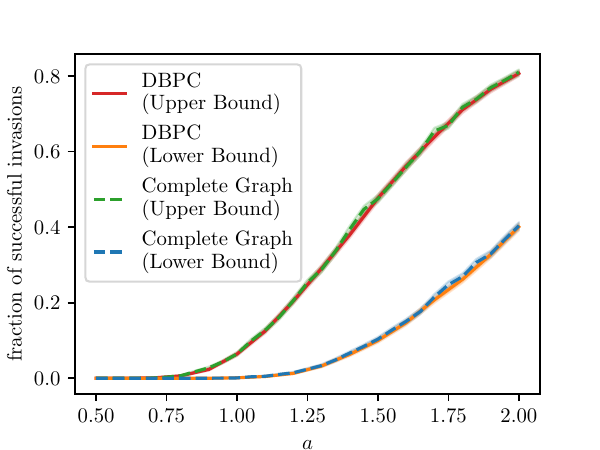}
    \hfill
    \caption{Simulated invasion probabilities with a host population structured by a complete graph for $N = 10^6$ and $\beta = 0.7$ as well as simulated survival probabilities $\pi\left(\frac{a}{\sqrt{2}}\right)$ and $\pi\left(a\right)$.}
    \label{fig:comparison-dbpc-cg}
\end{figure}

Finally we simulated invasion probabilities of the infection processes on the complete graphs that we use for a coupling from below. In \autoref{fig:comparison-dbpc-cg} one can observe that the simulated invasion  probabilities match very well with the probabilities $\pi\left(\frac{a}{\sqrt{2}}\right)$ and $\pi\left(a\right)$ of the corresponding DBPCs.\\
\\

\begin{figure}[t]
    \centering
    \includegraphics[scale=0.7]{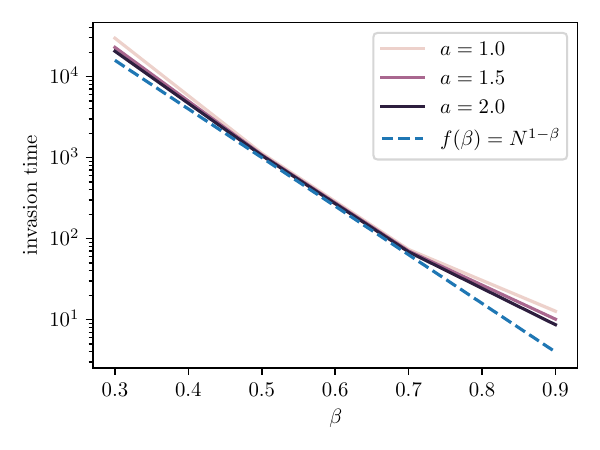}
    \hfill
    \includegraphics[scale=0.7]{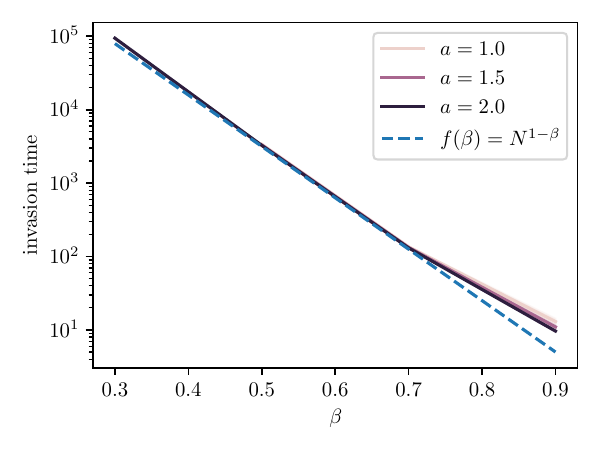}
    \caption{Invasion time on random geometric graphs on $[0, 1]$ for $N = 10^6$ (left) and $N = 10^7$ (right) with varying $a$ and $\beta$.}
    \label{fig:invasion-time-rgg1d}
\end{figure}

\begin{figure}[t]
    \centering
    \includegraphics[scale=0.7]{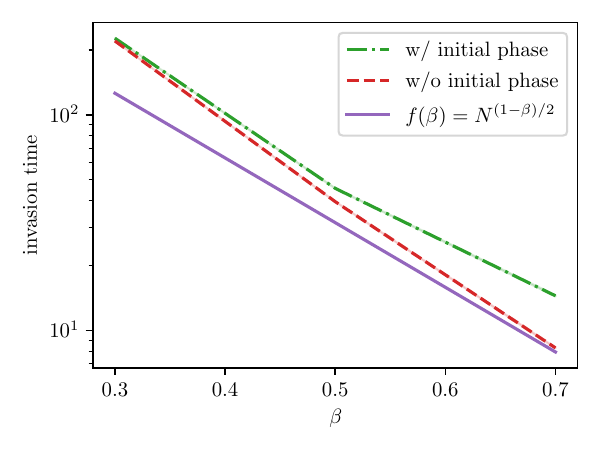}
    \hfill
    \includegraphics[scale=0.7]{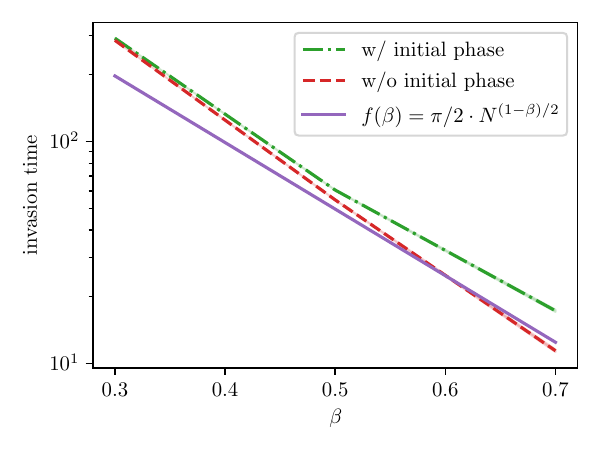}
    \caption{
        Invasion time on random geometric graphs on $[0, 1]^2$ (left) and $S^2$ (right) for $N = 10^6$ with $a = 2.0$.
    }
    \label{fig:invasion-time-rgg2d-s2}
\end{figure}

\noindent
\textbf{Invasion time} \\
\\
In \autoref{fig:invasion-time-rgg1d} and \autoref{fig:invasion-time-rgg2d-s2} we present the invasion time of simulated infection processes on the interval $[0,1]$, square $[0,1]^2$ and sphere $S^2$, respectively.
For reference we plot also the asymptotic order of the  invasion times derived in Theorem \ref{MainResult} and Remark \ref{remarkManifolds}.
In \autoref{fig:invasion-time-rgg1d} we observe a matching overlap that improves for increasing $N$ for all considered values of $a$ in the 1-dimensional case.

For large $\beta$ values the simulations showcase a higher invasion time than predicted.
This can be explained as follows: We show in Theorem \ref{MainResult} that the invasion time is asymptotically proportional to $N^{1-\beta}$. In particular,  the larger $\beta$ is the shorter is the invasion time. For $N\rightarrow \infty$ invasion is dominated by the time necessary to reach from a infinitesimally small neighbourhood of $x_0$ points close to the boundary of $[0,1]^n$ or in the setting of Remark \ref{remarkManifolds} the point that has the largest distance to the host that got initially infected. The initial phase until for the first time all direct neighbours of a vertex get infected is only of order $\log \log(N)$. For $\beta$ close to 1 and finite $N$ however both time frames are of approximately the same length, which explains the deviation from the theoretical prediction where the initial phase is ignored. 
In \autoref{fig:invasion-time-rgg2d-s2} we plotted the invasion time when the initial phase is removed. One observes that for intermediate and larger values of $\beta$ the gap between the predicted and simulated invasion time disappears. For larger values of $\beta$ the simulated invasion times lie slightly below the predicted invasion times. Probably this is caused by parasites spreading the infection further before the initial phase is over.

Also for small $\beta$ values we observe that simulated invasion times are generally higher than the predicted ones, even when the initial phases are removed.
This deviation is particularly pronounced for invasion on $[0,1]^2$, where the maximum metric is used. This can be explained as follows. As we pointed out in the sketch of the proof of Theorem \ref{MainResult} the parasite population expands furthest  due to parasites born at the boundary of an $r_N$ neighbourhood. When on the square the maximum distance is used $r_N$-squares on the diagonal can get infected fastest by parasites at the corners of neighbouring $r_N$-squares. However, when  $N$ is not large, the number of parasites located in the corners is pretty small, so that they might be not able to move the front forward as quickly as predicted for $N\rightarrow \infty.$

\begin{figure}[t]
    \centering
    \includegraphics[scale=0.7]{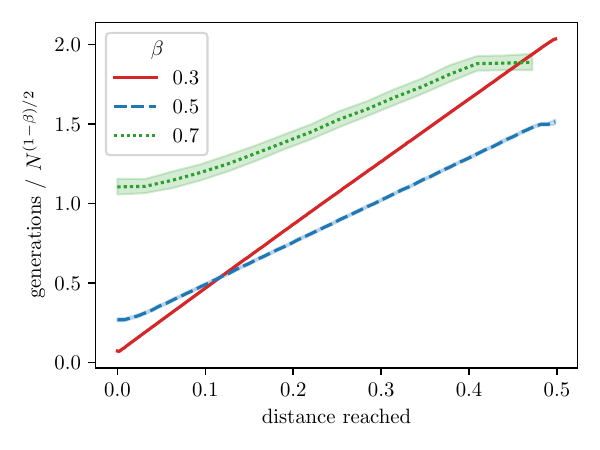}
    \caption{Distance reached for successful invasions on $[0,1]^2$ with $N=10^6$.}
    \label{fig:rgg2d-boxcoverage}
\end{figure}

This behavior is further studied in \autoref{fig:rgg2d-boxcoverage} where the progress of the infection process is tracked along boxes of radius $r_N/2$ on the unit-square $[0,1]^2$ for different values of $\beta$. The larger $\beta$ the more vertices are located at the corners. For $\beta=0.7$ and $\beta=0.5$ one observes that after the initial phase the parasite population expands linearly (almost) by a factor 1 (as predicted), while for $\beta=0.3$ (when in each box with edge length $r_N$ only $\approx N^{0.3}\approx 63$ vertices are contained) the population expands also linearly, but only by a  factor of (almost) 2.

\bigskip
In the following the manuscript is structured as follows. In Section \ref{lab: DBPC section} we show several properties (of sequences) of DBPCs that we will need in the subsequent section. Afterwards in Section \ref{Section: Infection on complete graph} we will prove Theorem \ref{Theorem:Invasion on complete graph}. Finally in the last section we will prepare and give the proof of Theorem \ref{MainResult}.

\section{Discrete branching processes with cooperation}
\label{lab: DBPC section}

In this section we collect properties of (sequences of) DBPC, that we need in the following. Some of the statements are  well-known or have been proven in \cite{BrouardEtAl2022} for Galton-Watson processes. As for these statements the proof techniques are similar, we do not give the proofs in the main text, but provide them in the supplementary material.

We start with the extinction-explosion principle, which is well-known for branching processes and also holds for DBPCs. 

\begin{Lem}(Extinction-explosion principle for DBPCs)
\label{Lemma: extinction-explosion principle DBPC}
Let $\ZZ$ be a DBPC satisfying $p_{1,o}\neq 1$ and $(p_{0,o},p_{1,c})\neq (1,1)$. Then there exists a nullset $\mathcal{N}$ such that
\begin{equation}
    \{\forall g \in \mathbb{N}_{0}: Z_{g}>0 \} \subseteq \{ \forall i \in \mathbb{N}, \exists g_0 \in \mathbb{N}_0, \forall g\geq g_0: Z_g \geq i\} \cup \mathcal{N}. 
\end{equation}
\end{Lem}

For the proof one shows that all states but 0 are transient states. Details of the proof can be found in the supplementary material.

 Before we proceed we introduce some useful notation. We denote the expectation of and the variance w.r.t.\ the offspring distribution by $\mu_o$ and $\nu_o$ and for the cooperation distributions by $\mu_c$ and $\nu_c$. 

In contrast to Galton-Watson process DBPC  we aim to show that except for pathological cases a DBPC always has a positive survival probability.
\begin{Lem}\label{lem:SurvivalDBPC}
     Let $\ZZ$ be a DBPC with $\mu_o,\mu_c,\nu_o,\nu_c>0$. Suppose $Z_0 = k>0$, $k\in \mathbbm{N}$, then $\ZZ$ has a positive survival probability, i.e.\
     \begin{equation*}
        \mP(  Z_g > 0 \,\forall\, g>0)>0.
    \end{equation*}
\end{Lem}

The proof of Lemma \ref{lem:SurvivalDBPC} is based on the next lemma, which we immediately formulate in a more general setting to be able to apply it also later in another context and which basically states that if a DBPC attains a certain level, then in subsequent generations the size will up to a constant factor (that does not depend on the size) be squared in the next generations  due to cooperation with a certain non-vanishing 
probability.

In the subsequent sections we will deal often with sequences of DBPCs rather, than a single process.
 We will often need the following assumption to be fulfilled.
\begin{Ass}\label{Ass:SequenceDBPC}
    Let $\big(\ZZ^{(N)}\big)_{N\in \N_0}$ be a sequence of DBPCs for which $\mu_o^{(N)},\mu_c^{(N)},\nu_o^{(N)},\nu_c^{(N)}$ denote the expectations and the variances of the offspring and cooperation distributions.
    We assume that 
    \begin{equation*}
        \mu_o^{(N)},\mu_c^{(N)},\nu_o^{(N)},\nu_c^{(N)}\underset{N \to \infty}{\longrightarrow} \mu_o,\mu_c,\nu_o,\nu_c>0.
    \end{equation*}
\end{Ass}

\begin{Lem}\label{lem:technicalBound}
    Let $\big(\ZZ^{(N)}\big)_{N\in \N}$ be a sequence of DBPCs which satisfies Assumption~\ref{Ass:SequenceDBPC}. Furthermore, set $f_i(k):= \frac{k^{2^i}\mu^{2^{i}-1}}{8^{2^{i}-1}}$ for $i\geq 1$, then there exists $N_0\in \N$ such that for any $N \geq N_0$, $g \in \mathbb{N}$, $M\in \N$ and 
    $k\geq L:= \lceil \mu^{-1}(8+\nu)^2 \rceil$ it holds that
    \begin{align*}
        \mP\left( \bigcap_{i=1}^{ M}\left\{ Z_{g+i}^{(N)} > f_i(k) \right\} \ \Big|\  Z_{g}^{(N)} =k\right) \geq \prod_{i=1}^M \left(1 -\frac{6\nu}{f_i(k) \mu} \right) \geq  \prod_{i=1}^{M} \bigg(1-\frac{3}{4(8+\nu)^{2^{i}-1} }\bigg).
    \end{align*}
\end{Lem}

The claim can be shown by several applications of Tchebychev's inequality, details can be found in the supplementary material.

Now we can proof Lemma \ref{lem:SurvivalDBPC}.
\begin{proof}[Proof of Lemma \ref{lem:SurvivalDBPC}]
    First of all $\mu_o,\mu_c>0$ implies that $p_{0,o},p_{0,c}<1$ and $\nu_o,\nu_c>0$ that $p_{1,o},p_{1,c}\neq 1$. This fact and since 
    $L:=\lceil \mu^{-1}(8+\nu))^2 \rceil$ is finite allow us to find a $p_1>0$ such that for any $x\in \{1, ..., L-1\}$ it exists an $g \in \mathbb{N}$ such that 
    \begin{align*}
        \mP_x( Z_g \geq L) \geq p_1.
    \end{align*}
  
    Applying Lemma~\ref{lem:technicalBound} for $Z^{(N)}\equiv Z$ one obtains 
    \begin{align*}
        \mP\left( \bigcap_{i=1}^{M} \{ Z_{g+i} > f_{i}(k) \} \Big| Z_{g} =k\right)  
         &\geq \prod_{i=1}^{M} \bigg(1-\frac{3}{4(8+\nu)^{2^{i}-1} }\bigg)
    \end{align*}
    for $k\geq L$. By continuity of $\mP$ we have
    \begin{align*}
        \mP\left( \bigcap_{i=1}^{\infty}\{ Z_{g+i} > f_{i}(L) \} \Big| Z_{g} =L\right) \geq
          \prod_{i=1}^{\infty}  \bigg(1-\frac{3}{4(8+\nu)^{2^{i}-1} }\bigg)>0,
    \end{align*}
    where it follows that the right-hand side is strictly positive by comparison with a geometric sum. We have that 
    \begin{equation*}
        \bigcap_{i=1}^{\infty}\{ Z_{g+i} > f_{i}(L) \} \cap \{Z_{g} \geq L\}  \subset \{Z_{g} >0 \,\,\forall g\geq 0\},
    \end{equation*}
    then by Markov property and monotonicity we get
    \begin{align*}
        \mP(Z_{g} >0 \,\forall g\geq 0)\geq \mP(Z_{g+i} >0 \,\forall i\geq 0|Z_{g} = L)  \mP_x( Z_{g} \geq L)>0. 
    \end{align*}
\end{proof}

The next lemma claims that  reaching a level $b_N$ that tends to $\infty$ (arbitrarily slowly) with $N$  implies that the DBPC survives whp.
\begin{Lem}\label{lem:UniformlyConvergence}
    Let $\big(\ZZ^{(N)}\big)_{N\in \N_0}$ be a sequence of DBPCs which satisfy Assumption~\ref{Ass:SequenceDBPC} and let $(b_N)_{N\in \mathbb{N}}$ be a sequence of positive numbers with $b_N\to \infty$ as $N\to \infty$. Then for all $g \in \mathbb{N}$ it holds that
    \begin{equation*}
        \mP\left( Z_{g+i}^{(N)} >0 \ \forall i\geq 0 \ \big|\  Z_{g}^{(N)} \geq b_N\right) \underset{N \to \infty}{\longrightarrow} 1.
    \end{equation*}
\end{Lem}

The proof relies on Lemma \ref{lem:technicalBound}, we provide a detailed proof in the supplementary material.

Finally we are able to derive results on the expansion speed in case of survival. The next lemma shows that for any sequence $K_N \to \infty$ reaching the level $K_N$ or dying out is at most of order $\log(\log(K_N))g(N)$ with $g(N)\to \infty$ arbitrarily slowly. 

\begin{Lem}\label{GrowFast0}
  Let $\big(\ZZ^{(N)}\big)_{N\in \N_0}$ be a sequence of DBPCs which satisfy Assumption~\ref{Ass:SequenceDBPC}. Let $\left(K_N\right)_{N \in \mathbb{N}}$ be a sequence satisfying $K_N \rightarrow \infty$ as $N\rightarrow \infty$. Assume $Z_0^{(N)} = x$ for some $x\in \{1, ..., K_N-1\}.$ Furthermore, assume that there exists an $N_0>0$ such that 
  \begin{equation} \label{Equation: assumption reaching L}
      \inf_{N\geq N_0} \mP \left(Z_1^{(N)}>L|Z^{(N)}_{0}=x\right)>0,
  \end{equation}
  for all $N>N_0$, where $L=:\lceil \mu^{-1}(8+\nu)^2 \rceil$.
  Denote by 
 \begin{align*}
  \tau_{K_N, 0} = \inf\{ g \in \mathbb{N}: Z_g^{(N)} =0 \text{ or }  Z_g^{(N)} \geq K_N\}.
 \end{align*}
Then there exist constants $q\in(0,1)$, $C>0$, and $N_{1} \in \mathbb{N}$ such that for all $N \geq N_1$ and for all $g(N) \to \infty$
 \begin{align*}
  \PP\left(\tau_{K_N,0}\leq Cg(N)\log(\log(K_N))\right) \geq 1-q^{\lfloor g(N)\rfloor}.
 \end{align*}
\end{Lem}

\begin{proof}
 We show that there exists $p_0>0, C, N_1 \in \N$ (independent of $N$), such that for all $N\geq N_1$ and all $x \in \{1, ..., K_N-1\}$
 \begin{align}\label{Geom}
  \mP_x\left(\tau_{K_N,0} \leq C \log (\log( K_N )) \right) \geq p_0.
 \end{align}
 From this follows, that for all $k\in \N$
 \begin{align*}
 \mP_x\left(\tau_{K_N, 0} \leq k C \log (\log( K_N ))  \right) \geq \sum_{i=1}^k (1-p_0)^{i-1} p_0.
 \end{align*}
 In particular for $k=g(N)$
  \begin{align*}
 \mP_x\left(\tau_{K_N, 0} \leq C g(N) \log(\log(K_N))\right) \geq \sum_{i=1}^{\lfloor g(N)\rfloor} (1-p_0)^{i-1} p_0= 1 - (1- p_0)^{\lfloor g(N)\rfloor}.
 \end{align*}
 It remains to show \eqref{Geom}. Since $L$ is finite and does not depend on $N$ we find a $p_1$ such that for any $x\in \{1, ..., L-1\}$ 
 \begin{equation*}
  \mP \left(Z_g^{(N)} =0 \text{ or } Z_g^{(N)} > L| Z_{g-1}^{(N)} =x\right) \geq p_1,
 \end{equation*}
 due to Equation \eqref{Equation: assumption reaching L}.\\
 Consequently, we can lower bound the time to reach the state 0 or a state $>L$ by a geometrically distributed random variable with success probability $p_1$ for any $x\in \{1, ..., L-1\}$. Reasoning as above shows that the waiting time to hit 0 or a state $>L$ is with probability $1- o(1)$ bounded by $\log \log (K_N)$.

We show next, that if $Z_g^{(N)}\geq L$ after $C_1 \log \log (K_N)$ further generations the level $K_N$ will be reached with some probability $p_2>0$ for any $N$ large enough. We use Lemma~\ref{lem:technicalBound} for this. Recall that $f_i(k)= \frac{k^{2^i}\mu^{2^{i}-1}}{8^{2^{i}-1}}$. Lemma~\ref{lem:technicalBound} implies that for $k\geq L$ and $M= C \log \log (K_N)$ it follows that 
 \begin{align}\label{lab:KeyInqualitySurv}
 \mP\left( \bigcap_{i=1}^{ C \log \log (K_N)}\left\{ Z_{g+i}^{(N)} > f_i(k) \right\} \Big| Z_{g}^{(N)} =k\right) 
  &\geq \prod_{i=1}^{C \log \log (K_N)} \bigg(1-\frac{3}{4(8+\nu)^{2^{i}-1} }\bigg)>p_2 >0.
 \end{align}
Now by \eqref{Equation: Lower bound in doubling steps DBPC} we know that 
$f_{i}(k)>\mu^{-1}8(8+\nu)^{2^{i}}$. 
Thus, by choosing $C=\frac{1}{\log(2)}$ we have that
\begin{equation*}
    f_{C \log \log (K_N)}(k)
    \geq \frac{8}{\mu} (K_N)^{\log(8)}.
\end{equation*}
Since $\log(8)>1$ we get that $f_{C \log \log (K_N)}(k) \geq K_N$ for $N$ large enough, which yields the claim. 
 \end{proof}

The next proposition improves the statement of the last lemma. It claims, that in at most an order of $\log(\log(K_N))$ generations whp the level $K_N$ is reached or 0 is hit.
 \begin{Prop}\label{GrowFast}
  Let $(\ZZ^{(N)})_{N\in \N_0}$ be a sequence of DBPCs which satisfy Assumption~\ref{Ass:SequenceDBPC}. Let $\left(K_N\right)_{N \in \mathbb{N}}$ be a sequence satisfying $K_N \rightarrow \infty$ as $N\rightarrow \infty$. Assume $Z_0^{(N)} = x$ for some $x\in \{1, ..., K_N-1\}.$ Furthermore, assume that there exists a $N_0>0$ such that 
  \begin{equation} 
      \inf_{N\geq N_0} \mP \left(Z_1^{(N)}>L|Z^{(N)}_{0}=x\right)>0,
  \end{equation}
  for all $N>N_0$, where $L=:\lceil \mu^{-1}(8+\nu)^2 \rceil$.
  Denote by 
 \begin{align*}
  \tau_{K_N, 0} = \inf\{ g \in \mathbb{N}: Z_g^{(N)} =0 \text{ or }  Z_g^{(N)} \geq K_N\}.
 \end{align*}
Then there exists a constant $C>0$ such that  
 \begin{align*}
  \PP\left(\tau_{K_N,0}\leq C\log(\log(K_N))\right) \to 1.
 \end{align*}
\end{Prop}

\begin{proof}
Let $b_N=\log(K_N)$ such that applying Lemma \ref{GrowFast0} to $b_N$ and $g(N)=\frac{\log(\log(K_N))}{C\log(\log(b_N))}$ we obtain that
\begin{equation}\label{lab: control reaching level to infty}
    \mathbb{P} \left( \tau_{b_N,0} \leq \log(\log(K_N))\right) \geq 1-q^{\lfloor g(N) \rfloor} \underset{N \to \infty}{\longrightarrow}1.
\end{equation}
For $N$ large enough such that $b_N \geq L$ we have according to \eqref{lab:KeyInqualitySurv0} that 
 \begin{align}
 \mP\left( \bigcap_{i=1}^{ C\log \log (K_N)}\left\{ Z_{g+i}^{(N)} > f_i(b_N) \right\} \Big| Z_{g}^{(N)} =b_N\right) 
  &\geq \prod_{i=1}^{ C\log \log (K_N)} \bigg(1-\frac{48 \nu}{f_{i-1}^{2}(b_N)\mu^{2}}\bigg).
 \end{align}
By definition we have for all $i \in \mathbb{N}$ that $f_i(b_N)= \frac{b_N^{2^i}\mu^{2^{i}-1}}{8^{2^{i}-1}}$. Taking $i=C\log(\log(N))$ with $C=\frac{1}{\log(2)}$ gives that
 \begin{align}
     f_{C\log(\log(K_N))}(a_N)=\frac{8}{\mu}\left(\frac{b_N \mu}{8}\right)^{2^{C\log(\log(K_N))}}&=\frac{8}{\mu}\left(\frac{b_N \mu }{8}\right)^{\log(K_N)}=\frac{8}{\mu} K_N^{\log\left(\frac{b_N \mu}{8}\right)},
 \end{align}
 and $\log\left(\frac{b_N\mu}{8}\right)>1$ for $N$ large enough, such that $f_{C\log(\log(K_N))} (b_N)\geq K_N$. 
 
And because $f_{0}(b_N)=b_N$ we also have that 
 \begin{align} \label{lab: proba to 1 loglog generation}
     \prod_{i=1}^{C\log(\log(K_N))} \left( 1-\frac{48 \nu}{f_{i-1}^{2}(b_N)\mu^{2}}\right) \geq \left(1-\frac{48 \nu}{b_N^2 \mu^{2}} \right)^{C\log(\log(K_N))}\sim\exp\left(-48\nu C\frac{\log(\log(K_N))}{\log^2(K_N)\mu^2} \right)\underset{N \to \infty}{\to}1.
 \end{align}

Using the strong Markov property at the stopping time $\tau_{b_N,0}$ gives that 
\begin{align}
    &\mathbb{P} \left( \tau_{K_N,0} \leq \left(1+\frac{1}{\log(2)}\right) \log(\log(K_N))\right)
    \geq \mathbb{P} \left(\tau_{b_N,0}\leq \log(\log(K_N)),Z_{\tau_{b_N,0}}^{(0)}=0\right) \\
    &\hspace{1.5cm}+\mathbb{P} \left( \tau_{b_N,0} \leq \log(\log(K_N)),Z_{\tau_{b_N,0}}^{(N)} \geq b_N\right) \mathbb{P} \left( \tau_{K_N,0} \leq \frac{1}{\log(2)}\log(\log(K_N))\vert Z_{0}^{(N)}=b_N\right) \\
    &\hspace{2cm}\underset{N\to \infty}{\longrightarrow} 1,
\end{align}
according to \eqref{lab: control reaching level to infty} and \eqref{lab: proba to 1 loglog generation}.

\end{proof}

 The following lemma states that the probability of reaching an arbitrary high level, that tends to $\infty$ as $N\rightarrow \infty$, at some generation or up to some
generation is asymptotically equal to the survival probability for a sequence of DBPCs.
\begin{Prop}
\label{Lemma: convergence survival probability DBPC}
 Consider a sequence of DBPC $\left(\ZZ^{(N)}\right)_{N \in \mathbb{N}}$ with offspring and cooperation distributions $\big(p_{k,o}^{(N)}\big)_{k \in \mathbb{N}_{0}}$ and $\big(p_{k,c}^{(N)}\big)_{k \in \mathbb{N}_{0}}$ respectively, which satisfies Assumption~\ref{Ass:SequenceDBPC}. Furthermore, let $\ZZ$ be a DBPC with offspring and cooperation distribution $\left(p_{k,o}\right)_{k \in \mathbb{N}_{0}}$ and $\left(p_{k,c}\right)_{k \in \mathbb{N}_{0}}$. Assume that $p_{k,o}^{(N)}\to p_{k,o}$ and $p_{k,c}^{(N)}\to p_{k,c}$ as $N\to \infty$ for all $k\geq 0$. 

Then for any $\mathbb{N}$-valued sequence $(a_N)_{N \in \mathbb{N}}$ with $a_N \to \infty$ it holds that 
 \begin{align}
 \label{Equation: total size DBPC converges surv prob}
     \lim_{N \to \infty} \mathbb{P} \left(\forall g \in \mathbb{N}_{0}: Z_{g}^{(N)}>0 \right)&=\lim_{N \to \infty} \mathbb{P} \left( \exists g \in \mathbb{N}_{0}: Z_{g}^{(N)} \geq a_N\right)\\
     &=
     \lim_{N \to \infty} \mathbb{P} \left( \exists g \in \mathbb{N}_{0}: \overline{Z}_{g}^{(N)} \geq a_N\right)\\
     &=\pi,
 \end{align}
 where $\pi$ denotes the survival probability of $\ZZ$.
\end{Prop}

We provide the proof in the supplementary material.

\section{Invasion of cooperative parasites in host populations structured by a complete graph}\label{Section: Infection on complete graph}

In this section we prepare and give the proof of Theorem \ref{Theorem:Invasion on complete graph}.

We will often use the inequalities 
$\exp(-x) \geq 1-x \geq \exp(-x) \exp(-x^2)$,
$\exp(-x) \leq  1-2x$ for $x\in [0,\frac{1}{2}]$ and Bernoulli's inequality $(1+x)^i \geq 1 +  i x$ for $i\in \mathbbm{N}$ and $x\geq -1$ in this and the next section.

Furthermore we will compare the infection dynamics happening within one generation often with balls-into-boxes experiments. The following lemma gives control about certain events happening in these experiments.

\begin{Lem} \label{lab: balls and boxes experiment}
Let $(m'_N), (V'_N), (h'_N), (\ell'_N), (D'_N)$ be non-negative sequences with $0\leq m'_N \leq h'_N\leq \ell'_N$ and assume $\frac{{\ell'_N}^4 {V'_N}^3}{{D'_N}^2} \in o\left(1\right)$.
Consider $D'_N-m'_N$ boxes and $m'_N V'_N$ balls. Assume the balls are put independently and purely at random into the boxes. 
Consider the event $C_k^{(h'_N)}$, that $k$ many of the first  $D'_N-h'_N$ boxes contain exactly two balls and the remaining boxes contain at most one ball.
We have for all $k \leq \ell'_N$
\begin{align*}
\mathbb{P}(W'=k) \exp\left( -\frac{{\ell'}_N^5 {V'}_N^3}{{D'}_N^2}\right) &\leq \mathbf{P}\left(C_k^{(h'_N)}\right) \\
&\leq \left( \frac{(m'_N V'_N)^2}{2 D'_N} \right)^k \frac{1}{k!} \exp\left(-\frac{(m'_N V'_N - 2 \ell'_N)^2}{2 D'_N} \right) \exp\left( \frac{{\ell'}_N^2 V'_N}{D'_N} \right),
\end{align*}
where $W'$ is Poisson distributed with parameter $\frac{\left(m'_N V'_N-2\ell'_N\right)^2}{2 D'_N}$.
\end{Lem}
A proof can be found in the supplementary material. \\
\\
In the following we will denote by $W^{(N)}_{o,m}$ a Poisson distributed random variable with parameter $\frac{m\left(V_N-2\ell_N\right)^2}{2 D_N}$ for any $m\in \mathbbm{N}$ and similarly by $W^{(N)}_{c,m}$ a Poisson distributed random variable with parameter $\frac{m\left( V_N-2\ell_N\right)^2}{D_N}$.

\subsection{Results to arrive at upper bounds for invasion probabilities}
To derive an upper bound on the invasion probability we estimate from above the total number of infected hosts by the total size of a branching process with cooperation with offspring and cooperation distributions that 
are approximately Poisson distributed. 
\begin{Def}(Upper DBPC) \\
\label{Definition: UBPI}
Let $\ell_{N} \underset{N \to \infty}{\to}\infty$ satisfying $\frac{\ell_N^7 V_{N}^{3}}{D_{N}^{2}}\in o(1)$. Let $\ZZ^{(N)}_{u}=\left(Z_{g,u}^{(N)}\right)_{g \in \mathbb{N}_0}$ be a branching process with cooperation with $Z_{0,u}^{(N)}=1$ almost surely, and offspring and cooperation distribution with probability weights $p_{u,o}^{(N)} = \left( p_{j,u, o}^{(N)}\right)_{j \in \mathbb{N}_{0}}$  and  $p_{u,c}=\left( p_{j,u,c}^{(N)}\right)_{j \in \mathbb{N}_{0}}$, respectively  with 
\begin{align}
    p_{j,u, o}^{(N)}:= \mathbb{P}(W^{(N)}_{o,1} =j)
 \exp \left( -\frac{\ell_N^{5}V_{N}^{3}}{D_{N}^{2}}\right),
\end{align}
for all $0 \leq j \leq \ell_N$ and 
\begin{align}
    p_{\ell_N+1,u,o}^{(N)}:=1-\sum_{j=0}^{\ell_N} p_{j,u,o}^{(N)},
\end{align}
as well as 
\begin{align}
    p_{j,u, c}^{(N)}:= \mathbb{P}(W^{(N)}_{c,1} =j) \exp \left( -2\frac{\ell_N^{5}V_{N}^{3}}{D_{N}^{2}}\right),
\end{align}
for all $0 \leq j \leq \ell_N$ and 
\begin{align}
    p_{\ell_N+1,u,c}^{(N)}:=1-\sum_{j=0}^{\ell_N} p_{j,u,c}^{(N)}.
\end{align}
Denote by $\overline{\ZZ}_{u}^{(N)}:= \big( \overline{Z}_{g,u}^{(N)}\big)_{g \in \mathbb{N}_{0}}$ where $\overline{Z}_{g,u}^{(N)}:=\sum_{i=0}^{g}Z_{i,u}^{(N)}$, that is $\overline{Z}_{g,u}^{(N)}$ gives the total size of $\ZZ^{(N)}_{u}$ accumulated till generation $g$.
\end{Def}

In the next proposition, we show that the total size of the infection process $\overline{\mathcal{J}}^{(N)}$ can be coupled whp from above with the total size of the DBPC of Definition \ref{Definition: UBPI} $\overline{\mathbf{Z}}^{(N)}_{u}$ up to the first random generation at which $\overline{\mathbf{Z}}^{(N)}_{u}$ reaches the size $\ell_N$ (for $\ell_N \to \infty$ not too fast) or the process $\mathbf{Z}^{(N)}_{u}$ dies out. 
\begin{Prop}\label{Prop: Coupling with upper DBPC}
Consider a sequence $\left(\ell_N \right)_{N \in \mathbb{N}}$ with $\ell_N \underset{N \to \infty}{\to}\infty$ satisfying $\ell_N^6\frac{V_N^3}{D_N^2}\in o(1)$. Introduce the stopping time
\begin{align}
    \tau^{(N)}_{\ell_N,0}:= \inf \left\{ g \in \mathbb{N}_{0}: \overline{Z}^{(N)}_{g,u} \geq \ell_{N} \text{ or } Z_{g,u}^{(N)}=0\right\}.
\end{align}
Then 
\begin{align}
    \lim_{N \to \infty} \mathbb{P} \left( \overline{J}^{(N)}_{g} \leq \overline{Z}^{(N)}_{g,u}, \forall n < \tau^{(N)}_{\ell_N, 0}\right)=1,
\end{align}
and 
\begin{align}
    \lim_{N \to \infty} \mathbb{P} \left( J^{(N)}_{\tau^{(N)}_{\ell_N, 0}} =0 \Big\vert  Z^{(N)}_{\tau^{(N)}_{\ell_N, 0}, u} =0\right)=1 \text{ and } \lim_{N \to \infty} \mP \left( \overline{J}^{(N)}_{\tau^{(N)}_{\ell_N,0}} \geq \ell_N \Big\vert \overline{J}^{(N)}_{\tau^{(N)}_{\ell_N,0}} \geq \overline{Z}^{(N)}_{\tau^{(N)}_{\ell_N,0}}\right)=1.
\end{align}
\end{Prop}

\begin{proof}
Up to generation $\sigma_{\ell_N,0}^{(N)}:=\inf \{g \in \mathbb{N}_{0}: \overline{I}_g^{(N)} \geq \ell_N \text{ or } I_{g}^{(N)}=0\}$ the total number of parasites that are moving in the graph is upper bounded by $\ell_N V_N$. Consider the following experiment with $D_N-\ell_N$ boxes, $\ell_N V_N$ balls. Assume that balls are thrown uniformly at random in the boxes. The probability that there exists a box with at least 3 balls in can be upper bounded as follows 
\begin{align}
    \mathbb{P}\left(\exists \text{ 1 box with more than 3 balls} \right) \leq D_N \left(\ell_NV_N^3\right)^3 \frac{1}{\left(D_N-\ell_N\right)^{3}} \sim \frac{\ell_N^3 V_N^3}{D_N^2} \to 0. 
\end{align}
This means that with the assumed scaling of $\ell_N$ it is unlikely that such an event occurs before generation $\tau_{\ell_N,0}^{(N)}$. Consequently for whp couplings, we can only focus on infections generated by pairs of parasites.\\
Now consider a complete graph with exactly $1\leq m_N <\ell_N$ infected vertices and at most $\ell_N-1 $ empty or infected vertices. The probability that $k$, $0\leq k \leq \ell_N$ infections are generated can be estimated from above by the probability that $k$ boxes are filled with at least two balls and the remaining boxes are filled by at most one ball in the following balls-into-boxes experiment: consider $D_N-m_N$ boxes and $m_N V_N$ balls. Place the balls uniformly at random into the boxes. Denote by $C_k^{(N)}$ the event that $k$ boxes contain exactly two balls and all other boxes contain at most one ball. \\

By Lemma \ref{lab: balls and boxes experiment} with $D'_N= D_N, m'_N= m_N, V'_N= V_N$ and $h'_N = m_N$ we can estimate for $m_N \geq 2$
\begin{align}\label{Eq for D_k}
    \mathbb{P} \left(C_k^{(N)}\right)
    \geq \left( \frac{m_N^2(V_{N}-2\ell_{N})^{2}}{2 D_N}\right)^{k}\frac{1}{k!} \exp \left(-m_N^{2}\frac{V_N^{2}}{2 D_N} \right) \exp \left( -\frac{\ell_N^{5}V_{N}^{3}}{D_{N}^{2}}\right) 
\end{align}
again for $N$ large enough, and for $m_N=1$
\begin{align}\label{Eq for D_k 2}
    \mathbb{P} \left(C_k^{(N)}\right) 
 \geq 
    \left( \frac{(V_{N}-2\ell_{N})^{2}}{2 D_N}\right)^{k}\frac{1}{k!} \exp \left(-\frac{(V_N-2\ell_N)^{2}}{2 D_N} \right) \exp \left( -\frac{\ell_N^{5}V_{N}^{3}}{D_{N}^{2}}\right) 
\end{align}    
In order to prove that $\overline{\mathcal{J}}^{(N)}$ can be coupled  with $\overline{\mathbf{Z}}^{(N)}_u$ such that $\overline{\mathbf{Z}}^{(N)}_u$ dominates $\overline{\mathcal{J}}^{(N)}$, we show that $\mathbb{P}  \big( Z_{n+1,u}^{(N)}  =k | Z_{n,u}^{(N)} =m_N\big) \leq$ \eqref{Eq for D_k}.

Consider independent random variables $\big(X_i^{(N)}\big)_{i \in \mathbb{N}}$ and $\big(Y_{(i,j)}^{(N)}\big)_{(i,j) \in \mathbb{N}^{2}}$ with probability weights $p^{(N)}_{u,o}$ and $p^{(N)}_{u,c}$ respectively. 
\begin{align}\label{Eq Z}
    &\mathbb{P}  \left( Z_{g+1,u}^{(N)}  =k | Z_{g,u}^{(N)} =m_N\right)  =
    \mathbb{P} \left( \sum_{i=1}^{m_N} X_i^{(N)} + \sum_{i,j=1, i>j}^{m_N} Y_{(i,j)}^{(N)} =k \right)
\notag \\ 
& = \sum_{\substack{k_o , k_c : \\ k_o + k_c=k, k_o, k_c \geq 0}}   
\mathbb{P} \left( \sum_{i=1}^{m_N} X_i^{(N)} =k_o \right)  \mathbb{P} \left( \sum_{i,j=1, i>j}^{m_N} Y_{(i,j)}^{(N)} =k_c \right) \notag \\
&= \sum_{\substack{k_o, k_c:\\ k_o + k_c=k, k_o, k_c \geq 0}} \mathbb{P}(W^{(N)}_{o, m_N} = k_o)  
\exp\left(-  m_N \frac{\ell_N^5 V_N^3}{D_N^2}  \right) \notag     \cdot
\mathbb{P}(W^{(N)}_{c, \binom{m_N}{2}} = k_c)
\exp\left(- 2\binom{m_N}{2} \frac{\ell_N^5 V_N^3}{D_N^2} \right) \notag \\  &= \exp\left(-  m^2_N \frac{\ell_N^5 V_N^3}{D_N^2}  \right) \sum_{\substack{k_o, k_c:\\ k_o + k_c=k, k_o, k_c \geq 0}} \mathbb{P}(W^{(N)}_{o, m_N} = k_o)  
  \cdot
\mathbb{P}(W^{(N)}_{o, m_N(m_N-1)} = k_c) \notag \\
&=  \exp\left(- m_N^2 \frac{\ell_N^5 V_N^3}{D_N^2} \right) \mathbbm{P}(W^{(N)}_{o, m_N^2} = k)  
\end{align}
where we have used that 

\[
	\begin{split}
		 \PP\left(\sum_{i=1}^{m_N} X_i^{(N)} = k\right) 
		&=\sum_{(k_1,..., k_{m_N}): k_1+\dots+k_{m_N}=k} 
  \left[\prod_{i=1}^{m_N} 
  \PP(W_{o,1}^{(N)}= k_i)
  \exp \left( -\frac{\ell_N^{5} V_{N}^{3}}{D_{N}^{2}}\right)\right]\\
		&=
  \exp \left( -m_N\frac{\ell_N^{5} V_{N}^{3}}{D_{N}^{2}}\right) \PP(W^{(N)}_{o, m_N}=k)
	\end{split}
\]
and a similar reasoning for $Y_{(i,j)}^{(N)}$.

Since $\exp \left(- \frac{m_N^2 V_N  ^2}{2 D_N}\right) \exp\left(- \frac{\ell_N^5 V_N^3}{D_N^2} \right) \geq  \exp\left(- m_N^2 \frac{\ell_N^5 V_N^3}{D_N^2}\right) \exp \left(- \frac{(m_N)^2 (V_N  - 2 \ell_N)^2}{2 D_N}\right)$  for $1\leq  m_N\leq \ell_N$ and $N$ large enough, we have that \eqref{Eq for D_k} $\geq$ \eqref{Eq Z} for $m_N \geq 2$ and \eqref{Eq for D_k 2}$\geq$ \eqref{Eq Z} for $m_N=1$ for $N$ large enough. Thus, because of the Markov property we can successively couple the two processes until $\overline{\ZZ}^{(N)}_u$ reaches the level $\ell_N$.
\end{proof}

Using the previous Section \ref{lab: DBPC section} we will show that for the upper DBPC defined in Definition \ref{Definition: UBPI} the probability of reaching an arbitrary high number of individuals up to a generation is asymptotically the same as the survival probability of a DBPC whose offspring and cooperation distributions are respectively $\text{Poi}\left(\frac{a^2}{2}\right)$ and $\text{Poi}(a^{2})$ distributed. 
\begin{Prop}\label{Proposition: asymptotic probability reaching certain level for BPI}(Probability for the total size of the upper DBPC to reach a level $b_N$). \\
Consider a sequence $\left(b_N\right)_{N \in \mathbb{N}}$ with $b_N \underset{N \to \infty}{\to}\infty$ and assume that $V_N \sim a \sqrt{D_N}$ for $0<a<\infty$. Then, we have 
\begin{equation}
    \lim_{N \to \infty} \mathbb{P} \left( \exists g \in \mathbb{N}_{0}: \overline{Z}^{(N)}_{g,u} \geq b_{N}\right)=\pi(a). 
\end{equation}
\end{Prop}

\begin{proof}
    The claim follows as an application of Proposition~\ref{Lemma: convergence survival probability DBPC}. Thus, we need to check that the sequence $\big(\ZZ^{(N)}_{u}\big)_{N\in \N}$ satisfies the assumption of Proposition~\ref{Lemma: convergence survival probability DBPC}. Let us first consider the convergence of $p_{j,u, o}^{(N)}\to p_j\big(\tfrac{a^2}{2}\big)$ for every $j\in\N_{0}$. Note that for a given $j$ we can choose $N$ large enough, such that $j\leq \ell_N$ and hence
    \begin{align}
        p_{j,u, o}^{(N)}= \PP(W^{(N)}_{o,1}=j) 
        \exp \left( -\frac{\ell_N^{5}V_{N}^{3}}{D_{N}^{2}}\right).
    \end{align}
    We now set $u_N:= \frac{(V_N - 2 \ell_N)^2}{2 D_N}$. By the choice of $\ell_N$ and since we assumed that $V_N\sim a \sqrt{D_N}$ we have $u_N\to \frac{a^2}{2}$. Thus, by continuity it follows that
    \begin{align}
       ( u_N)^{j} \frac{1}{j!} \exp ( -u_N) \exp \left( -\frac{\ell_N^{5}V_{N}^{3}}{D_{N}^{2}}\right)\underset{N \to \infty}{\to} \left( \frac{a^{2}}{2}\right)^{j} \frac{1}{j!} \exp \left( -\frac{a^2}{2}\right)=p_j\left(\frac{a}{\sqrt{2}}\right).
   \end{align}
 
Thus, we showed that $ p_{j,u, o}^{(N)}\to p_j\big(\tfrac{a^2}{2}\big)$ as $N\to \infty$ for every $j\in \N$. Analogously one can show that $ p_{j,u, c}^{(N)}\to p_j(a^2)$ as $N\to \infty$ for every $j\in \N$. Next we need to check that the first and second moment of the offspring and cooperation distribution converge.

Let $X^{(N)}$ be distributed according to the offspring distributions $\left(p^{(N)}_{j,u,o}\right)_{j\geq 0}$ of the upper DBPC $(Z_{n,u}^{(N)})$. Then 
\begin{align}
\E\left[X^{(N)}\right]& = \sum_{j=0}^{\ell_N+1} j  p_{j,u,o}^{(N)} \\
& = \sum_{j=1}^{\ell_N}\left(u_N \right)^j \frac{1}{(j-1)!} e^{-u_N} \exp \left( -\frac{\ell_N^{5}V_{N}^{3}}{D_{N}^{2}}\right)  + (\ell_{N} +1)\bigg( 1- \sum_{j=0}^{\ell_N} p_{j,u,o}^{(N)}\bigg). 
\end{align}

Since $l_N\to \infty$,
we have
\begin{equation*}
     \lim_{N\to \infty}e^{-u_N} \sum_{j=0}^{\ell_N-1}\frac{(u_N)^j }{j!}=  1.
\end{equation*}
It follows that 
\begin{align}
&\exp \left( -\frac{\ell_N^{5}V_{N}^{3}}{D_{N}^{2}}\right) \sum_{j=1}^{\ell_N} \left(u_N \right)^j \frac{1}{(j-1)!} e^{- u_N} = u_N\exp \left( -\frac{\ell_N^{5}V_{N}^{3}}{D_{N}^{2}}\right) \left(\sum_{j=0}^{\ell_N-1} \frac{\left(u_N \right)^j }{j!} e^{- u_N}\right)\rightarrow \frac{a^2}{2}
\end{align}
and
\begin{align}
(\ell_{N} +1)\left( 1- \sum_{j=0}^{\ell_N} p_{j,u,o}^{(N)} \right) & = 
(\ell_N+1) \left(1 - \exp \left( -\frac{\ell_N^{5}V_{N}^{3}}{D_{N}^{2}}\right) \mP\left(W_{o,1}^{(N)} \leq \ell_N\right)\right). 
\end{align} Now by Markov's inequality follows that 
\begin{equation*}
    \mP\left(W_{o,1}^{(N)} \leq \ell_N\right)=\mP\left((W_{o,1}^{(N)})^2 \leq \ell_N^2\right)\geq 1-\frac{\E[(W_{o,1}^{(N)})^2]}{\ell_N^2}=1-\frac{u_N+u_N^2}{\ell_N^2}.
\end{equation*}
Hence 
\begin{align*}
    (\ell_{N} +1)\left( 1- \sum_{j=0}^{\ell_N} p_{j,u,o}^{(N)} \right)
    &\leq (\ell_N+1) \left(1 - \exp \left( -\frac{\ell_N^{5}V_{N}^{3}}{D_{N}^{2}}\right) \right) +(u_N+u_N^2)\frac{\ell_N+1}{\ell^2_N}\\
    &\leq  \frac{(\ell_N+1)\ell_N^{5}V_{N}^{3}}{D_{N}^{2}} +(u_N+u_N^2)\frac{\ell_N+1}{\ell^2_N} \rightarrow 0.
\end{align*}

Consequently 
\begin{equation*}
    \E\left[X^{(N)}\right] \rightarrow\frac{a^2}{2}. 
\end{equation*}
as $N\to \infty$. Similarly, we have for the second moment
\begin{align}
\E\Big[\big(X^{(N)}\big)^2\Big]
= \sum_{j=0}^{\ell_N} j^2\left(u_N \right)^j \frac{1}{j!} e^{-u_N} \exp \left( -\frac{\ell_N^{5}V_{N}^{3}}{D_{N}^{2}}\right)  + (\ell_{N} +1)^2\bigg( 1- \sum_{j=0}^{\ell_N} p_{j,u,o}^{(N)}\bigg).
\end{align}
The second term again vanishes in the limit  by the same argument as before just that we use the Markov inequality for the third moment such that $ \PP(W^{(N)}_{o,1} \leq \ell_N)\geq 1-\frac{u_N+3u_N^2+u_N^3}{\ell_N^3}$, which yields that
\begin{align*}
    (\ell_{N} +1)^2\left( 1- \sum_{j=0}^{\ell_N} p_{j,u,o}^{(N)} \right)
    \leq  \frac{(\ell_N+1)^2\ell_N^{5}V_{N}^{3}}{D_{N}^{2}}+(u_N+3u_N^2+u_N^3)\frac{(\ell_N+1)^2}{\ell^3_N}\to 0,
\end{align*}
as $N\to \infty$ and 
\begin{align*}
   \exp \left( -\frac{\ell_N^{5}V_{N}^{3}}{D_{N}^{2}}\right) \sum_{j=0}^{\ell_N} j^2 \frac{\left(u_N \right)^j}{j!}
   e^{-u_N}
   =\exp \left( -\frac{\ell_N^{5}V_{N}^{3}}{D_{N}^{2}}\right) u_N \left(u_N \sum_{j=0}^{\ell_N-2} \frac{ \left(u_N \right)^j}{j!}e^{-u_N}+\sum_{j=0}^{\ell_N-1} \frac{\left(u_N \right)^j}{j!}e^{-u_N}\right)
\end{align*}
Now one can show analogously as before that  
\begin{equation*}
    \E\Big[\big(X^{(N)}\big)^2\Big] \rightarrow \frac{a^2}{2} + \frac{a^4}{4}.
\end{equation*}

For the expectations and the second moments of the cooperation distributions one argues analogously, except that one shows convergence to $a^2$  and $a^2 + a^4 $, respectively.
\end{proof}

\subsection
{Lower bound on the invasion probability on a complete graph}

\subsubsection{Lower bound on the probability to infect at least $N^\eps$ many hosts}

We first aim to show that the total number of infected hosts until the parasite population dies out or $N^\varepsilon$ hosts are infected for $\varepsilon >0$ small enough can be lower bounded by the total size process of a  DBPC. This DBPC we introduce next.

\begin{Def}(Lower discrete branching process with cooperation) \label{Definition:LDBPC} \\
Let $\ell_{N} \underset{N \to \infty}{\to}\infty$ and $\frac{1}{2}<\delta <1$ satisfying $\frac{\ell_N^4 V_N \log(\log(N))}{D_N^\delta}\in o(1)$. Let $\ZZ^{(N)}_{\ell}=\left(Z_{g,\ell}^{(N)}\right)_{g \in \mathbb{N}_0}$ be a branching process with interaction with $Z_{0,\ell}^{(N)}=1$ almost surely, and offspring and cooperation distributions  $p_{\ell,o}^{(N)} = \left( p_{j,\ell,o}^{(N)}\right)_{j \in \mathbb{N}_{0}}$ and $p_{\ell,c}^{(N)} = \left( p_{j,\ell,c}^{(N)}\right)_{j \in \mathbb{N}_{0}}$ with 
\begin{align}
    p_{j,\ell,o}^{(N)}:= \PP(W^{(N)}_{o,1}=j)
    \exp \left( -3\frac{ \ell_N V_{N}}{D_N^\delta}\right),
\end{align}
for all $0 < j \leq \ell_N$ 
and 
\begin{align}
p_{0,\ell,o}^{(N)}:=1-\sum_{j=1}^{\ell_N} p_{j,\ell,o}^{(N)},
\end{align}
as well as 
\begin{align}
    p_{j,\ell,c}^{(N)}:= \PP(W_{c,1}^{(N)}=j)  
    \exp \left( -3 \frac{\ell_N V_N}{D_N^\delta}\right),
\end{align}
for all $0 < j \leq \ell_N$ and 
\begin{align}
    p_{0,\ell,c}^{(N)}:=1-\sum_{j=1}^{\ell_N} p_{j,\ell,c}^{(N)}. 
\end{align}

Denote by $\overline{\ZZ}_{\ell}^{(N)}:= \left( \overline{Z}_{g,\ell}^{(N)}\right)_{g \in \mathbb{N}_{0}}$ where $\overline{Z}_{g,\ell}^{(N)}:=\sum_{i=0}^{g}Z_{i,\ell}^{(N)}$, that is $\overline{Z}_{g,\ell}^{(N)}$ gives the total size of $\ZZ^{(N)}_{\ell}$ accumulated till generation $g$.
\end{Def}

\begin{Prop}\label{Coupling Lower DBPC}
Consider a sequence $\left(\ell_N \right)_{N \in \mathbb{N}}$ and $\frac{1}{2}< \delta <1$ with $\ell_N=N^{\varepsilon} \underset{N \to \infty}{\to}\infty$ for $\varepsilon>0$ such that $\frac{\ell_N^5 V_N \log(\log(N))}{D_N^{\delta}}\in o(1)$. Introduce the stopping time
\begin{align}
    \sigma^{(N)}_{\ell_N,0}:= \inf \left\{ g \in \mathbb{N}_{0}: \overline{Z}^{(N)}_{g,\ell} \geq \ell_{N} \text{ or } Z_{g,\ell}^{(N)}=0\right\}.
\end{align}
Then 
\begin{align}\label{CouplingLowerDBPC}
    \lim_{N \to \infty} \mathbb{P} \left( \overline{J}^{(N)}_{g} \geq \overline{Z}^{(N)}_{g,\ell}, \forall g \leq \sigma^{(N)}_{\ell_N, 0}\right)=1. 
\end{align}
\end{Prop}

\begin{proof}
As in the proof of Proposition \ref{Prop: Coupling with upper DBPC} in the scaling $\frac{\ell_N^{3}V_N^3}{D_N^2}=o(1)$ it is unlikely, that up to the first generation, at which the total infection process reaches size $\ell_N$, an infection occurs due to more than a pair of parasites. Consequently for a coupling whp we can only focus on infections generated by pairs of parasites, and do not need to treat infections generated by at least 3 parasites. \\
Now consider a complete graph with exactly $m_N <\ell_N$ infected vertices and with at most $\ell_N-1$ empty or infected vertices. The number of new infections generated on such a graph in the next generation can be lower bounded by the number of infections arising in the following experiment: Consider $D_N$ boxes and $m_N V_N$ balls. Distribute the balls uniformly at random into the boxes. Assume a new infection is created for each of the first $D_N-\ell_N$ boxes that contains at least two balls. Let $A_{m_N}$ be the number of infections generated in this experiment and let $B_k$ be the event that exactly $k$ of the $D_N-\ell_N$ first boxes contain exactly two balls and all other boxes have at most one ball. We have 
\begin{align*}
\mathbbm{P}\left( J_g^{(N)} \geq A_m | J_{g-1}^{(N)} =m  \right)=1
\end{align*}
and 
\begin{align*}
\mathbbm{P}(A_m =k) \geq \mathbbm{P}(B_k).
\end{align*}

We will show below that there exists a constant $C_1>0$ such that for $N$ large enough 
\begin{align}\label{ComparisonCZ}
&\left| \mathbbm{P}\left(B_k\right) - 
\mathbbm{P}\left(Z^{(N)}_{g+1, \ell} =k | Z_{g,\ell}^{(N)} =m \right)\right| \\
&\hspace{3cm}\leq \left(\frac{m_N^{2} V_N^{2}}{2 D_N}\right)^k \frac{1}{k!} \exp\left(- \frac{m_N^2 V_N^2}{2D_N}\right) 
\left( \frac{C_1 V_N \ell_N^4}{D_N^\delta} \right) + \frac{C_1 \ell_N^4 V_N}{D_N^\delta}
\end{align}
for any $k$ with $0\leq k \leq \ell_N$.

Since 
\[ 
\sum_{k=0}^{\ell_N} \left( \left(\frac{m_N^{2} V_N^{2}}{2 D_N}\right)^{k} \frac{1}{k!} \exp\left(- \frac{m_N^2 V_N^2}{2D_N}\right) 
\left( \frac{C_1 V_N \ell_N^4}{D_N^\delta} \right) + \frac{C_1 \ell_N^4 V_N}{D_N^\delta}\right)  \leq  \frac{2 C_1 \ell_N^5 V_N}{D_N^\delta} \rightarrow 0
\]
we can couple the balls into boxes experiment with the lower DBPC $\mathbf{Z}^{(N)}_\ell$, such that given that $\{Z_{g, \ell}^{(N)}=m\}$ the event $B_k^{(N)}$ occurs whp, if $\{Z_{g+1, \ell}^{(N)}=k\}$ for any $k\in \{0, ... \ell_N\}$ and vice-versa.
We can repeat this argument till $\sigma^{(N)}_{\ell_N,0}$ whp. Indeed, by Proposition~\ref{GrowFast} it exists $C_2>0$ such that $\mP\left(\sigma^{(N)}_{\ell_N, 0}\leq C_2\log(\log(\ell_N))\right) \rightarrow 1 $, as by analogous arguments as in the proof of Proposition \ref{Proposition: asymptotic probability reaching certain level for BPI} it can be shown that the first and second moment of $Z^{(N)}_{1,\ell}$ is uniformly bounded in $N$. 
Since 
\[\left(1-\frac{2 C_1 \ell_N^5 V_N}{D_N^\delta} \right)^{C_2 \log \log(\ell_N)} \rightarrow 1\]
it follows that we can couple whp subsequently performed balls into boxes problems and $\mathbf{Z}_\ell^{(N)}$ for any generation $g \leq \sigma^{(N)}_{\ell_N,0}$, which implies \eqref{CouplingLowerDBPC}. 

So to finish the proof it remains to show \eqref{ComparisonCZ}.
 
We start by controlling the probabilities of the events $B_k$. By Lemma \ref{lab: balls and boxes experiment} with $D'_N = D_N$, $m'_N= m_N$, $V'_N= V_N$ and $h'_N= \ell_N$
we can estimate 
\begin{align}\label{Estimate B_k upper}
\mathbb{P} \left(B_k\right) \leq 
 \left( \frac{(m_N V_N)^2}{2 D_N} \right)^k \frac{1}{k!} \exp\left(-\frac{(m_N V_N - 2 \ell_N)^2}{2 D_N} \right) \exp\left( \frac{{\ell}_N^2 V_N}{D_N} \right),
\end{align}
and 
\begin{align}\label{Estimate B_k lower}
\mathbb{P} \left(B_k\right)\geq 
\mathbb{P}(W^{(N)}_B=k)\exp \left( -\frac{\ell_N^5 V_N^3}{D_N^2}\right) 
\end{align}

for a Poisson distributed random variable $W^{(N)}_B$ with parameter $\frac{(m_N V_N - 2 \ell_N)^2}{2 D_N}$ and  $N$ large enough. 

Next we control the transition probabilities of $\mathbf{Z}_\ell^{(N)}$.
Consider independent random variables $\left(X_i^{(N)}\right)_{i \in \mathbb{N}}$ and $\left(Y_{(i,j)}^{(N)}\right)_{i<j}$ with probability weights $p_{\ell,o}^{(N)}$ and $p_{\ell,c}^{(N)}$ respectively.

We have 
\begin{align*} 
& \mP \left( X_1^{(N)}=0\right) = 1 - \sum_{j=1}^{\ell_{N}} p_{j, \ell,o}^{(N)} \\ &
=  1- \sum_{j=1}^{\ell_{N}} \exp\left(-3 \frac{\ell_N V_N}{D_N^{\delta}} \right) \mP\left(W_{o,1}^{(N)} =j\right)\\  &
= 1 + \exp\left(-3 \frac{\ell_N V_{N}}{D_N^\delta} \right) \exp(-u_N) -  \sum_{j=0}^{\ell_N} \exp\left(-3 \frac{\ell_N V_N}{D_N}\right) \mP\left(W_{o,1}^{(N)} =j\right)  \\ 
& = \exp\left(-3 \frac{\ell_N V_N }{D_N^\delta }\right) \exp(-u_N)  + \exp\left(-3 \frac{\ell_N V_N}{D_N^\delta}\right) \mP\left(W_{o,1}^{(N)}\geq \ell_N +1\right)+1-\exp\left(-3 \frac{\ell_N V_N}{D_N^\delta}\right).
\end{align*}

We define the constant $c^{(N)}_o:=  \exp\big(-3 \frac{\ell_N V_N}{D_N^\delta}\big) \mP\big(W_{o,1}^{(N)}\geq \ell_N +1\big)+1-\exp\big(-3 \frac{\ell_N V_N}{D_N^\delta}\big)$. We have $\mP\big(W_{o,1}^{(N)} \geq \ell_N +1\big) \in  \mathcal{O}\Big(\frac{u_N^{\ell_N}}{(\ell_N)!} \Big) $ so $c_o^{(N)} \in \Theta\Big(\frac{\ell_N V_N}{D_N^\delta}\Big)$, where we used that $\ell_N\sim N^{\varepsilon}$, and thus $\frac{u_N^{\ell_N}}{(\ell_N)!}$ decays exponentially fast in $N$.

Let us recall that by definition  
\begin{align*}
    \mP\left(X_i^{(N)}= k\right)= \frac{\left(u_N\right)^{j}}{j!} e^{-u_N} \exp \left( -3\frac{ \ell_N V_{N}}{D_N^\delta}\right)+c_o^{(N)}\1_{\{k=0\}}.
\end{align*}
$0\leq k \leq \ell_N$. We see that for $0\leq k_o \leq \ell_N$ it holds that
\begin{align*} 
 & \mP \left( \sum_{i=1}^{m_N} X_i^{(N)} =k_o \right)   =\sum_{\stackrel{k_1, ..., k_m:}{ k_1 +...+ k_m=k_o}} \prod_{i=1}^{m_N} \mP\left(X_i^{(N)}= k_i\right),
\end{align*}
which allows is to derive the lower and upper bound
\begin{align}\label{EstimateX}
\mP \left( \sum_{i=1}^{m_N} X_i^{(N)} =k_o \right)  &\geq \exp\left(-3 m_N \frac{\ell_N V_N}{D_N^{\delta}}\right)
\mP \left( W_{o, m_N}^{(N)} = k_o \right) \quad \text{ and }\\
\mP \left( \sum_{i=1}^{m_N} X_i^{(N)} =k_o \right)  &\leq \exp\left(-3 m_N \frac{\ell_N V_N}{D_N^{\delta}}\right)
\mP \left(  W_{o, m_N}^{(N)} = k_o \right) + m_N c_{o}^{(N)}.
\end{align}
where $ W_{o, m_N}^{(N)}\sim \text{Poi}(m_N u_N)$ and if $k_1+ ...+ k_{m_N} =k_o$, then the number of $k_i$ with $k_i=0$ is at most $m_N$. 

Now we obtain analogously as before that 
\begin{align*}
\PP\left(Y_{1,2}^{(N)} =0\right) = \exp\left(- 3\frac{\ell_N V_N}{D_N^\delta}\right) \exp(-u_N) + c_c^{(N)}
\end{align*}
where the constant is defined as
\[ c^{(N)}_c:= \exp\left(-3 \frac{\ell_N V_N}{D_N^\delta}\right) \PP\left(\tilde{Y}^{(N)} \geq \ell_N +1\right) + 1- \exp\left(-3 \frac{\ell_N V_N}{D_N^\delta}\right)\in \Theta\left(\frac{\ell_N V_N}{D_N^\delta}\right), \]
with $\tilde{Y}^{(N)}\sim \text{Poi}(2u_N)$. Similarly as before we arrive at the lower and upper bound

\begin{align}\label{estimateY}
     \mP\left( \sum_{i,j=1, i<j}^{m_N} Y_{(i,j)}^{(N)} =k_c \right ) 
    &\geq \mP\left(W^{(N)}_{c, \binom{m_N}{2}}= k_c \right)  \exp\left(-3 \binom{m_N}{2} \frac{\ell_N V_N}{D_N^\delta}\right)\\ 
    \mP\left( \sum_{i,j=1, i<j}^{m_N} Y_{(i,j)}^{(N)} =k_c \right ) 
    &\leq \mP\left(W^{(N)}_{c, \binom{m_N}{2}}= k_c \right)  \exp\left(-3 \binom{m_N}{2} \frac{\ell_N V_N}{D_N^\delta}\right) 
    + \binom{m_N}{2}c_c^{(N)}.
\end{align}

 So in summary we have
\begin{align}
\mathbb{P}\left(Z_{g+1,\ell}^{(N)} =k | Z_{g,\ell}^{(N)} =m_N\right) & =
    \mathbb{P} \left( \sum_{i=1}^{m_N} X_i^{(N)} + \sum_{i,j=1, i<j}^{m_N} Y_{(i,j)}^{(N)} =k \right)
\notag \\ 
& = \sum_{k_o , k_c: k_o + k_c=k, k_o, k_c \geq 0}   
\mathbb{P} \left( \sum_{i=1}^{m_N} X_i^{(N)} =k_o \right)  \mathbb{P} \left( \sum_{i,j=1, i<j}^{m_N} Y_{(i,j)}^{(N)}
=k_c \right) 
\end{align}
and hence using  \eqref{EstimateX} and \eqref{estimateY}
 for any $0\leq k\leq \ell_N$ we have for an appropriate constant $C>0$ 
\begin{align}\label{EstimateZ}
 \PP(W^{(N)}_{o, m_N^2} =k) \exp\left(-3 m_N^2 \frac{\ell_N V_N}{D_N^\delta} \right) & \leq \mathbb{P}\left(Z_{g+1,\ell}^{(N)} =k | Z_{g,\ell}^{(N)} =m_N\right) \\ & \leq \PP(W^{(N)}_{o, m_N^2} =k) \exp\left(-3 m_N^2 \frac{\ell_N V_N}{D_N^\delta} \right) + \ell_N^3 c_c^{(N)} + \ell_N^2 c_o^{(N)} +  \ell_N^4 c_c^{(N)} c_o^{(N)} \\
 & \leq  \PP(W^{(N)}_{o, m_N^2} =k) \exp\left(-3 m_N^2 \frac{\ell_N V_N}{D_N^\delta} \right) + \frac{C \ell_N^4 V_N}{D_N^\delta}
\end{align}
 Subtracting upper and lower, resp., bounds of the transition probabilities of $\mathbf{Z}^{(N)}_\ell$ from the lower and upper, resp. bounds of $\mP(B_k)$ and taking the modulus yields \eqref{ComparisonCZ}. Indeed,
by \eqref{Estimate B_k lower} and \eqref{EstimateZ}
we have for a constant $C>0$ that may change from line to line
\begin{align*}
&\mP(Z_{g+1, \ell}^{(N)} =k | Z_{g, \ell}^{(N)} = m_N) - \mP(B_k)\\& \leq 
\mP(W^{(N)}_{o, m_N^2}=k)\left(\exp\left(- 3m_N^2 \frac{\ell_N V_N}{D_N^\delta}\right) - \exp\left(- \frac{\ell_N^5 V_N^3}{D_N^2} \right)\right) + \frac{C \ell_N^4 V_N}{D_N^\delta} \\
&\leq \left(\frac{(m_N V_N)^2}{2D_N}\right)^k \frac{1}{k!} \exp\left(-\frac{ (m_N V_N)^2}{2D_N} \right)\frac{\ell_N^4 V_N}{D_N^\delta} +  \frac{C \ell_N^4 V_N}{D_N^\delta}
\end{align*}
since $\exp\left(-\frac{m_N^2(V_N - 2\ell_N)^2}{2D_N} \right) \leq \exp\left(- \frac{m_N^2 V_N^2}{2 D_N}\right) \exp\left(\frac{2m_N^2 \ell_N V_N}{D_N} \right) \leq  \exp\left(- \frac{m_N^2 V_N^2}{2 D_N}\right)\left(1 + \frac{ \ell_N^4 V_N}{D_N} \right) $ for $N$ large enough
 and for an appropriate constant $C>0$ (that may differ from the constant $C$ used above) for $N$ large enough. Furthermore, we see that
 \begin{align*}
\mP(W^{(N)}_{o, m_N^2}=k) &= \left(\frac{m_N^2 (V_N-2\ell_N)^2}{2 D_N} \right)^k \frac{1}{k!} \exp\left( - \frac{m_N (V_N-2\ell_N)^2}{2D_N}\right)\\ 
&\geq \left(\frac{m_N^2 V_N^2}{2 D_N} \right)^k \left(1-\frac{2\ell_N}{ V_N} \right)^{2k}\frac{1}{k!} \exp\left( - \frac{m_N (V_N-2\ell_N)^2}{2D_N}\right)\\
& \geq 
\left(\frac{m_N^2 V_N^2}{2 D_N} \right)^k \frac{1}{k!} \exp\left( - \frac{m_N (V_N-2\ell_N)^2}{2D_N}\right) - \frac{C \ell_N^4 V_N}{D_N^\delta},
\end{align*}
where we used again Bernoulli's inequality in the second inequality. Now we have for $N$ large enough by \eqref{Estimate B_k upper} and \eqref{EstimateZ}
\begin{align*}
&\mP(B_k) -\mP(Z_{g+1, \ell}^{(N)} =k | Z_{g, \ell}^{(N)} = m_N)\\ & \leq
\left(\frac{m_N^2 V_N^2}{2 D_N}\right)^k \frac{1}{k!} \exp\left( - \frac{m_N^2 (V_N-2\ell_N)^2}{2D_N}\right) \left( \frac{2\ell_N^2 V_N}{D_N} + \frac{6 \ell_N^3 V_N}{D_N^\delta}\right) \\ 
& \quad  \quad \quad + \frac{C \ell_N^4 V_N}{D_N^\delta} \\
& \leq  \left(\frac{m_N^2 V_N^2}{2 D_N}\right)^k \frac{1}{k!} \exp\left( - \frac{m_N^2 V_N^2}{2D_N}\right)  \frac{7 \ell_N^3 V_N}{D_N^\delta} + \frac{C \ell_N^4 V_N}{D_N^\delta}. 
\end{align*}
This yields the claim.

\end{proof}

As a counterpart of Proposition \ref{Proposition: asymptotic probability reaching certain level for BPI}, we show that the total size of the lower DBPC of Definition \ref{Definition:LDBPC} reaches a level tending to infinity with asymptotically the survival probability of a DBPC whose offspring and cooperation distributions are respectively $\text{Poi}\left(\frac{a^2}{2}\right)$ and $\text{Poi}(a^2)$ distributed.
\begin{Prop}\label{Asymptotic survival probability lower DBPC}(Probability for the Total Size of the Lower BPI to Reach a Level $k_N$). \\
Consider a sequence $\left(k_N\right)_{N \in \mathbb{N}}$ with $k_N \underset{N \to \infty}{\to}\infty$ and $V_N\sim a\sqrt{D_N}$. Then, we have 
\begin{align}
    \lim_{N \to \infty} \mathbb{P} \left( \exists g \in \mathbb{N}_{0}: \overline{Z}^{(N)}_{g,l} \geq k_{N}\right)=\pi(a). 
\end{align}
\end{Prop}

\begin{proof}
This Proposition is shown by the same line of argument as Proposition~\ref{Proposition: asymptotic probability reaching certain level for BPI}, i.e.\ basically one applies Lemma~\ref{Lemma: convergence survival probability DBPC}.
\end{proof}

\subsubsection{Final phase of an epidemic on a complete graph}
\label{Subsection: Final phase Epidemic complete graph}
In this subsection we are going to show that if the total size of the infection process reaches the level $N^{\varepsilon}$ for $\varepsilon>0$, then in a finite number of generations, all the hosts are killed. For that we will intensively use the following Lemma. 
\begin{Lem}
\label{Lem: experiment boxes and balls}
Let $\varphi_1(N), \varphi_2(N)$ and $\varphi_{3}(N)$ such that $\varphi_i(N)=o(D_N)$ for $i \in \{1,2\}$ and $\varphi_3(N)=o(V_N)$. Let $H(N)$ such that it exists $\varepsilon>0$ such that $N^{\varepsilon}=o(H(N))$ and $H(N)=o\left(\sqrt{D_N}\right)$. Consider the following experiment: Assume $H(N)(V_N-\varphi_{3}(N))$ balls are distributed purely at random into $D_N-\varphi_1(N)$ boxes. Denote by $G^{(N)}$ the number of boxes among the first $D_N-\varphi_1(N)-\varphi_2(N)$ boxes that contain at least 2 balls. Then it holds:
\begin{itemize}
    \item[(i)]  Define $\ell:= \inf \left\{i \geq 2: H(N)^{i+1}=o\left(\sqrt{D_N}^{i-1}\right) \right\}$. Assume $f_1(N)$ satisfies $\log(H(N)V_N)=o(f_1(N))$. Then we have 
    \begin{align}
    \label{Eq: min number pairs of balls}
        \mathbb{P} \left(\frac{H^{2}(N)}{f_1(N)} \leq G^{(N)}\right) \geq 1-\Theta\left(\frac{H(N)^{\ell+1}}{\sqrt{D_N}^{\ell-1}}\right). 
    \end{align}
    
    \item[(ii)] Let $f_2(N) \to \infty$. Then we have 
\begin{align}
\label{Eq: max number pair of balls}
    \mathbb{P} \left( G^{(N)} \leq H^{2}(N)f_2(N)\right) \geq 1-\Theta \left(f_{2}^{-1}(N)\right). 
\end{align}
\end{itemize} 
\end{Lem}
A proof can be found in the supplementary material. \\ 
Now, introduce
\begin{align}
    &\overline{\tau}^{(N)}:= \inf \left\{g \in \mathbb{N}_{0}: \overline{J}_{g}^{(N)} \geq N^{\varepsilon}\right\}, \\
    &\tau^{(N)}:= \inf \left\{g \in \mathbb{N}_0: J^{(N)}_{g} \geq \frac{N^{\varepsilon}}{\log(N)}\right\}.
\end{align}
\begin{Prop}
\label{Prop: reached certain level if total size reaches one}
\begin{align}
    \mathbb{P} \left( \tau^{(N)} \leq \overline{\tau}^{(N)} \vert \overline{\tau}^{(N)} <\infty \right) \underset{N \to \infty }{\to}1. 
\end{align}
\end{Prop}
\begin{proof}
If for any generation before $\overline{\tau}^{(N)}$ the number of infected vertices is strictly smaller than $\frac{N^{\varepsilon}}{\log(N)}$ then this would mean that the number of generations until the total size of the infection process reaches the level $N^{\varepsilon}$  is at least $\log(N)$. But this contradicts the fact that it exists a constant $C>0$ such that 
\begin{align}
    \mathbb{P} \left( \overline{\tau}^{(N)} \leq C\log(\log(N)) \vert \overline{\tau}^{(N)}<\infty\right) \to 1, 
\end{align}
which follows from coupling from below with the DBPC of Definition \ref{Definition:LDBPC} and Proposition \ref{GrowFast}.  
\end{proof}

\begin{Lem}
\label{Lem: first step final phase}
We have 
\begin{align}
    \mathbb{P} \left( J^{(N)}_{\tau^{(N)}}\geq \frac{N^{ \varepsilon}}{\log(N)}, \overline{J}^{(N)}_{\tau^{(N)}}\leq \frac{N^{2 \varepsilon}}{\log(N)} \Big{\vert} \overline{\tau}^{(N)}<\infty\right) \underset{N \to \infty}{\longrightarrow}1. 
\end{align}
\end{Lem}
\begin{proof}
By definition at generation $\tau^{(N)}$ the number of infected vertices is at least $\frac{N^{\varepsilon}}{\log(N)}$. At generation $\tau^{(N)}-1$, whp we have $\overline{J}^{(N)}_{\tau^{(N)}-1} \leq N^{\varepsilon}$ because otherwise we have a contradiction to Proposition \ref{Prop: reached certain level if total size reaches one}. Then to bound from above the total number of infected vertices up to generation $\tau^{(N)}$, it suffices to add to $N^{\varepsilon}$ an upper bound on the number of new infections generated at generation $\tau^{(N)}$. This upper bound is obtained by an application of Lemma \ref{Lem: experiment boxes and balls} with $H(N)=\frac{N^{\varepsilon}}{\log(N)}$, $\varphi_1(N)=N^{\varepsilon}$, $\varphi_2(N)=0$, $\varphi_3(N)=0$, $f_2(N)=\frac{1}{2}\log(N)$ and an arbitrary function $f_1$  satisfying the condition of Lemma \ref{Lem: experiment boxes and balls}. \\
Indeed, since before generation $\tau^{(N)}$ the total number of parasites on the graph is at most $\frac{N^{ \varepsilon}}{\log(N)}v_N$, the number of new infections generated is controlled from above using the previous experiment.
\end{proof}
Next choose $\varepsilon>0$ such that for all $k \in \mathbb{N}$,  $2^{k}\varepsilon\neq\frac{\beta}{2}$. Then define $\overline{k}$ as the largest $k \in \mathbb{N}$ satisfying  $2^{k+1}\varepsilon<\beta$. In particular it holds $2^{\overline{k}+1}\varepsilon >\frac{\beta}{2}$ because otherwise $2^{\overline{k}+2} \varepsilon <\beta$  which contradicts the definition of $\overline{k}$.

\begin{Lem}
\label{Lemma: control infected and empty vertices final phase}
Let $k \in \llbracket 0,\overline{k} \rrbracket$. We have 
\begin{align}
    \mathbb{P} \left( J^{(N)}_{\tau^{(N)}+k} \geq \frac{N^{2^{k}\varepsilon}}{\log^{\alpha_k}(N)}, \overline{J}^{(N)}_{\tau^{(N)}+k} \leq \frac{N^{2^{k+1}\varepsilon}}{\log(N)}\Big{\vert} \overline{\tau}^{(N)}<\infty\right) \to 1,
\end{align}
where we set $\alpha_0:=1$, and for all $k \geq 1$ we set $\alpha_k:=2 \alpha_{k-1}+2$. 
\end{Lem}
\begin{proof}
We prove the claim via induction over $k$. For $k=0$ the claim follows by Lemma \ref{Lem: first step final phase}.
Next we prove the claim for $k+1$ assuming the claim holds for all $0 \leq j \leq k$.  \\ \\
For the lower bound on the number of infected vertices at generation $\tau^{(N)}+k+1$, apply Lemma \ref{Lem: experiment boxes and balls} with $H(N)=\frac{N^{2^{k}\varepsilon}}{\log^{\alpha_k}(N)}$, $\varphi_1(N)=0$, $\varphi_2(N)=\frac{N^{2^{k+1}\varepsilon}}{\log(N)}$, $\varphi_3(N)=0$, $f_1(N)=\log^{2}(N)$ and an arbitrary function $f_2$ with $f_2(N)\to \infty$, which yields that the number of infected vertices at generation $\tau^{(N)}+k+1$ is whp at least of order $\frac{1}{\log^{2}(N)}\left(\frac{N^{2^{k}\varepsilon}}{\log^{\alpha_k}(N)}\right)^{2}=\frac{N^{2^{k+1}\varepsilon}}{\log^{\alpha_{k+1}}(N)}$. \\
Indeed by considering $D_N$ boxes we lower bound the probability for a parasite to attack an occupied vertex, which is $\frac{1}{D_N-1}$ in the case of the complete graph. According to the induction hypothesis we have considered whp  by Lemma \ref{Lem: experiment boxes and balls} the minimal number of parasites which is $\frac{N^{2^{k}\varepsilon}}{\log^{\alpha_k}(N)}V_N$. In the balls into boxes experiment new infections are (only) counted when reaching one of the $D_N-\frac{N^{2^{k+1}\varepsilon}}{\log(N)}$ first boxes whereas in the original process there are at least this number of occupied vertices. 
\\ \\ 
To arrive at the upper bound on the number of empty vertices, apply Lemma \ref{Lem: experiment boxes and balls} with $H(N)=\frac{N^{2^{k+1}\varepsilon}}{\log(N)}$, $\varphi_1(N)=H(N)$, $\varphi_2(N)=0$, $\varphi_3(N)=0$, $f_2(N)=\log(N)$ and an arbitrary function $f_1$ that satisfies the conditions of Lemma \ref{Lem: experiment boxes and balls}, since in the previous upper bound the number of empty vertices is bounded by $\frac{N^{2^{k+1}\varepsilon}}{\log(N)}$. So according to Lemma \ref{Lem: experiment boxes and balls}  the number of empty vertices at generation $\tau^{(N)}+k+1$ is whp at most of order $\log(N)\left(\frac{N^{2^{k+1}\varepsilon}}{\log(N)}\right)^{2}=\frac{N^{2^{k+2}\varepsilon}}{\log(N)}$.
\end{proof}

Applying Lemma \ref{Lemma: control infected and empty vertices final phase} with $k=\overline{k}$ we obtain 
\begin{align}
    \mathbb{P}\left(J_{\tau^{(N)}+\overline{k}}^{(N)} \geq \frac{N^{2^{\overline{k}} \varepsilon}}{\log^{\alpha_{\overline{k}}}(N)}, \overline{J}_{\tau^{(N)}+\overline{k}}^{(N)}\leq \frac{N^{2^{\overline{k}+1}\varepsilon}}{\log(N)} \vert \overline{\tau}^{(N)}<\infty \right) \to 1.
\end{align}
Define $\delta=\frac{1}{2}\left(2^{\overline{k}+1}\varepsilon-\frac{\beta}{2}\right)>0$. In the next Lemma we show that at generation $\overline{k}+1$ the number of infected vertices is at least of order $N^{\frac{\beta}{2}+\delta}$. 
\begin{Lem}
\label{Lem: beta over 2 + delta}
\begin{align}
    \mathbb{P}\left(J_{\tau^{(N)}+\overline{k}+1}^{(N)} \geq N^{\frac{\beta}{2}+\delta} \vert \overline{\tau}^{(N)}<\infty \right) \to 1.
\end{align}
\end{Lem}
\begin{proof}
Here we apply again Lemma \ref{Lem: experiment boxes and balls} to obtain this lower bound. More precisely with the following set of parameters: $H(N)=\frac{N^{2^{\overline{k}}\varepsilon}}{\log^{\alpha_{\overline{k}}}(N)}$, $\varphi_1(N)=0$, $\varphi_2(N)=\frac{N^{2^{\overline{k}+1}\varepsilon}}{\log(N)}$, $\varphi_3(N)=0$, $f_1(N)=\log^{2}(N)$ and an arbitrary function $f_2$. We obtain that whp $I_{\tau^{(N)}+\overline{k}+1}\geq \frac{N^{2^{\overline{k}+1}\varepsilon}}{\log^{\alpha_{\overline{k}+1}}(N)}\geq N^{\frac{\beta}{2}+\delta}$, by definition of $\delta$. 
\end{proof}

In the next lemma we show that in one more generation whp any vertex will be reached by at least 2 parasites, in other words each of the remaining hosts gets infected whp. 
\begin{Lem}
\label{Lem: every vertices are killed}
\begin{align}
    \mathbb{P}\left(\overline{J}_{\tau^{(N)}+\overline{k}+2}^{(N)}= D_N \vert \overline{\tau}^{(N)}<\infty \right) \to 1.
\end{align}
\end{Lem}

\begin{proof}
 We aim to show that all hosts that have not been infected so far, get infected whp in generation $\tau^{(N)}+\overline{k}+2$. 
 According to Lemma \ref{Lem: beta over 2 + delta} we have whp $J^{(N)}_{\tau^{(N)}+\overline{k}+1}\geq  N^{\frac{\beta}{2}+\delta}$. Hence we have whp at least $m_N:=N^{\frac{\beta}{2}+ \delta} V_N$ parasites that may infect the remaining hosts.
So, the probability that an up to generation $\tau^{(N)} +\overline{k}+2$ uninfected host gets attacked by at most one of the $m_N$ parasites (and hence with high probability remains uninfected) can be estimated from above by
 \begin{align*}
&\left(1- \frac{1}{D_N-1}\right)^{m_N}   + \left(1-  \frac{1}{D_N-1 }\right)^{m_N- 1} m_N \frac{1}{D_N-1 }   \\ & \quad \quad =\Theta( N^{\delta} \exp(-a N^{\delta})),
 \end{align*}
because $\frac{m_N}{D_N}=\frac{N^{\beta/2+\delta}V_N}{D_N}\sim \frac{N^{\beta/2}N^{\delta}a\sqrt{D_N}}{D_N}=a N^{\delta}\frac{N^{\beta/2}}{\sqrt{D_N}}=\Theta(N^{\delta}).$\\
The number of uninfected hosts at the beginning of generation $\tau^{(N)}+\overline{k}+2$ is at most $D_N$. Consequently, the probability that at least one of these hosts remains uninfected till the end of generation $\tau^{(N)}+\overline{k}+2$ can be estimated from above by a probability proportional to 
 \begin{align*}
      D_N N^{\delta} \exp(-a N^{\delta})  = o(1),
 \end{align*}
 which yields the claim of Lemma \ref{Lem: every vertices are killed}.
\end{proof}

\subsection{Proof of Theorem \ref{Theorem:Invasion on complete graph} (ii)}
Now we have all necessary ingredients to prove Theorem \ref{Theorem:Invasion on complete graph} (ii). \\
The first step is to show
\begin{align}\label{lab: upper bound invasion probability complete graph}
    \limsup_{N \to \infty} \mP \left(F_{u}^{(N)}\right) \leq \pi(a).
\end{align}
For a sequence $\left(\ell_N\right)_{N \in \mathbb{N}}$ introduce the event 
\begin{align}
    F^{(N)}_{\ell_N}:= \left\{\exists g \in \mathbb{N}_{0}: \overline{J}^{(N)}_{g} \geq \ell_N\right\}. 
\end{align}
Then it follows that for all $0<u\leq 1$ and any sequence $\ell_N \leq u D_N$ we have 
\begin{align} \label{lab: global infection means local infection Complete Graph 2}
    \mathbb{P} \left(F_{u}^{(N)}\right) \leq \mP \left(F^{(N)}_{\ell_N}\right).
\end{align}

Taking a sequence $\ell_N$ satisfying $\ell_N \to \infty$ and $\frac{\ell_N^4 V_{N}^{3}}{D_{N}^{2}}\in o(1)$ we have by Proposition \ref{Prop: Coupling with upper DBPC} that 
\begin{align} \label{lab: coupling upper DBPC}
    \mathbb{P} \left( F^{(N)}_{\ell_N}\right) \leq \mP \left(\exists g \in \mathbb{N}_{0}: \overline{Z}^{(N)}_{g,u} \geq \ell_N\right)+o(1).
\end{align}
Proposition \ref{Proposition: asymptotic probability reaching certain level for BPI} gives that 
\begin{align} \label{lab: control survival probability upper DBPC Complete Graph}
    \lim_{N \to \infty} \mP \left(\exists g \in \mathbb{N}_{0}: \overline{Z}^{(N)}_{g,u} \geq \ell_N\right)=\pi(a).
\end{align}
In summary combining \eqref{lab: global infection means local infection Complete Graph 2}, \eqref{lab: coupling upper DBPC} and \eqref{lab: control survival probability upper DBPC Complete Graph} gives exactly \eqref{lab: upper bound invasion probability complete graph}. \\ 

The second step is to show
\begin{align}\label{lab: lower bound invasion probability complete graph}
    \liminf_{N \to \infty} \mP \left(F_{u}^{(N)}\right) \geq \pi\left(\frac{a}{\sqrt{2^n}}\right).
\end{align}
Proposition \ref{Coupling Lower DBPC} combined with Proposition \ref{Asymptotic survival probability lower DBPC} gives that 
\begin{align} \label{lab: proba initial phase complete graph}
    \liminf_{N \to \infty}\mathbb{P} \left(\exists g \in \mathbb{N}, \overline{I}^{(N)}_{g} \geq N^{\varepsilon}\right) \geq \pi\left(\frac{a}{\sqrt{2^{n}}}\right),
\end{align}
for $\varepsilon>0$ small enough. 
Then Lemma \ref{Lem: every vertices are killed} yields that conditioned on the event $\left\{\exists g \in \mathbb{N}, \overline{I}^{(N)}_{g} \geq N^{\varepsilon}\right\}$ whp all the vertices on the graph finally get infected. Combined with \eqref{lab: proba initial phase complete graph} the claim of \eqref{lab: lower bound invasion probability complete graph} follows.

\subsection{Proof of Theorem \ref{Theorem:Invasion on complete graph}(i)}
\label{lab: proof subcritical case complete graph}

In this subsection we prove Theorem \ref{Theorem:Invasion on complete graph}(i). Recall that in this case $V_N\in o(\sqrt{D_N})$.

    We initially start with one individual, i.e.\ $I_0^{(N)}=1$. We determine the probability that the parasite population gets extinct after one generation. For that we consider the following experiment, where we distribute uniformly at random $V_N$ balls into $D_N-1$ boxes. The probability of extinction after one generation is the same as the probability of the event that all boxes contain at most one ball. Thus, we get that 
    \begin{equation*}
        \mathbb{P} \big(J^{(N)}_{1} = 0 \big)= \frac{(D_N-1)!}{(D_N-1-V_N)!(D_N-1)^{V_N}}\geq \exp\left(-\frac{V_N^{2}}{2(D_N-1)}\right),
    \end{equation*}
    where the inequality was proven in the proof of Proposition~\ref{Prop: Coupling with upper DBPC}. We assumed that $V_N \in o(\sqrt{D_N})$ which implies that $\frac{V_N^2}{2 (D_N-1)}\to 0$ as $N\to \infty$, and thus the right hand side converges to $1$. On the other hand for any $u\in(0,1]$ and $N$ large enough an obvious upper bound for the invasion probability is $\mathbb{P} \left(F_{u}^{(N)}\right)\leq 1-\mathbb{P} \left(J^{(N)}_{1} = 0 \right)$. This implies that 
    \begin{equation*}
        \lim_{N \to \infty} \mathbb{P} \left( F_{u}^{(N)}\right)\leq 1-\lim_{N\to \infty}\mathbb{P} \left(J^{(N)}_{1} = 0 \right)=0.
    \end{equation*}

\subsection{Proof of Theorem 2.7 (iii)}\label{Section: Proof of Theorem 2.7 (iii)}
In this subsection we are going to prove Theorem \ref{Theorem:Invasion on complete graph} (iii). In this case $\sqrt{D_N} \in o(V_N)$. 
The proof is based on using a coupling from below of the total size of the infection process with the total size of a Galton-Watson process whose offspring distribution is close to a $\text{Poi}\left(\frac{a^2}{2}\right)$ distribution until  a level $N^{\alpha}$, with $0<\alpha<\beta$ is reached or until the process dies out. This coupling is possible  for any $a>0$ which yields that the total size of the infection process reaches the level $N^{\alpha}$ with asymptotically probability 1. Then  by choosing $\alpha>\beta/2$ one shows that there exists a generation in which there are at least $N^{\widetilde{\alpha}}$ infected individuals, for some $\widetilde{\alpha}>\beta/2$. In the subsequent generation, all remaining hosts are infected, in the same manner as in Subsection \ref{Subsection: Final phase Epidemic complete graph}.  

We will show 
\begin{align}
    \lim_{N \to \infty} \mathbb{P} \left( \exists g \in \mathbb{N}_{0}: \overline{J}^{(N)}_{g}=D_N\right)=1, 
\end{align}
which shows the claim of Theorem \ref{Theorem:Invasion on complete graph}(iii).

The first step is to couple $\left(\mathcal{S}^{(N)},\mathcal{J}^{(N)},\mathcal{R}^{(N)}\right)$  to an infection process $\big(\widetilde{\mathcal{S}}^{(N)},\widetilde{\mathcal{J}}^{(N)},\widetilde{\mathcal{R}}^{(N)}\big)$, in which infections are only generated by pairs of parasites originating from the same vertex, but not if a host gets infected only by parasites stemming from different vertices. We will show  that \begin{equation} \label{lab: coupling only CoSame}
    \widetilde{\mathcal{J}}_g^{(N)}\cup \widetilde{\mathcal{R}}_g^{(N)}\subset \mathcal{J}_g^{(N)}\cup \mathcal{R}_g^{(N)}, \forall g \in \mathbb{N}_{0}.
\end{equation}

For every vertex $x$ we only need to determine once to which neighbours the $V_N$ offspring parasites move, since afterwards the vertex cannot be used anymore. We denote by $\mathcal{H}_x^{(N)}\subset \{1,\dots, D_N\}\backslash\{x\}$ the set of all vertices which are occupied by at least two or more of the $V_N$ offspring parasites generated on $x$ after their movement. 
With this we can build the coupling of the two processes step by step. We consider for both processes the initial configuration where only vertex $1$ is currently infected and all other vertices are susceptible, i.e.\
\begin{equation*}
    \big(\mathcal{S}^{(N)}_0,\mathcal{I}^{(N)}_0,\mathcal{R}^{(N)}_0\big)=\big(\widetilde{\mathcal{S}}_0^{(N)},\widetilde{\mathcal{I}}_0^{(N)},\widetilde{\mathcal{R}}_0^{(N)}\big)= (\{2,\dots,D_N\},\{1\},\emptyset).
\end{equation*}
Then assume that we constructed the process until generation $g\geq 0$. Then from $g$ to $g+1$ the dynamics are as follows
\begin{align*}
    \widetilde{\mathcal{J}}_{g+1^{(N)}}=\bigcup_{x\in \widetilde{\mathcal{J}}_g^{(N)}}\mathcal{H}_x^{(N)}\backslash(\widetilde{\mathcal{I}}_g^{(N)}\cup \widetilde{\mathcal{R}}_g^{(N)}),\quad \widetilde{\mathcal{S}}_{g+1}^{(N)}=\widetilde{\mathcal{S}}_{g}^{(N)}\backslash \widetilde{\mathcal{J}}_{g+1}^{(N)} \quad \text{ and }\quad  \widetilde{\mathcal{R}}_{g+1}^{(N)}=\widetilde{\mathcal{R}}_{g}^{(N)}\cup \widetilde{\mathcal{J}}_{g}^{(N)}.
\end{align*}
In words every vertex $y\in\mathcal{H}_x^{(N)}$ which is attacked by at least two parasites that are originating from a single vertex $x\in \widetilde{\mathcal{J}}_g^{(N)}$ is added to $\widetilde{\mathcal{J}}_{g+1}^{(N)}$, except for vertices which were already attacked at a previous generation, i.e.\ $y\in\widetilde{\mathcal{J}}_g^{(N)}\cup \widetilde{\mathcal{R}}_g^{(N)}$. Furthermore, all previously infected hosts $\widetilde{\mathcal{J}}_g^{(N)}$ are declared as removed and all vertices which were infected in this generation $\widetilde{\mathcal{J}}_{g+1}^{(N)}$ are removed from the set of susceptible vertices. 

In the process $\left(\mathcal{S}^{(N)},\mathcal{J}^{(N)},\mathcal{R}^{(N)}\right)$ cooperation from different infected vertices for the spread of the epidemic is allowed. Since we defined movement parasites independent from the generation at which vertices get infected, we have
\begin{equation} \label{lab: coupling only CoSame}
    \widetilde{\mathcal{J}}_g^{(N)}\cup \widetilde{\mathcal{R}}_g^{(N)}\subset \mathcal{J}_g^{(N)}\cup \mathcal{R}_g^{(N)}, \forall g \in \mathbb{N}_{0},
\end{equation}
almost surely. 
As by cooperation only more infections are generated, it is not possible that a vertex $x$ which is susceptible for both processes at a generation $g$ gets infected at generation $g+1$ for the process $\big(\widetilde{\mathcal{S}}^{(N)},\widetilde{\mathcal{J}}^{(N)},\widetilde{\mathcal{R}}^{(N)}\big)$ but not for the process $\big(\mathcal{S}^{(N)},\mathcal{J}^{(N)},\mathcal{R}^{(N)}\big)$.

The infection process $\big(\widetilde{\mathcal{S}}^{(N)},\widetilde{\mathcal{J}}^{(N)},\widetilde{\mathcal{R}}^{(N)}\big)$ is monotone with respect to the parameter $V_N$, in contrast to the original process $\big(\mathcal{S}^{(N)},\mathcal{J}^{(N)},\mathcal{R}^{(N)}\big)$. Now let $a>0$ and consider $V_N^{(a)}=a\sqrt{D_N}$ as well as $\big(\widetilde{\mathcal{S}}^{(N,a)},\widetilde{\mathcal{J}}^{(N,a)},\widetilde{\mathcal{R}}^{(N,a)}\big)$ to be the analogously defined infection process. Infections are only generated by pairs parasites originating from the same vertex as well as the number of parasites generated at an infection event is $V_N^{(a)}$. Since we assume that $D_N\in o(V_N)$ it follows for $N$ large enough that $V_N^{(a)}\leq V_N$. Thus, by monotonicity it follows that we can couple $\big(\widetilde{\mathcal{S}}^{(N,a)},\widetilde{\mathcal{J}}^{(N,a)},\widetilde{\mathcal{R}}^{(N,a)}\big)$ and $\big(\widetilde{\mathcal{S}}^{(N)},\widetilde{\mathcal{J}}^{(N)},\widetilde{\mathcal{R}}^{(N)}\big)$, such that 
\begin{align}\label{Couple}
|\widetilde{\mathcal{R}}^{(N,a)}_g \cup \widetilde{\mathcal{J}}^{(N,a)}_g | \leq  |\widetilde{\mathcal{R}}^{(N)}_g \cup \widetilde{\mathcal{J}}^{(N)}_g| .\end{align}

For the sequence of processes $\big(\widetilde{\mathcal{S}}^{(N,a)},\widetilde{\mathcal{J}}^{(N,a)},\widetilde{\mathcal{R}}^{(N,a)}\big)$ we can show (by a coupling with Galton-Watson processes) that the probability to infect eventually $N^\alpha$ host is asymptotically lower bounded by the survival probability $\varphi_a$ of a Galton-Watson process with Poi$\left(\tfrac{a^2}{2}\right)$ offspring distribution. The proof of this statement can be found in the proof of Lemma \ref{Lemma: InvasionTildeProcess}, where this statement is formulated, in the supplementary material (since it can be shown by very similar arguments that have been used to show Proposition 4.7 in \cite{BrouardEtAl2022}). 

Because this result is true for any $a>0$, taking the limit when $a$ goes to $\infty$ gives, together with \eqref{Couple} and  \eqref{lab: coupling only CoSame}
\begin{equation}\label{lab: reaching high level proba 1}
    \lim_{N \to \infty} \mP \left( \exists g \in \mathbb{N}_{0}:  \overline{J}^{(N)}_{g} \geq N^{\alpha}\right)=1.
\end{equation}
Now let $\frac{\beta}{2}<\alpha<\beta$ and introduce 
\begin{align}
    &\overline{\tau}^{(N)}_{N^{\alpha}}:= \inf \left\{ g \in \mathbb{N}_{0}: \overline{J}^{(N)}_{g} \geq N^{\alpha}\right\}, \\
    &\tau^{(N)}_{N^{\alpha}}:=\inf \left\{g \in \mathbb{N}_{0}:  J^{(N)}_{g}  \geq \frac{N^{\alpha}}{\log^{2}(N)}\right\}.
\end{align}
Then one can show as in Proposition \ref{Prop: reached certain level if total size reaches one} that 
\begin{equation}
    \mP \left( \tau^{(N)}_{N^{\alpha}} \leq \overline{\tau}^{(N)}_{N^{\alpha}} \vert \overline{\tau}^{(N)}_{N^{\alpha}}<\infty\right) \to 1.
\end{equation}
Indeed, if for any generation before $\overline{\tau}^{(N)}_{N^{\alpha}}$ the number of infected vertices is strictly smaller than $\frac{N^{\alpha}}{\log^{2}(N)}$ then this would mean that the number of generations to reach the level $N^{\alpha}$ for the total size of the infection process is at least $\log^{2}(N)$. But this is in contradiction with the couplings of \eqref{lab: coupling only CoSame} and \eqref{lab: coupling below GWP} and Lemma 5.5 from \cite{BrouardEtAl2022}. \\ 
Then using a similar approach as in the proof of Lemma \ref{Lem: every vertices are killed}, one shows that 
\begin{align}\label{lab: everyone is killed}
    \mP \left( \overline{J}^{(N)}_{\tau^{(N)}_{N^{\alpha}}+1}  =D_N \vert \overline{\tau}^{(N)}_{N^{\alpha}}<\infty\right) \to 1.
\end{align}
Finally combining \eqref{lab: reaching high level proba 1} and \eqref{lab: everyone is killed} it follows that 
\begin{align}
    \mP \left( \exists g \in  \mathbb{N}_{0}: \overline{J}^{(N)}_{g}=D_N\right) \to 1,
\end{align}
which completes the proof.

\section{Invasion on a random geometric graph on $[0,1]^{n}$}\label{Sec:RGG}

To start with we show  that $\mathcal{G}^{(N)}$ is with high probability fully connected and is fairly dense in the sense that the number of vertices contained in every ball of radius $r_N$ is of order $N^{\beta}$.

\begin{Lem}\label{Remark: uniform bound number vertices}
    \begin{enumerate}
        \item The graph $\mathcal{G}^{(N)}=(\mathcal{V}^{(N)},\mathcal{E}^{(N)})$ is fully connected with high probability as $N\to \infty$.
        \item Let $\frac{2}{n+2}\beta<\gamma<\beta$. 
        Then, it holds that 
        \begin{align*}
            \lim_{N\to \infty} &\mP\big(|\mathcal{V}^{(N)}\cap B_{r_N}(x)|\geq N^{\beta}-(n+1)N^{\frac{(n-1)\beta+\gamma}{n}} \,\forall\, x\in \big[\tfrac{1}{2} N^{\beta-1},1-\tfrac{1}{2} N^{\beta-1}\big]^n)=1\text{ and}\\
            \lim_{N\to \infty} &\mP(|\mathcal{V}^{(N)}\cap B_{r_N}(x)|\leq N^{\beta}+(2n+1)N^{\frac{(n-1)\beta+\gamma}{n}}\,\forall\, x\in \big[\tfrac{1}{2} N^{\beta-1},1-\tfrac{1}{2} N^{\beta-1}\big]^n)=1.
        \end{align*}
   \end{enumerate}
\end{Lem}

\begin{proof}[Proof of Lemma \ref{Remark: uniform bound number vertices}]
  Choose $0<\gamma<\beta$ and $0<\varepsilon<\gamma/2$. The idea of the proof is to define disjoint boxes $B(l)$ with $l=(l_1,\dots,l_n)\in \N^n$ with side length $N^{\frac{\gamma-1}{n}}$ which cover the whole unit box, i.e.\ $[0,1]^n\subset \bigcup B(l)$. In the second step we  gain control on the asymptotic number of Poisson points contained in every box simultaneously, i.e.\ we will show with the help of Lemma~\ref{Lemma: moderate deviation number of vertices} that every box contains $N^{\gamma}\pm N^{\gamma/2+\varepsilon}$ many points with high probability. A technical problem is  that we defined our Poisson point set $\mathcal{V}^{(N)}$ only on $[0,1]^n$. Not for every $N$ are we able to perfectly cover the unit box with our boxes $B(l)$ such that $[0,1]^n= \bigcup B(l)$. Thus, we need to extend our Poisson point set. This can be easily done by sampling independent Poisson points with intensity measure $N\mathsf{d}x$ on $[0,2]^n \backslash [0,1]^n$. We denote this Poisson point set by $\mathcal{V}_N'$. Now we set $\mathcal{V}''_N:=\mathcal{V}^{(N)}\cup \mathcal{V}'_N$, so $\mathcal{V}''_N$ is a Poisson point set on $[0,2]^n$ with intensity measure $N \mathsf{d}x$.
  
  Let us set $\overline{M}:=\{0,\dots, \lceil N^{\frac{1-\gamma}{n}}\rceil\}^n$ and $\underline{M}:=\{0,\dots, \lfloor N^{\frac{1-\gamma}{n}}\rfloor\}^n$ . Define boxes of side length $N^{\frac{\gamma-1}{n}}$ by setting $B(l):=[0,N^{\frac{\gamma-1}{n}})^n+l$, where $l\in \overline{M}$. Set $\overline{k}:=|\overline{M}|$ and $\underline{k}:=|\underline{M}|$. For these boxes we have
  \begin{equation*}
      \bigcup_{l\in \underline{M}}B(l)\subset [0,1]^n\subset \bigcup_{l\in \overline{M}}B(l)
  \end{equation*}
  Set $X_N:=|B(\mathbf{0})\cap \mathcal{V}''_N|$, where $\mathbf{0}=(0,\dots,0)\in\R^n$ then
  \begin{align*}
   \mP\Big(\bigcap_{l\in \overline{M}}\{N^{\gamma}-N^{\frac{\gamma}{2}+\varepsilon} \leq |B(l)\cap \mathcal{V}''_N|\leq  N^{\gamma}+N^{\frac{\gamma}{2}+\varepsilon}\}\Big)
      =\mP\big(N^{-\tfrac{\gamma}{2}-\varepsilon}|X_N-N^{\gamma}|\leq    1\big)^{\overline{k}}.
  \end{align*}
  According to  Lemma~\ref{Lemma: moderate deviation number of vertices}, where we control the size of Poisson random variables via moderate deviations,
  \begin{align*}
    \lim_{n\to \infty} -\frac{2}{N^{2\varepsilon}}\log\big(\mP\big(N^{-\tfrac{\gamma}{2}-\varepsilon}|X_N-N^{\gamma}|>1\big)\big)= 1.
   \end{align*}
   This implies that
   \begin{align*}
  \log\big(\mP\big(N^{-\tfrac{\gamma}{2}-\varepsilon}|X_N-N^{\gamma}|>1\big)\big)=  -\frac{N^{2\varepsilon}}{2}(1+h(N^{2\varepsilon}))
   \end{align*}
   where $h(x)\in o(1)$ as $x\to \infty$. Since $\overline{k}=\lceil N^{\frac{1-\gamma}{n}}\rceil^n$ with Bernoulli's inequality 
    \begin{align*}
        \big(1-\mP\big(N^{-\tfrac{\gamma}{2}-\varepsilon}|X_N-N^{\gamma}|>   1\big)\big)^{\overline{k}}\geq 1-\lceil N^{\frac{1-\gamma}{n}}\rceil^n\exp\Big(-\frac{N^{2\varepsilon}}{2}(1+h(N^{2\varepsilon}))\Big) \rightarrow 1
    \end{align*}
    as $N\to \infty$. Thus, we have shown that all boxes $(B(l))_{l\in \overline{M}}$ simultaneously contain with high probability $N^{\gamma} \pm N^{\frac{\gamma}{2}+\varepsilon}$ many Poisson points as $N\to \infty$.
    \begin{enumerate}
        \item The first claim is a direct consequence of what we just showed. Let $l\in \underline{M}$, i.e. we consider a box $B(l)\subset [0,1]^n$, then  it follows that every vertex $x\in \mathcal{V}^{(N)}$ contained in $B(l)$ is connected to every other vertex contained in the same box $B(l)$ since $\gamma<\beta$. This means that the vertices in a box $B(l)$ form a complete graph for every $l\in \underline{M}$. Furthermore, for $N$ large enough it holds that $2N^{\frac{\gamma-1}{n}}<N^{\frac{\beta-1}{n}}$, and thus every vertex contained in a box $B(l)$ is connected to every vertex contained in all adjacent boxes $B(l')$.
        Thus, we have shown that the random geometric graph with vertex set $\mathcal{V}^{(N)}\cap \bigcup_{l\in \underline{M}}B(l)$ forms a connected graph with high probability. It only remains to argue that every vertex $x\in \mathcal{V}^{(N)}\cap  [0,1]^n\backslash \bigcup_{l\in \underline{M}}B(l)$ is connected to its neighbouring box.  Note that a vertex  $x\in \mathcal{V}^{(N)}\cap  [0,1]^n\backslash \bigcup_{l\in \underline{M}}B(l)$ it holds $B_{r_N}(x)\cap \bigcup_{l\in M}B(l)\neq \emptyset$ since $\gamma<\beta$. Hence, for $N$ large enough  these vertices are connected to its closest box $B(l)$ with high probability, since with high probability every box $B(l)$ is non-empty for $N\to \infty$.
        \item Every ball $B_{r_N}(x)\subset[0,1]$ contains $\lfloor N^{\frac{\beta-\gamma}{n}}\rfloor^n$ many boxes of length $N^{\frac{1-\gamma}{n}}$. This means that with high probability every ball $B_{r_N}(x)\subset[0,1]$ contains at least $\lfloor N^{\frac{\beta-\gamma}{n}}\rfloor^n(N^{\gamma}-N^{\frac{\gamma}{2}+\varepsilon})$ many vertices. Note that
        \begin{align*}
            \lfloor N^{\frac{\beta-\gamma}{n}}\rfloor^n(N^{\gamma}-N^{\frac{\gamma}{2}+\varepsilon})\geq& (N^{\frac{\beta-\gamma}{n}}-1)^n(N^{\gamma}-N^{\frac{\gamma}{2}+\varepsilon})
            =N^{\beta}-nN^{\frac{(n-1)\beta+\gamma}{n}}+R_{-}(N),
        \end{align*}
        where $R_-(N)= -N^{\beta-\frac{\gamma}{2}+\varepsilon}+o(N^{\beta-\frac{\gamma}{2}+\varepsilon})$. Since
        \begin{equation*}
            \beta-\frac{\gamma}{2}<\frac{(n-1)\beta+\gamma}{n} \Leftrightarrow 
            \frac{2}{(n+2)}\beta\leq \gamma
        \end{equation*}
        for all $n\geq 1$, and thus we can choose $\varepsilon$ small enough such that $R_-(N)$ consists only of lower order terms with the leading order term having a negative sign. This means that for $N$ large enough it follows that
        \begin{equation*}
            \lfloor N^{\frac{\beta-\gamma}{n}}\rfloor^n(N^{\gamma}-N^{\frac{\gamma}{2}+\varepsilon})\geq N^{\beta}-(n+1)N^{\beta-\frac{\gamma}{2}+\varepsilon}.
        \end{equation*}
        If $B_{r_N}(x)\subset [0,1]$, then it is covered by $(\lfloor N^{\frac{\beta-\gamma}{n}}\rfloor+2)^n$ many boxes. Thus, we obtain similar as before that 
        \begin{equation*}
            (\lfloor N^{\frac{\beta-\gamma}{n}}\rfloor+2)^n(N^{\gamma}-N^{\frac{\gamma}{2}+\varepsilon})\leq N^{\beta}+2nN^{\frac{(n-1)\beta+\gamma}{n}}+R_{+}(N),
        \end{equation*}
        where $R_{+}(N)=N^{\beta-\frac{\gamma}{2}+\varepsilon}+o(N^{\beta-\frac{\gamma}{2}+\varepsilon})$ such that we again for $N$ large enough we get that
        \begin{equation*}
            (\lfloor N^{\frac{\beta-\gamma}{n}}\rfloor+2)^n(N^{\gamma}-N^{\frac{\gamma}{2}+\varepsilon})\leq N^{\beta}+(2n+1)N^{\frac{(n-1)\beta+\gamma}{n}}.
        \end{equation*}
    \end{enumerate}
\end{proof}
\begin{remark}\label{Remark: Choice of gamma}
The optimal choice of $\gamma$ to minimize the order of the error term is to choose $\gamma$ close to $\frac{2}{2+n}\beta$, which leads to an order close to
\begin{equation*}
    \frac{(n-1)\beta+\gamma}{n}
    =\Big(\frac{n+1}{n+2}\Big)\beta.
\end{equation*}
But the result of Lemma~\ref{Remark: uniform bound number vertices} does not allow for this choice. Thus, one reasonable choice would for example be $\gamma=\frac{3}{3+n}\beta$, then we get that the order of the error term is 
\begin{equation*}
    \frac{(n-1)\beta+\gamma}{n}
    =\Big(\frac{n+2}{n+3}\Big)\beta,
\end{equation*}
which yields for $n=1$ the value $\frac{3}{4}\beta$.
\end{remark}

\subsection{Upper bound  on the invasion probability} 

To derive an upper bound on the invasion probability we couple whp the total number of infected hosts from above  with the total size of a DBPC whose laws are approximately Poisson distributed until the DBPC dies out or reaches at least the level $\ell_N$. 

Let $\delta_{N,\ell}: = N^\beta -  (n+1)N^{\frac{n+2}{n+3}\beta}$ and $\delta_{N,u}: = N^\beta +  (2n+1)N^{\frac{n+2}{n+3}\beta}$ . According to Lemma \ref{Remark: uniform bound number vertices} and Remark \ref{Remark: Choice of gamma} whp every ball $B_{r_N}$ contains at least $\delta_{N, \ell}$ and at most $\delta_{N,u}$ vertices $x\in \mathcal{V}^{(N)}$. 

\begin{Def}(Upper DBPC) \\
\label{Definition: UDBPC}
Let $\ell_{N} \underset{N \to \infty}{\to}\infty$ satisfying $\ell_N\in o(\log \log N)$. Let $\ZZ^{(N)}_{u}=\left(Z_{g,u}^{(N)}\right)_{g \in \mathbb{N}_0}$ be a branching process with interaction with $Z_{0,u}^{(N)}=1$ almost surely, and offspring and cooperation distribution with probability weights $p_{u,o}^{(N)} = \left( p_{j,u, o}^{(N)}\right)_{j \in \mathbb{N}_{0}}$  and  $p_{u,c}=\left( p_{j,u,c}^{(N)}\right)_{j \in \mathbb{N}_{0}}$, respectively  with 
\begin{align}
    p_{j,u, o}^{(N)}:= \left( \frac{(v_N-\ell_N^2)^{2}}{2 \delta_{N,u}}\left( 1 -\frac{  2^{\ell_N} \ell_N^2}{\delta_{N,\ell}}\right) \right)^{j} \frac{1}{j!} \exp \left( -\frac{(v_{N}- \ell_N^2)^{2}}{2 \delta_{N,u}} \left( 1 -\frac{  2^{\ell_N} \ell_N^2}{\delta_{N,\ell}}\right) \right),
\end{align}
for all $0 \leq j \leq \ell_N$ and 
\begin{equation}
    p_{\ell_N+1,u,o}^{(N)}:=1-\sum_{j=0}^{\ell_N} p_{j,u,o}^{(N)},
\end{equation}
as well as 
\begin{align}
    p_{j,u, c}^{(N)}:= \left( \frac{(v_N-\ell_N^2)^{2}}{\delta_{N,u}} \left( 1 -\frac{  2^{\ell_N} \ell_N^2}{\delta_{N,\ell}}\right) \right)^{j} \frac{1}{j!} \exp \left( -\frac{(v_{N}- \ell_N^2)^{2}}{\delta_{N,u}} \left( 1 -\frac{  2^{\ell_N} \ell_N^2}{\delta_{N,\ell}}\right) \right) ,
\end{align}
for all $0 \leq j \leq \ell_N$ and 
\begin{equation}
    p_{\ell_N+1,u,c}^{(N)}:=1-\sum_{j=0}^{\ell_N} p_{j,u,c}^{(N)}.
\end{equation}
Denote by $\overline{Z}_{u}^{(N)}:= \left( \overline{Z}_{g,u}^{(N)}\right)_{g \in \mathbb{N}_{0}}$ where $\overline{Z}_{g,u}^{(N)}:=\sum_{i=0}^{g}Z_{i,u}^{(N)}$, that is $\overline{Z}_{g,u}^{(N)}$ gives the total size of $\ZZ^{(N)}_{u}$ accumulated till generation $g$.
\end{Def}

\begin{Prop}\label{SurvivalProb Upper DBPC RGG}(Probability for the total size of the upper DBPC to reach a level $a_N$). \\
Consider a sequence $\left(a_N\right)_{N \in \mathbb{N}}$ with $a_N \underset{N \to \infty}{\to}\infty$ and assume that $v_N \sim a \sqrt{d_N}$ for $0<a<\infty$. Then, we have 
\begin{align}
    \lim_{N \to \infty} \mathbb{P} \left( \exists g \in \mathbb{N}_{0}: \overline{Z}^{(N)}_{g,u} \geq a_{N}\right)=\pi(a). 
\end{align}
\end{Prop}
 \begin{proof}
This proposition is proven analogously as Proposition~\ref{Proposition: asymptotic probability reaching certain level for BPI}, i.e.\ one can show  that the assumptions of Lemma~\ref{Lemma: convergence survival probability DBPC} are satisfied and then the claim is a consequence of this lemma.
 \end{proof}

\begin{Prop}
Consider a sequence $\left(\ell_N \right)_{N \in \mathbb{N}}$ with $\ell_N \underset{N \to \infty}{\to}\infty$ satisfying $\ell_N\in o(\log \log (N))$. Introduce the stopping time
\begin{align}
    \tau^{(N)}_{\ell_N,0}:= \inf \left\{ g \in \mathbb{N}_{0}: \overline{Z}^{(N)}_{g,u} \geq \ell_{N} \text{ or } Z_{g,u}^{(N)}=0\right\}.
\end{align}
Then 
\begin{equation}
    \lim_{N \to \infty} \mathbb{P} \left( \overline{I}^{(N)}_{g} \leq \overline{Z}^{(N)}_{g,u}, \forall g < \tau^{(N)}_{\ell_N, 0}\right)=1 
\end{equation}
and 
\begin{equation}
    \lim_{N \to \infty} \mathbb{P} \left( I^{(N)}_{\tau^{(N)}_{\ell_N, 0}} =0 \Big\vert  Z^{(N)}_{\tau^{(N)}_{\ell_N, 0}, u} =0\right)=1.
\end{equation}
\end{Prop}\label{UpperBoundRGG}

\begin{proof}
For the proof couple the infection process with another model, that uses the same infection rules but assumes that every generation empty vertices are reoccupied by an host. This increases only the number of infections. Denote by $\tilde{I}^{(N)}=(\tilde{I}_g^{(N)})_{g \in \mathbb{N}_{0}}$ the corresponding process that counts the number of infections generated in this modified model. We have $ I_g^{(N)} \leq \tilde{I}_g^{(N)}$ a.s.
Next we show that $\tilde{I}_g^{(N)} \leq Z_{g,u}^{(N)}$ whp for all $g < \tau^{(N)}_{\ell,0}$. We say that in generation $g$ we have $k$ $\tilde{I}^{(N)}$ infections, if $\tilde{I}_g^{(N)}=k$ and we say that in generation $g$ we have  $k$ $\ZZ_u^{(N)}$ infections, if $Z_{g,u}^{(N)}=k$.
Start with generation $g=0$. 
Since initially only a single vertex is infected, in the first generation only CoSame infections are possible. As in \cite{BrouardEtAl2022} we can couple $\tilde{I}_1^{(N)}$ with $Z_{1,u}^{(N)}$, such that $\mP\left(\tilde{I}_1^{(N)} \leq Z_{1,u}^{(N)}\right)\geq 1- o(\frac{1}{\log N})$ for any $\eps>0$, see Proposition 3.2 in \cite{BrouardEtAl2022}. 
Next we proceed iteratively. Assume in generation $g$ $m=m_N$ vertices are $\tilde{I}^{(N)}$-infected. 
If $m=1$ we can use the coupling  as in generation 0 and add independently additional CoSame and CoDiff infections according to the distribution DBPC in $\ZZ_{u}^{(N)}$, if $Z_{1,u}^{(N)}>1$. 

If $m >1$, let $v_1, ..., v_{m}$ be the infected vertices and denote by $\overline{\mathcal{D}}_i$ the set and by $\overline{D}_i$ the number of vertices in the ball of radius $r_N$ around vertex $v_i$ for $i=1, ..., m$. For $y\in \{0,1\}^{m}$ denote by $\mathcal{D}_{y}$ the set and by $D_y$ the number of vertices that are contained in all balls that are centered around vertices $v_i$, $i\in 1, ..., m$ which have a 1 at the $i$-th position of the vector $y$ and are not contained in the other balls. So for example for $m=3$ $D_{001}$ gives the number of vertices that are contained only in the ball around vertex 3, but not in the balls centered around vertex 1 and 2. 

For a vector $x=(x_1, ..., x_m, x_{1,2}, x_{1,3}, ..., x_{m-1,m}) \in \N^{m + \binom{m}{2}}$ denote by $E_x$ the event that (i) in the next generation $x_i$ CoSame infections occur caused by exactly two parasites generated on vertex $v_i$ for $i=1, ..., m$, $x_{i,j}$ Codiff infections occur caused by exactly two parasites being generated on vertex $v_i$ and vertex $v_j$ for $i,j\in \{1, ..., m\}$ with $i<j$ and all other vertices get attacked by at most one parasite.

To determine the probability of the event $E_x$ we distinguish different cases. Let for $y \in \{0,1\}^m$ denote by $x^o_{y,i}$ the number of CoSame infections caused by parasites generated on vertex $i$ attacking vertices in $\mathcal{D}_y$ as well as by $x^c_{y,i,j}$ the number of CoDiff infections generated by parasites from vertex $i$ and $j$ that are attacking vertices in $\mathcal{D}_y$ as well as $x^r_{i,j}$ the number of parasites originating from vertex $j$ and attacking a vertex without any other parasite in $\overline{\mathcal{D}}_i$. The probability of $E_x$
is given by the sum of the probabilities of infection patterns corresponding to vectors $(x^o_{y,i})_y$, $(x^c_{y,i,j})_y$ and $(x_{r,i,j})_{i>j}$
with $\sum_y x^o_{y,i}= x_i$, $\sum_y x^c_{y,i,j}= x_{i,j}$ and $x^c_{y,i,j}=0$, if the $i$th or $j$th coordinate of $y$ is 0, such that $D_y- \underline{x}_{y} >0$ with  $\underline{x}_{y} = \sum_i x^o_{y,i} + \sum_{i=1}^{m-1}\sum_{j=i+1}^{m} x^c_{y,i,j}$ for all $y\in \{0,1\}^m$. The probability of an infection pattern according to the vectors 
$x^o=(x^o_{y,i})_{y,i}$ and $x^c=(x^c_{y,i,j})_{y,i,j}$ is given by the product of the three factors $p_o(x^o,x^c)$, $p_c(x^o, x^c)$ and $p_r(x^o, x^c)$ representing the CoSame, CoDiff and single infections with
\begin{align*}
p_o= \prod_{i=1}^m \binom{v_N}{2} \cdots \binom{v_N -2\left( x_i-1\right)}{2}  \frac{1}{D_i^{2 x_i}}
\prod_y  \frac{(D_y- \overline{x}_{y,i})!} {(D_y- \overline{x}_{y,i}- x_{y,i}^{o})!}\frac{1}{x_{y,i}^{o}!} 
\end{align*}

where $\overline{x}_{y,i} = \sum_{j=1}^{i-1}  x_{y,j}^{o} $. The factor $\binom{v_N}{2} \cdots \binom{v_N -2\left( x_i-1\right)}{2}$ gives the number of possibilities to choose $x_i$ pairs of parasites from the $v_N$ parasites generated on vertex $i$, for $i=1,..., m$,  $\frac{(D_y- \overline{x}_{y,i})!} {(D_y- \overline{x}_{y,i}- x_{y,i}^{o})!}\frac{1}{x_{y,i}^{o}!}$ gives the number of possibilities to choose for $x^o_{y,i}$ pairs of parasites a location in $D_y$, when we already distributed the pairs of parasites generated on vertices $j=1, ..., i-1$ on $D_y$. $\frac{1}{D_i^{2x_i}}$ is the probability to place the pairs of parasites exactly on these locations in $D_y$.

\begin{align*}
p_c=& \prod_{i=1}^{m} \prod_{j=i+1}^m \frac{ (v_N- \underline{x}_{1,i,j} )!}{(v_N - \underline{x}_{1,i,j} - x_{i,j})!} \frac{  (v_N - \underline{x}_{2,i,j})!}{(v_N - \underline{x}_{2,i,j} - x_{i,j})!} \\ &  \cdot  \prod_y \frac{(D_y-\underline{x}_{y,i,j})!}{(D_y - \underline{x}_{y,i,j}- x^c_{y,i,j})!} \left(\frac{1}{D_i}\right)^{x_{y,i,j}}   \left(\frac{1}{D_j}\right)^{x_{y,i,j}} \frac{1}{(x_{y,i,j})!} 
\end{align*}
with $\underline{x}_{1,i,j} = 2x_i  + \sum_{\ell=i+1}^{j-1} x_{i,\ell} $ and 
$\underline{x}_{2,i,j} = 2x_j  + \sum_{\ell=1}^{i-1} x_{\ell,j} $,
$\underline{x}_{y,i,j} = \sum_{k=1}^m x^o_{y,k} + \sum_{k=1}^{i-1}\sum_{\ell=k+1}^{m} x^c_{y,k,\ell} + \sum_{k=i+1}^{j-1} x^c_{y,i,k}$.
The factor $\frac{(v_N- \underline{x}_{1,i,j} )!}{(v_N - \underline{x}_{1,i,j} - x_{i,j})!}$ 
gives the number of possibilities to choose $x_{i,j}$ parasites from the parasites generated on vertex $i,$ when the parasites for the CoSame infections as well as the parasites for the CoDiff infections of the vertex pairs $(i,i+1), ..., (i, j-1)$ have already been determined. The factor $\frac{(D_y-\underline{x}_{y,i,j})!}{(D_y - \underline{x}_{y,i,j}- x^c_{y,i,j})!}  \frac{1}{(x_{y,i,j})!} $ gives the number of possibilities to choose in $D_y$ the $x^c_{y,i,j}$ locations for the pairs of parasites generating a CoDiff infection from vertex $i$ and $j$, when the locations for the CoSame infections as well as for the CoDiff infections of vertex pairs $(1,2),  ..., (i, j-1)$ have already been determined. 
Finally, the factor $\left(\frac{1}{D_i}\right)^{x_{y,i,j}}   \left(\frac{1}{D_j}\right)^{x_{y,i,j}}$ is the probability to place the pairs of parasites generating the CoDiff infections on exactly these locations.

\begin{align*}
p_r = \prod_{i=1}^{m} \frac{D_i -  2 x_i - \sum_{j=i+1}^m x_{i,j} - \sum_{j=1}^{i-1} x^r_{i,j}  }{D_i} ... \frac{D_i  -\sum_{j=1}^{i-1} x^r_{i,j}  -v_N +1}{D_i}. 
\end{align*}
$p_r$ is the probability to place the remaining parasites all onto different vertices.
 
To analyse the above probabilities, consider only those configurations $(x^o_{y,i})_{y,i}, (x^c_{y,i,j})_{y,i,j}$ with positive entries for vectors $y$ for which $1/D_y \in o(1/(d_N)^{1-\varepsilon})$ for some $\varepsilon >0$ (independently of $N$) and only values $x_i, x_{i,j}\leq \ell_N$, because the sum of the remaining probabilities is $O(d_N^{-\varepsilon}).$ Under this assumption we can estimate

\begin{align}\label{po}
p_o \geq \prod_{i=1}^m \prod_y \left(\frac{(v_N- \ell_N^2)^2}{2d_{N,u} } \frac{D_y-\ell_N^2}{D_i} \right)^{x^o_{y,i}}   \frac{1}{x^o_{y,i}!}. 
\end{align}
Then by setting $t_N^y= \frac{(v_N -  \ell_N^2)^2}{2 \delta_{N,u}} \frac{D_y - \ell_N^2}{D_i}$ we can write
\begin{align*}
\eqref{po} & =\prod_{i=1}^m \prod_y \left(\frac{(v_N- \ell_N^2)^2}{2\delta_{N,u} } \frac{D_y-\ell_N^2}{D_i} \right)^{x^o_{y,i}}   \frac{1}{x^o_{y,i}!}
\\ &= \prod_{i=1}^m \prod_y (t_N^y)^{x_{y,i}^o} \frac{1}{x_{y,i}^o!} 
\exp(-t_N^y) \exp(t_N^y) \\
&= \prod_{i=1}^m \left(\prod_y \left(t_N^y\right)^{x_{y,i}^o} \frac{1}{x_{y,i}^o!} \exp(-t_N^y)  \right) 
\left(\prod_y \exp(t_N^y) \right) 
\\ & \geq \prod_{i=1}^m \PP(Y_i= x_i) \exp(a_{N,i}), 
\end{align*}
for $\text{Poi}(a_{N,i})$ distributed random variables $Y_i$ with $a_{N,i}= \frac{(v_N- \ell_N^2)^2}{2 \delta_{N,u}} \left( 1 -\frac{  2^{m} \ell_N^2}{D_i}\right) $. Since $a_{N,i} \geq  \frac{(v_N- \ell_N^2)^2}{2\delta_{N,u}}\left( 1 -\frac{  2^{m} \ell_N^2}{\delta_{N,\ell}}\right)=: a_N $
we have $\eqref{po} \geq \prod_{i=1}^m \PP(Y =x_i) \exp(a_N)$
with $Y \sim \text{Poi}(a_N)$. 

Similarly, we have 

\begin{align*}
p_c & \geq \prod_{i=1}^{m} \prod_{j=i+1}^m  (v_N- \ell_N^{2})^{2x_{i,j}} \prod_y \left(\frac{D_y - \ell_N^2}{\delta_{N,u}}\right)^{x_{y,i,j}} \left(\frac{1}{\delta_{N,u}}\right)^{x_{y,i,j}} \frac{1}{x_{y,i,j}!} \\ &\geq \prod_{i=1}^m \prod_{j=i+1}^m  \PP(Y_2= x_{i,j} ) \exp(2 a_N)
\end{align*} 

with $Y_2 \sim \text{Poi}(2a_{N})$. 

Furthermore  $p_r \geq \exp \left(- m^2 a_N \right)$ for $N$ large enough, since $m\geq 2$.
Consequently, we have 
$$\PP(E_x) \geq \left(\prod_{i=1}^m \PP(Y_1=x_i) \prod_{k=1}^{m} \prod_{\ell=k+1}^{m} \PP(Y_2=x_{k,\ell})\right).$$ 

Furthermore, we have
$$\PP(\cup_{x, x\leq \ell_N} E_x)= 1 - o(1),$$ where we write $x\leq \ell_N$, if $x_i \leq \ell_N$ and $x_{i,j}\leq \ell_N$ for all $i,j\in \{1, ...,m\}$ with $j>i$. 

Since for any $0\leq  k \leq \ell_N$ 
\begin{align*}
&\PP\left(Z_{g,u}^{(N)} =k| Z_{g-1,u}^{(N)} =m\right)  \\ & \quad \quad = \sum_{x \in \mathbbm{N}^{m+ \binom{m}{2}} : \sum_{i=1}^m x_i + \sum_{i=1}^m \sum_{j=i+1}^m x_{i,j} = k  }  
\prod_{i=1}^m \PP(Y_1=x_i) \prod_{k=1}^{m} \prod_{\ell=k+1}^{m} \PP(Y_2=x_{k,\ell}),
\end{align*}
we can couple $(\tilde{I}_g)$ and $(Z_{g,u})_g$ such that $\tilde{I}_g \leq Z_{g,u}$ whp for any $g \leq \tau_{\ell_N, 0}$.
\end{proof}

\subsection{Lower bound on the invasion probability}
\subsubsection{Establishing invasion}
In this section we show that in the random geometric graph  the level $N^{\varepsilon}$ is reached with at least the probability with which a well chosen lower DBPC reaches this level. \\
Consider the ball of radius $\frac{r_N}{2}$ centered in the initial infected vertex. According to  Lemma~\ref{Remark: uniform bound number vertices} and \ref{Remark: Choice of gamma}  this ball contains whp at least  $\frac{N^{\beta}}{2^{n}}- (2n+1)N^{\frac{n+2}{n+3}\beta}$ and at most $\frac{N^{\beta}}{2^{n}}+ (2n+1)N^{\frac{n+2}{n+3}\beta}$ vertices. We are interested in the probability that at least  $N^{\varepsilon}$ get infected for some $0<\varepsilon< \frac{\beta}{4}$.

\begin{Def} (Sub infection process) \\
Let $x_{c}$ be the center point of $[0,1]^n$ and $x_0\in \mathcal{V}^{(N)}$ be the vertex with the smallest distance to $x_c$. The sub-infection process $\mathcal{J}^{(N)}=(\mathcal{J}_g)_{g\geq0}$ is defined on the complete neighbourhood $\mathcal{C}(x_c)$ of $x_c$, introduced in the sketch of the proof of Theorem \ref{MainResult} that is contained in the ball $B_{r_N/2}(x_c)$. We set $\mathcal{J}^{(N)}_0=\{x_0\}\cap B_{r_N/2}(x_c)\subset \mathcal{I}^{(N)}_0$. Assume the process is defined up to generation $g\geq 0$. Then conditional on $\sigma(\mathcal{S}^{(N)}_m, \mathcal{I}^{(N)}_m, \mathcal{R}^{(N)}_m, \mathcal{J}^{(N)}_m: m\leq g)$ let $\mathcal{J}^{(N)}_{g+1}\subset \mathcal{I}^{(N)}_{g+1}$ be the set of all infected hosts contained in $B_{r_N/2}(x_c)$ generated by previously infected hosts $x\in\mathcal{J}^{(N)}_{g}$. We again set $J_n^{(N)}:=|\mathcal{J}^{(N)}_{g}|$ for all $g\geq 0$.
\end{Def}

Note that by Lemma~\ref{Lemma: moderate deviation number of vertices} the the ball $B_{r_N/2}(x_c)$ is not empty.

For any sequence $(a_{N})_{N\in \mathbb{N}}$ define 
\begin{align}
    &\tau_{a_N}^{(N)}:=\inf \left\{g \in \mathbb{N}: I_{g}^{(N)} \geq a_N\right\}, \\ 
    &\overline{\tau}^{(N)}_{a_N}:=\inf \left\{g \in \mathbb{N}: \overline{I}^{(N)}_{g} \geq a_N\right\}.
\end{align}
\begin{Lem} (Coupling between $J^{(N)}$ and an infection process on a complete graph)\\ 
\label{Lem: coupling sub infection process and complete graph infection process}
Let $\left(H_{g}^{(N)}\right)_{g \in \mathbb{N}_{0}}$ be an infection process on a complete graph as in Section \ref{Section: Infection on complete graph} with $\widetilde{d}_{N}$ vertices where $\widetilde{d}_{N}$ is the number of vertices in $B\left(\frac{r_N}{2}\right)$, and with $\widetilde{v}_N:=\frac{v_N}{2^{n}}-N^{\frac{n+2}{n+3}\frac{\beta}{2}}$ offspring parasites. In particular we have $\widetilde{v}_N \sim \frac{a}{\sqrt{2^n}}\sqrt{\widetilde{d}_{N}}.$ Then one can show that $\mathcal{J}^{(N)}$ can be coupled whp from below with $H^{(N)}$ until $\overline{I}^{(N)}$ reaches the size $N^{\varepsilon}$, more precisely 
\begin{align}
    \mathbb{P} \left( \forall g \leq \overline{\tau}_{N^{\varepsilon}}, J_{g}^{(N)}\geq H^{(N)}_{g} \right) \underset{N \to \infty}{\longrightarrow}1. 
\end{align}
\end{Lem}

\begin{proof}
The complete graph has the same number of vertices as the number of vertices contained in $B\left(\frac{r_N}{2}\right)$.
Every time that an infection in the process $J^{(N)}$ occurs, $v_N$ parasites are generated, but only those that are moving to vertices in $B\left(\frac{r_N}{2}\right)$ will create infections that are counted in the sub infection process $J^{(N)}$. We need to control from below  how many among the $v_N$ parasites will move to a vertex in $B\left(\frac{r_N}{2}\right)$. This control needs to hold for at most $N^{\varepsilon}$ infections.  Using a similar approach as in Lemma~\ref{Lemma: moderate deviation number of vertices} and Remark~\ref{Remark: uniform bound number vertices} we can show that for each of at most $N^{\varepsilon}$ infections, the number of parasites, that are moving to vertices within $B(\frac{r_N}{2})$, is bounded from below by $\frac{v_N}{2^{n}}- N^{\frac{n+2}{n+3}\frac{\beta}{2}}$ 
whp, where we used that there are whp between $\frac{d_N}{2^{n}}\pm (2n+1)N^{\frac{n+2}{n+3}\beta}$ vertices in $B\left( \frac{r_N}{2}\right)$ and between $d_N \pm (2n+1)N^{\frac{n+2}{n+3}\beta}$ many neighbors for each of the infected vertices. Consequently, as we have set $\widetilde{v}_N=\frac{v_N}{2^{n}}-N^{\frac{n+2}{n+3}\frac{\beta}{2}}$  only less less infection con be created  in the infection process $H^{(N)}$ compared to $J^{(N)}$. \\ 
The only problem that could happen to make the coupling fail is that one infection generated from $I^{(N)}$ is on an empty vertex (generated due to $J^{(N)}$ or due to the global process $I^{(N)}$). This creates no infection in $J^{(N)}$ but it is actually be counted in $H^{(N)}$. However, such an event is possible only if in the original process, two pairs of parasites are pointing to the same vertex in $B \left( \frac{r_N}{2}\right)$ (at the same or at different generations). In particular this event is contained in the event that at least one vertex in $B\left( \frac{r_N}{2}\right)$ is hit by at least four parasites up to generation $\overline{\tau}_{N^{\varepsilon}}$. But until generation $\overline{\tau}_{N^{\varepsilon}}$  less than $N^\varepsilon$  vertices get infected cumulative over all generations. So it is possible to estimate from above the probability that such an event happens before generation $\overline{\tau}_{N^{\varepsilon}}$ by estimating the probability of the  event $A$ in the following experiment: Assume $N^{\varepsilon}v_N$ balls placed uniformly at random into  $\widetilde{d}_N$ boxes and we are interested in the event $A$ that it exists (at least) one box containing at least four balls. 
Indeed, the probability of the event $A$ gives an upper bound, all balls (corresponding to parasites in the original process) are put into $\tilde{d_N}$ boxes (parasites have a large choice of vertices to move to). 
This increases the probability for one box to contain four balls.

We upper bounded the probability of $A$ as follows 
\begin{align}
    &\mathbb{P} \left(A\right) \\
    &\hspace{1cm}\leq \left(\frac{d_N}{2^{n}}+(2n+1)N^{\frac{n+2 }{n+3}\beta}\right) \left(N^{\varepsilon}v_N\right)^4 \left(\frac{1}{\frac{d_N}{2^{n}}-(2n+1)N^{\frac{n+2}{n+3}\beta}} \right)^{4} =\mathcal{O}\left( \frac{N^{4 \varepsilon}}{d_N}\right)\underset{N \to \infty}{\longrightarrow}0, 
\end{align}
because $\varepsilon<\frac{\beta}{4}$.

So $J^{(N)}$ can be coupled  with $H^{(N)}$ such that whp $H^{(N)}_g\leq J^{(N)}_g$ for all $g \leq \overline{\tau}_{N^\varepsilon}$. 
\end{proof}

\begin{Lem} \label{Lem: reaching local infection RGG}
Using Lemma \ref{Lem: coupling sub infection process and complete graph infection process} and the results of Section \ref{Section: Infection on complete graph}, one can show that 
\begin{align}
    \liminf_{N \to \infty}\mathbb{P} \left(\exists g \in \mathbb{N}, \overline{I}^{(N)}_{g} \geq N^{\varepsilon}\right) \geq \pi\left(\frac{a}{\sqrt{2^{n}}}\right),
\end{align}
where $\pi \left(\frac{a}{\sqrt{2^n}}\right)$ is the survival probability of a DBPC with offspring and cooperation distribution $Poi \left(\frac{a^{2}}{2^{n+1}}\right)$ and $Poi \left( \frac{a^{2}}{2^n}\right)$. 
\end{Lem}

\subsubsection{Increasing from a total number of  $N^{\varepsilon}$ infections  to $N^{\frac{\beta}{2}+\widetilde{\varepsilon}}$ infections within a single box}

 Cover the space with non-overlapping boxes, such that all boxes except those having with the boundary of $[0,1]^n$ a non-empty intersection, have an edge length $r_N$ and such that one of the boxes is centered around $x_0^{(N)}$. Label the boxes and denote by $\mathbf{K}$ the set of labels and by $\mathcal{V}^{(N)}_{k}$ the set of vertices in box $k$. Furthermore, denote by $I_{g}^{(N)}(k)$ the number of infected vertices in box $k$ in generation $g$.
\begin{Lem}
\label{Lem: time control RGG}
\begin{align}
    \mathbb{P} \left( \exists g \leq \overline{\tau}_{N^{\varepsilon}}^{(N)}, \exists k \in \mathbf{K}, I_{g}^{(N)}(k) \geq \frac{N^{\varepsilon}}{\log(N)} \Big{\vert} \overline{\tau}^{(N)}_{N^{\varepsilon}}<\infty\right) \to 1. 
\end{align}
\end{Lem}
\begin{proof}
For any sequence $(a_N)_{N \in \mathbb{N}}$ introduce the following set
\begin{align}
    B^{(N)}_{a_N}:=\{k \in \mathbf{K}: \exists g \leq a_N, I^{(N)}_{g}(k)\geq 1\}.
\end{align}
At each generation a parasite may move a distance of at most $r_N$. So in dimension $n$, in $a_N$ generations, the number of balls of diameter $r_N$ that can be reached is $(2a_N+1)^{n}$. This gives that $\vert B_{a_N}^{(N)} \vert \leq (2a_N+1)^{n}$. \\ 
Using the coupling from below with the DBPC until generation $\overline{\tau}_{N^{\varepsilon}}^{(N)}$ and applying Proposition \ref{GrowFast} to the DBPC we obtain that it exists $C>0$ such that 
\begin{align} \label{lab: control time intial phase}
    \mathbb{P} \left(\overline{\tau}_{N^{\varepsilon}}^{(N)} \leq C\log(\log(N))\vert \overline{\tau}_{N^{\varepsilon}}^{(N)}<\infty\right) \to 1. 
\end{align}
Combining these two results we obtain that whp $\Big{\vert} B^{(N)}_{\overline{\tau}^{(N)}_{N^{\varepsilon}}} \Big{\vert} \leq \left(2C \log(\log(N))+1\right)^{n}$. \\ 
If for any generation before $\overline{\tau}^{(N)}_{N^{\varepsilon}}$ the number of infected individuals in any of the balls of $B^{(N)}_{\overline{\tau}^{(N)}_{N^{\varepsilon}}}$ is smaller than $\frac{N^{\varepsilon}}{\log(N)}$, this would mean that the total number of infected individuals up to generation $\overline{\tau}_{N^{\varepsilon}}^{(N)}$ would be upper bounded by $\frac{N^{\varepsilon}}{\log(N)}\left(2C \log(\log(N))+1\right)^{n}\log^{\frac{1}{n+2}}(N)=o(N^{\varepsilon})$, which gives a contradiction.
\end{proof}
Next let $\sigma^{(N)}$ be the stopping time, at which for the first time in one of the balls at least $\frac{N^{ \varepsilon}}{\log(N)}$ hosts get infected
\begin{align}
\sigma^{(N)}:= \inf \left\{g \in \mathbb{N}: \exists k \in \mathbf{K}, I_{g}^{(N)}(k) \geq \frac{N^{\varepsilon}}{\log(N)}\right\}. 
\end{align}
The last lemma exactly states that 
\begin{align}
    \mathbb{P} \left( \sigma^{(N)} \leq \overline{\tau}^{(N)}_{N^{\varepsilon}} \vert \overline{\tau}^{(N)}_{N^{\varepsilon}}< \infty \right) \underset{N \to \infty}{\longrightarrow}1. 
\end{align}
Now we will show that after a finite number of generations after generation $\sigma^{(N)}$, there is whp one box with at least $N^{\frac{\beta}{2}+\delta}$ infected vertices, where $\delta>0$ is small enough.\\
To achieve this goal, we will argue in the same way as we have done in Subsection \ref{Subsection: Final phase Epidemic complete graph}. \\
Choose $\varepsilon>0$ such that for all $g \in \mathbb{N}$,  $2^{g}\varepsilon\neq\frac{\beta}{2}$. Then define $\overline{g}$ as the largest $g \in \mathbb{N}$ satisfying  $2^{g+1}\varepsilon<\beta$. In particular we have $2^{\overline{g}+1}\varepsilon >\frac{\beta}{2}$ because otherwise we would have $2^{\overline{g}+2} \varepsilon <\beta$  which is a contradiction with the definition of $\overline{g}$. \\ 
Define 
\begin{align} \label{lab: S}
    S^{(N)}:= \left\{k \in \mathbf{K}: I_{\sigma^{(N)}}^{(N)}(k) \geq \frac{N^{\varepsilon}}{\log(N)}\right\}
\end{align}
the set of boxes that contain at least $\frac{N^{ \varepsilon}}{\log(N)}$ infected vertices in generation $\sigma^{(N)}$ . By definition of $\sigma^{(N)}$ if $\sigma^{(N)}<\infty$ then $S^{(N)}$ is not empty almost surely.
\begin{Lem}
\label{Lem: first step of the before last phase}
We have 
\begin{align}
    \mathbb{P} \left( \overline{I}^{(N)}_{\sigma^{(N)}}\leq \frac{N^{2 \varepsilon}}{\log(N)} \Big{\vert} \overline{\tau}^{(N)}_{N^{\varepsilon}}<\infty\right) \underset{N \to \infty}{\longrightarrow}1. 
\end{align}
\end{Lem}

\begin{proof}
By definition, at generation $\sigma^{(N)}-1$ the number of infected vertices in each box $i$ is at most $\frac{N^{\varepsilon}}{\log(N)}$ and the total number of boxes that have been infected is at most whp $\left(2C \log(\log(N))+1\right)^{n}$. To show that $\overline{I}^{(N)}_{\sigma^{(N)}} \leq \frac{N^{2 \varepsilon}}{\log(N)}$, we will control the number of infections in each box by applying a similar argument as in Lemma \ref{Lem: first step final phase} in the context of the complete graph. \\
At generation $\sigma^{(N)}-1$ whp we have $\overline{I}^{(N)}_{\sigma^{(N)}-1} \leq N^{\varepsilon}$ because otherwise we would have a contradiction to Lemma \ref{Lem: time control RGG}. Then to bound from above at generation $\sigma^{(N)}$ the total number of infected vertices up to
this generation, it suffices to add to $N^{\varepsilon}$ an upper bound on the number of new infections generated in generation $\sigma^{(N)}$. \\ 
In each box, there are at most $\frac{N^{ \varepsilon}}{\log(N)}v_N$ parasites that will move. Because of the sizes of the boxes, each box can receive infections from outside only due to its $3^{n}-1$ neighbouring boxes. To arrive at an upper bound on the number of new infections generated in a box, one can compare the situation with the following balls-into-boxes experiment: Consider $d_N:=\frac{N^{\beta}}{2^n}-(2n+1)N^{\frac{n+2}{n+3}\beta}-N^{\varepsilon}$ boxes. Put $\left(3^{n}-1\right)\frac{N^{\varepsilon}}{\log(N)}v_N$ balls into the boxes purely at random and count the number of boxes that contain at least two balls. Applying Lemma \ref{Lem: experiment boxes and balls} Equation \eqref{Eq: max number pair of balls} with $H(N)=\left(3^{n}-1\right) \frac{N^{\varepsilon}}{\log(N)}$, $\varphi_1(N)=N^{\varepsilon}$, $\varphi_2(N)=0$, $\varphi_3(N)=0$ $f_2(N)=\frac{1}{\left(3^{n}-1\right)^{2}}\frac{\log(N)}{\left(2C\log(\log(N))+1\right)^{n}}$, we obtain the result for each box. Choosing $f_{2}(N)=\Theta\left(\frac{\log(N)}{\log(\log(N))}\right)$ and because whp at most $\left(2C\log\log((N))+1\right)^{n}=\Theta \left(\log(\log(N))\right)$ boxes got so far infected, whp the result is true for all the boxes, according to Lemma \ref{Lem: experiment boxes and balls}, see Equation \eqref{Eq: max number pair of balls}. This implies that the number of empty vertices at generation $\sigma^{(N)}$ is whp at most $\left(3^{n}-1\right)^{2}\frac{N^{2\varepsilon}}{\log^{2}(N)}\left(2 C\log(\log(N))+1\right)^{n}f_{2}(N)=\frac{N^{2\varepsilon}}{\log(N)}$. 
\end{proof}
\begin{Lem}
\label{Lem: control infected and empty vertices}
Let $g \in \llbracket 0,\overline{g} \rrbracket$ we have 
\begin{align}
    \mathbb{P} \left( \forall k \in S^{(N)}, I^{(N)}_{\sigma^{(N)}+g}(k) \geq \frac{N^{2^{g} \varepsilon}}{\log^{\alpha_{g}}(N)}, \text{ and } \overline{I}^{(N)}_{\sigma^{(N)}+g} \leq \frac{N^{2^{g+1}\varepsilon}}{\log(N)} \Big{\vert} \overline{\tau}^{(N)}<\infty\right) \to 1,
\end{align}
where $\alpha_0=1$ and for all $g \geq 1, \alpha_{g}=2\alpha_{g-1}+2$. 
\end{Lem}
\begin{proof}
The proof is obtained by induction. First for $g=0$ the result is given by Lemma \ref{Lem: first step of the before last phase}.
Then let $g\leq \overline{g}-1$, assume the result is obtained for $0 \leq j \leq g$. Now we will show the result for $g+1$. \\ \\
To derive the lower bound on the number of infected vertices in a box $k \in S^{(N)}$ at generation $\sigma^{(N)}+g+1$, one can consider only the infections generated due to infected vertices inside this box. According to the induction hypothesis  there are at least $\frac{N^{2^{g}\varepsilon}}{\log^{\alpha_{g}}(N)}$ infected vertices in the box. Among the parasites generated on these vertices, at least $\frac{N^{2^{g}\varepsilon}}{\log^{\alpha_{g}}(N)}\left(\frac{1}{2^{n}}v_N-N^{\frac{\beta}{3}}\right)$ of them will move the vertices in the box. Then it suffices to apply Lemma \ref{Lem: experiment boxes and balls} Equation \eqref{Eq: min number pairs of balls} where $d_N$ in this Lemma is equal to $\frac{1}{2^n}N^{\beta}+(2n+1)N^{\frac{n+2}{n+3}\beta}$, with $H(N):=\frac{1}{2^{n}}\frac{N^{2^{g}\varepsilon}}{\log(N)}$, $\varphi_1(N)=0$, $\varphi_2(N)=\frac{N^{2^{g+1}\varepsilon}}{\log(N)}$, $\varphi_3(N)=2^{n}N^{\frac{\beta}{3}}$, $f_1(N)=\frac{\log^{2}(N)}{2^{n+1}}$, which gives that the number of infected vertices at generation $\sigma^{(N)}+g+1$ is whp at least of order $\frac{2^{n+1}}{\log^{2}(N)}\left(\frac{1}{2^{n}}\frac{N^{2^{g}\varepsilon}}{\log^{\alpha_{g}}(N)}\right)^{2}=\frac{N^{2^{g+1} \varepsilon}}{\log^{\alpha_{g+1}}(N)}$. Because there are whp at most  $\Theta \left(\log(\log(N))\right)$ boxes in $S^{(N)}$ and by Equation \eqref{Eq: min number pairs of balls} of Lemma \ref{Lem: experiment boxes and balls}, the statement holds whp for all boxes in $S^{(N)}$. \\
Indeed considering $\frac{1}{2^n}N^{\beta}+(2n+1)N^{\frac{n+2}{n+3}\beta}$ ball lower bounds the probability for a parasite to move to an occupied vertex, because whp there are at most $\frac{1}{2^n}N^{\beta}+(2n+1)N^{\frac{n+2}{n+3}\beta}$ many vertices in the box. Furthermore, according to the induction assumption we have considered the minimal number of parasites which is $\frac{N^{2^{g}\varepsilon}}{\log^{\alpha_{g}}(N)}\left(\frac{1}{2^{n}}v_N-N^{\frac{\beta}{3}}\right)$ and new infections are counted when reaching one of the $\frac{1}{2^n}d_N-\frac{N^{2^{g+1}\varepsilon}}{\log(N)}$ first boxes whereas in the original process there are at least this number of occupied vertices. 
\\ \\ 
To derive the upper bound on the number of empty vertices, we control for each box the number of new infections generated in generation $\sigma^{(N)}+g+1$. 
Since by induction  the number of empty vertices in generation $\sigma^{(N)}+g$ is  $\frac{N^{2^{g+1}\varepsilon}}{\log(N)}$ whp, 
we apply Lemma \ref{Lem: experiment boxes and balls} Equation \eqref{Eq: max number pair of balls} with $H(N)=\left(3^{n}-1\right)\frac{N^{2^{g+1}\varepsilon}}{\log(N)}$, $\varphi_1(N)=\frac{N^{2^{g+1}\varepsilon}}{\log(N)}$, $\varphi_2(N)=0$, $\varphi_3(N)=0$, $f_2(N)=\frac{1}{(3^{n}-1)^{2}}\frac{\log(N)}{\left(2\log^{\frac{1}{n+2}}(N)+g+2\right)^{n}}$ to estimate the number of new infection in generation $\sigma^{(N)} +g +1$ in each box in $B_{\sigma^{(N)}+g+1}^{(N)}$. The lemma yields that in each box there are at most $\frac{N^{2^{g+2}\varepsilon}}{\log^{2}(N)}\frac{\log^{(N)}}{\left(2 \log^{\frac{1}{n+2}}(N) +g+2 \right)^{n}}$ new infections whp. Since there are whp at most $\left(2C\log(\log(N))+g+2\right)^{n}$ boxes and since $f_{2}(N)=\Theta \left(\log(N)^{\frac{2}{n+2}}\right)$ whp for all boxes the number of new infections is bounded from above by $\frac{N^{2^{g+2}\varepsilon}}{\log^{2}(N)}\frac{\log^{(N)}}{\left(2 \log^{\frac{1}{n+2}}(N) +g+2 \right)^{n}}$, see Equation \eqref{Eq: max number pair of balls} of Lemma \ref{Lem: experiment boxes and balls}. 
Consequently, the total number of empty vertices at generation $\sigma^{(N)}+g+1$ is whp at most 
$\frac{N^{2^{g+2}\varepsilon}}{\log^{2}(N)}\frac{\log^{(N)}}{\left(2 \log^{\frac{1}{n+2}}(N) +g+2 \right)^{n}} (2 C  \log(\log(N)) + g+2)^n \leq \frac{N^{2^{g+2}\varepsilon}}{\log(N)}$.  
\end{proof}

Applying Lemma \ref{Lem: control infected and empty vertices} for $g=\overline{g}$ 
\begin{align}
    \mathbb{P}\left(\forall k \in S^{(N)}, I_{\sigma^{(N)}+\overline{g}}^{(N)}(k) \geq \frac{N^{2^{\overline{g}} \varepsilon}}{\log^{\alpha_{\overline{g}}}(N)}, \text{ and } \overline{I}_{\sigma^{(N)}+\overline{g}}^{(N)}\leq \frac{N^{2^{\overline{g}+1}\varepsilon}}{\log(N)} \Big{\vert} \overline{\tau}^{(N)}<\infty \right) \to 1.
\end{align}
Define $\delta=\frac{1}{2}\left(2^{\overline{g}+1}\varepsilon-\frac{\beta}{2}\right)>0$. 
In the next lemma we show that at generation $\sigma^{(N)}+\overline{g}+1$ the number of infected vertices in each box of $S^{(N)}$ is at least of order $N^{\frac{\beta}{2}+\delta}$. 
\begin{Lem} \label{Lem: boxes with sufficient infection}
\begin{align}
    \mathbb{P}\left(\forall k \in S^{(N)}, I_{\sigma^{(N)}+\overline{g}+1}^{(N)}(k) \geq N^{\frac{\beta}{2}+\delta} \vert \overline{\tau}^{(N)}<\infty \right) \to 1.
\end{align}
\end{Lem}
\begin{proof}
Here we apply again Lemma \ref{Lem: experiment boxes and balls} to obtain this lower bound. More precisely with the following set of parameters: $H(N)=\frac{1}{2^{n}}\frac{N^{2^{\overline{g}}\varepsilon}}{\log^{\alpha_{\overline{g}}}(N)}$, $\varphi_1(N)=0$, $\varphi_2(N)=\frac{N^{2^{\overline{g}+1}\varepsilon}}{\log(N)}$, $\varphi_3(N)=2^{n}N^{\frac{\beta}{3}}$, $f_1(N)=\frac{\log^{2}(N)}{2^{n+1}}$. We obtain that whp $I_{\tau^{(N)}+\overline{g}+1}\geq \frac{N^{2^{\overline{g}+1} \varepsilon}}{\log^{\alpha_{\overline{g}+1}}(N)} \geq N^{\frac{\beta}{2}+\delta}$, by definition of $\delta$. 
\end{proof}

\subsubsection{"Pulled travelling wave" epidemic spread:} 
\label{section; travelling wave}
We start with a general lemma that we will use multiple times in this subsection. It says that when a box of diameter $\varepsilon_N$ is fully infected, then in the next generation, all the vertices in the neighboring area of diameter $2r_N-\varepsilon$ are visited by at least two parasites whp.
\begin{Lem} \label{Lem: traveling wave}
    Consider a box of diameter $\varepsilon_N$ centered around a point $x \in [0,1]^{n}$, denoted by $B_2$, where $\varepsilon_N \in \Theta \left( N^{-\frac{1-\beta/2}{n}+\delta}\right)$ for $\delta>0$ small enough. Assume that the proportion of currently infected vertices in this box is asymptotically $1$. Then in the next generation it follows that all the vertices on the box centered around $x$ with diameter $2r_N-\varepsilon_N$, denoted by $B_1$, are visited by at least $2$ parasites with probability at least $1-\mathcal{O} \left(\varepsilon_N^n N v_N \exp \left(-\frac{\varepsilon_N^n}{ (2r_{N})^{n}}v_N\right)\right)$.
\end{Lem}
\begin{proof}
We estimate from above the probability that at least one vertex is visited by at most 1 parasite in the following generation. Denote by $K_1$ the number of vertices in $B_1$. One shows with the same kind of arguments as in Lemma \ref{Remark: uniform bound number vertices} that $(2r_N-\varepsilon_N)^{n}N-(2n+1) N^{\frac{n+2}{n+3}\beta} \leq K_1 \leq (2r_N-\varepsilon_N)^{n}N+(2n+1) N^{\frac{n+2}{n+3}\beta}  \in \Theta(d_N)$ whp. Denote by $\underline{K}_1=(2r_N-\varepsilon_N)^{n}N-(2n+1) N^{\frac{n+2}{n+3}\beta}$. Similarly, within $B_2$ whp at least $\varepsilon_N^{n}Nv_N/2\in \Theta(\varepsilon^n N v_N)$ parasites are generated, denote this number by $K_2$. Since $B_2$ is contained in $B_1$ and $B_1$ has diameter $2r_N-\varepsilon_N$ every vertex in $B_2$ is connected over an edge to any vertex in $B_1$. Hence, the probability $p$ that a particular vertex in $B_1$ gets attacked by at most one parasite generated in $B_2$ can be estimated from above by
\begin{align*}
p&  \leq \left(1 -\frac{1}{\underline{K}_1} \right)^{K_2} +  K_2 \left(1 -\frac{1}{d_{N,l}} \right)^{K_2-1} \frac{1}{\underline{K}_1} \\ & \leq  
\left(1 -\frac{1}{\underline{K}_1} \right)^{\varepsilon_N^n N v_N /2}
+ 2 \varepsilon_N^n Nv_N \left(1 -\frac{1}{\underline{K}_1} \right)^{\varepsilon_N^n Nv_N /2-1} \frac{1}{\underline{K}_1}. 
\end{align*}
Consequently, we can estimate from above the probability that at least one of the vertices gets attacked by at most one parasite by 
\begin{align*}
& K_1 \left(\left(1 -\frac{1}{\underline{K}_1} \right)^{\varepsilon_N^n Nv_N /2}
+ 2 \varepsilon_N^n N v_N \left(1 -\frac{1}{\underline{K}_1} \right)^{\varepsilon_N^n Nv_N/2-1} \frac{1}{\underline{K}_1}\right) \\ & \in \mathcal{O}\left(\varepsilon_N^n N v_N \exp\left(-\frac{\varepsilon_N^n}{(2 r_N)^n} v_N \right)\right),
\end{align*}
from which follows the claim.
\end{proof}
This lemma implies that within one generation all vertices in the box $B_1$ get attacked by at least two parasites whp from parasites generated on vertices in $B_2$, since $\varepsilon_N$ is chosen such that $\varepsilon_N^n v_N/ r_N^n \sim N^{n\delta}$. In particular this means that if in $B_2$ almost all the vertices that are not contained in $B_1$ (which is asymptotically completely infected) are still occupied by a host, then all these hosts in $B_2$ get infected in the next generation. This allows us to repeat the same argument subsequently. Due to the exponential decay of the error term $\mathcal{O} \left(\varepsilon_N^n N v_N \exp \left(-\frac{\varepsilon_N^n}{ (2r_{N})^{n}}v_N\right)\right)$ we can apply this argument for many boxes, in particular for the $2^{n}N^{1-\beta}$ many boxes of diameter $r_N$. In particular, this implies that we can show that a pulled traveling wave in any direction is created by repeating the argument, as long as the invasion is not stopped by a region in which no susceptible hosts are available anymore and which cannot be crossed by parasites. \\
However, such a region with a non-trivial proportion of hosts killed and with a diameter of at least $r_N-\varepsilon_N$ (such that it cannot be crossed by a travelling wave whp) cannot arise, because if in a box of size $r_N-\varepsilon_N$ at least $k_N$, with $k_N \rightarrow \infty$ arbitrarily slowly, hosts get infected by couplings with  DBPCs (which have a positive survival probability) we can show that in this region either a new infection wave is started or a travelling wave is hitting the box. Consequently boxes cannot be slowly depleted and we will reach the boundary of $[0,1]^n$ whp after at most
\begin{align} \label{lab: control time invasion}
    \frac{1}{2(r_N - \varepsilon_N)}=\frac{1}{2r_N}\left(1+\mathcal{O} \left(\frac{\varepsilon_N}{r_N}\right)\right),
\end{align}
many generations.
\subsection{
Proof of Theorem \ref{MainResult} 1) (ii):}

Now we have all necessary materials to prove Theorem \ref{MainResult} 1) (ii). \\
The first step is to show 
\begin{align}\label{lab: upper bound invasion probability random geometric graph}
    \limsup_{N \to \infty} \mP \left(E_{u}^{(N)}\right) \leq \pi(a).
\end{align}
For a sequence $\left(\ell_N\right)_{N \in \mathbb{N}}$ introduce the event 
\begin{align}
    E^{(N)}_{\ell_N}:= \left\{\exists g \in \mathbb{N}_{0}: \overline{I}^{(N)}_{g} \geq \ell_N\right\}. 
\end{align}
Then it follows that for all $0<u\leq 1$ and any sequence $\ell_N \leq u \vert \mathcal{V}^{(N)}\vert$ we have 
\begin{align} \label{lab: global infection means local infection}
    \mathbb{P} \left(E_{u}^{(N)}\right) \leq \mP \left(E^{(N)}_{\ell_N}\right).
\end{align}
Taking a sequence $\ell_N$ satisfying $\ell_N \to \infty$ and $\ell_N\in \log(\log(N))$ we have by Proposition \ref{UpperBoundRGG} that 
\begin{align} \label{lab: coupling upper DBPC RGG}
    \mathbb{P} \left( E^{(N)}_{\ell_N}\right) \leq \mP \left(\exists g \in \mathbb{N}_{0}: \overline{Z}^{(N)}_{g,u} \geq \ell_N\right)+o(1).
\end{align}
Proposition \ref{SurvivalProb Upper DBPC RGG} gives that 
\begin{align} \label{lab: control survival probability upper DBPC}
    \lim_{N \to \infty} \mP \left(\exists g \in \mathbb{N}_{0}: \overline{Z}^{(N)}_{g,u} \geq \ell_N\right)=\pi(a).
\end{align}
In summary combining \eqref{lab: global infection means local infection}, \eqref{lab: coupling upper DBPC RGG} and \eqref{lab: control survival probability upper DBPC} gives exactly \eqref{lab: upper bound invasion probability random geometric graph}.

The second step is to show 
\begin{align}\label{lab: lower bound invasion probability random geometric graph}
    \liminf_{N \to \infty} \mP \left(E_{u}^{(N)}\right) \geq \pi\left(\frac{a}{\sqrt{2^n}}\right).
\end{align}
Lemma \ref{Lem: reaching local infection RGG} gives that 
\begin{align} \label{lab: proba initial phase RGG}
    \liminf_{N \to \infty}\mathbb{P} \left(\exists g \in \mathbb{N}, \overline{I}^{(N)}_{g} \geq N^{\varepsilon}\right) \geq \pi\left(\frac{a}{\sqrt{2^{n}}}\right).
\end{align}
Then Lemma \ref{Lem: boxes with sufficient infection} allows to get that for all boxes in the set $S^{(N)}$, see \eqref{lab: S} for its definition, in a random generation the infection level in these boxes is at least of order $N^{\beta/2+\delta}$ as mentioned below 
\begin{align}
    \mathbb{P}\left(\forall k \in S^{(N)}, I_{\sigma^{(N)}+\overline{g}+1}^{(N)}(k) \geq N^{\frac{\beta}{2}+\delta} \vert \overline{\tau}^{(N)}<\infty \right) \to 1.
\end{align}
Then arguing as in the proof of Lemma \ref{Lem: every vertices are killed} one can show that all the vertices in the boxes of the set $S^{(N)}$ get killed in one more generation, that is to say 
\begin{align}
    \mP \left(\forall k \in S^{(N)}, \overline{I}^{(N)}_{\sigma^{(N)}+\overline{g}+2}(k)=\vert \mathcal{V}_{k}^{(N)}\vert \right)\to 1. 
\end{align}
And finally using the results from Subsection \ref{section; travelling wave} one shows that whp every vertex eventually gets infected conditioned on the event $\{\overline{\tau}^{(N)}<\infty\}$, which combined with \eqref{lab: proba initial phase RGG} gives \eqref{lab: lower bound invasion probability random geometric graph}.

\subsection{Proof of Theorem \ref{MainResult} 1) (i) }

Assume $v_N \in o\left(\sqrt{d_N}\right)$. Then using a similar approach as in Subsection \ref{lab: proof subcritical case complete graph} which is to show that whp there are no infected individuals at generation $1$, one obtains the result.

\subsection{Proof of Theorem \ref{MainResult} 1) (iii) }

In this section we assume $\sqrt{d_N} \in o(v_N)$. We will prove that 
\begin{align}
    \lim_{N \to \infty} \mathbb{P} \left( \exists g \in \mathbb{N}_{0}: \overline{I}^{(N)}_{g}=\vert \mathcal{V}^{(N)}\vert\right)=1. 
\end{align}

\begin{proof}
The proof is split into two parts. First we argue that we can reach with high probability a level $N^{\alpha}$ for any $\alpha<\beta$ in a time of order $\log^2(N)$. In the second part we show that similar as in the critical scaling the host population is killed by a traveling wave.

For the first part we closely follow the proof of Lemma \ref{Lemma: InvasionTildeProcess} and the proof strategy in Subsection \ref{Section: Proof of Theorem 2.7 (iii)}.  We build in the same way an infection process $\big(\widehat{\mathcal{S}}^{(N)},\widehat{\mathcal{I}}^{(N)},\widehat{\mathcal{R}}^{(N)}\big)$, in which  infections are only transmitted due to parasites originating from the same vertex and  $v_N^{(a)}=a \sqrt{d_N}=a N^{\frac{\beta}{2}}$ many parasites are generated. This means a host is only infected if at least two parasites which originate from the same vertex attack the host simultaneously. Note that $d_N\in o(v_N)$, which means that for every $a>0$ there exists an $N$ large enough such that  $v_N^{(a)}\leq v_N$. Thus, analogously as we showed in Subsection \ref{Section: Proof of Theorem 2.7 (iii)}, for every $a>0$ we can couple this process to the original infection process $\big(\mathcal{S}^{(N)},\mathcal{I}^{(N)},\mathcal{R}^{(N)}\big)$ such that 
\begin{equation} \label{lab: coupling only CoSame RGG}
    \widehat{\mathcal{I}}_g^{(N)}\cup \widehat{\mathcal{R}}_g^{(N)}\subset \mathcal{I}_g^{(N)}\cup \mathcal{R}_g^{(N)}, \forall g \in \mathbb{N}_{0},
\end{equation}
for $N$ large enough. Denote by 
$H_x^{(N)}$ the  number of vertices which get attacked by at least two parasites  originating from $x$.   
Denote by $\text{deg}(x)$ the degree of  vertex $x\in \mathcal{V}^{(N)}$. Then 
\begin{align*}
    \mP\left(H_x^{(N)}=k\right)\geq & \frac{\prod_{i=1}^{k}\binom{v_N-2(i-1)}{2}}{k!\text{deg}(x)^{k}}\frac{\deg(x)!}{(\deg(x)-k- v_N)!\text{deg}(x)^{v_N-k}},
\end{align*}
where we only consider infections resulting from cooperation from the same edge and ignore infections generated by groups of 3 or more parasites, since these events happen with a negligible probability. Recall that $\delta_{N,\ell} = N^\beta -(n+1)N^{\frac{n+2}{n+3}\beta}$,  $\delta_{N,u} = N^\beta +  (2n+1)N^{\frac{n+2}{n+3}\beta}$ and set  
\begin{equation*}
    A^{(N)}:=\big\{\delta_{N,\ell}\leq \text{deg}(x)\leq \delta_{N,u} \,\forall \, x\in \mathcal{V}^{(N)}\big \}.
\end{equation*}
By Lemma~\ref{Remark: uniform bound number vertices} it follows that $\mP\big(A^{(N)}\big)\to 1$ as $N\to\infty$. Thus, $\delta_{N,\ell} $ and $\delta_{N,u}$ act as a uniform lower and upper bound on $\text{deg}(x)$ for all $x\in \mathcal{V}^{(N)}$ with high probability and we can again conclude  analogously as in Proposition 3.5 in \cite{BrouardEtAl2022} that 
\begin{align*}
    &\frac{\prod_{i=1}^{k}\binom{v_N-2(i-1)}{2}}{k!\text{deg}(x)^{k}}\frac{\deg(x)!}{(\deg(x)-k- v_N)!\text{deg}(x)^{v_N-k}}\\
    \geq &\left(\dfrac{(v_{N}-2a_{N})^{2}}{2\text{deg}(x)} \right)^{k}\frac{1}{k!}\exp \left(-\frac{v_{N}^{2}}{2\text{deg}(x)}\right) \left(1-\frac{1}{\text{deg}(x)^{\delta}} \right)\\
    \geq &\left(\dfrac{(v_{N}-2a_{N})^{2}}{2\delta_{N,u}} \right)^{k}\frac{1}{k!}\exp \left(-\frac{v_{N}^{2}}{2\delta_{N,\ell}}\right) \bigg(1-\frac{1}{\delta_{N,\ell}^{\delta}} \bigg)
\end{align*}
for $0\leq k\leq a_N$. This suggests that we can  couple the process $\big( |\widehat{\mathcal{I}}^{(N)}_g\cup \widehat{\mathcal{R}}^{(N)}|_g\big)_{g\geq 0}$ with an appropriately chosen branching process until we reach a level $N^{\alpha}$ with $\alpha<\beta$.
\begin{Def}(Modified lower Galton-Watson Process)\label{modifcation lower branching process}
Let $0<\delta<\frac{1}{2}$ and $(a_N)_{N \in \N}$ be a sequence with $a_{N} \rightarrow \infty$ and $a_N \in o\left( N^{\frac{\beta}{2}}\right)$. Furthermore assume $(\vartheta_N)_{N \in \N}$ is a $[0,1]$-valued sequence with
$\vartheta_{N} \rightarrow 0$ as $N\to \infty$. 
Let $\mathbf{Y}_{l}^{(N)}=\big( Y^{(N)}_{g,l}\big)_{g \in \N_0}$ be a Galton-Watson process with mixed binomial offspring distribution $\text{Bin}\big(\widehat{Y}^{(N)},1-\vartheta_{N}\big)$, where the probability weights $\big(\hat{p}_k^{(N)}\big)_{k \in \N_0}$ of $\widehat{Y}^{(N)}$ are
for all   $ 1 \leq j \leq a_{N}$
\begin{linenomath*}
    \begin{equation}
    \hat{p}_j^{(N)}:=\left(\dfrac{(v_{N}-2a_{N})^{2}}{2\delta_{N,u}} \right)^{k}\frac{1}{k!}\exp \left(-\frac{v_{N}^{2}}{2\delta_{N,\ell}}\right) \bigg(1-\frac{1}{\delta_{N,\ell}^{\delta}} \bigg),
    \end{equation}
    \end{linenomath*}
    and 
    \begin{linenomath*}
    \begin{equation}
    \hat{p}_0^{(N)}:=1-\sum_{j=1}^{a_N}\hat{p}_j^{(N)}.
    \end{equation}
\end{linenomath*}
Denote by 
$\overline{Y}^{(N)}_{g,l} := \sum_{i=0}^{g} Y_{i,l}^{(N)}$ the total size of the Galton-Watson process until generation $g$ and by $\overline{\mathbf{Y}_{l}} = \big(\overline{Y}_{g,l}^{(N)}\big)_{g \in \mathbb{N}_0} $ the corresponding process.  
\end{Def}
Now let $0<\alpha<\beta$ and define $\overline{\sigma}^{(N)}_{N^{\alpha}}=\inf \left\{ g \in \mathbb{N}_{0}: \vert \widehat{\mathcal{I}}^{(N)}_{g} \cup \widehat{\mathcal{R}}^{(N)}_{g}\vert \geq N^{\alpha}\right\}$. One can show similarly as in proof of Lemma \ref{Lemma: InvasionTildeProcess} Equation \eqref{lab: coupling below GWP} that
\begin{equation}
    \mP \left( \vert \widetilde{\mathcal{I}}^{(N)}_{g} \cup \widetilde{\mathcal{R}}^{(N)}_{g} \vert \geq \overline{X}^{(N)}_{g,l} \, \forall \, g \leq \overline{\sigma}^{(N)}_{N^{\alpha}}\right) \to 1
\end{equation}
as $N\to \infty$. Indeed, as in the proof of Lemma \ref{Lemma: InvasionTildeProcess} essentially we need to control the probability that a) an already empty vertex is re-attacked by at least two parasites moving along the same edge or b) a vertex gets simultaneously attacked by several pairs of parasites moving along the same edge.  

In the following we will call pairs of parasites moving along the same edge packs of parasites.
Similar as before we need to determine that each pack of parasites generated by an infected vertex before generation $\overline{\sigma}^{(N)}_{N^{\alpha}}$ is involved in one of the events a) or b) (independently of the other packs of parasites) with probability at most $\vartheta_N$. In this case we can remove packs of parasites with probability $\vartheta_N$ such that the number of new infections generated by an infected hosts can with high probability be bounded from below by the number of offspring drawn according to the distribution with weights $(p_{k, l}^{(N)})_{k\in \mathbb{N}_0}$ from Definition \ref{modifcation lower branching process} for any generation $n \leq\overline{\sigma}^{(N)}_{N^{\alpha}}$.

Next we determine an upper bound on the probabilities of the events a) and b).
\begin{itemize}
    \item[a)] Before generation $\overline{\sigma}^{(N)}_{N^{\alpha}}$ the probability that a pack of parasites originating from a vertex $x$ attacks an already empty vertex is bounded from above by 
    \begin{equation*}
        \frac{N^{\alpha}}{\text{deg}(x)}\leq \frac{N^{\alpha}}{\delta_{N,\ell}^{\delta}}=\frac{N^{\alpha}}{N^{\beta} -(d+1)N^{\frac{d+2}{d+3}\beta}}\in \Theta\Big(\frac{1}{N^{\beta-\alpha}}\Big).
    \end{equation*}
    \item[b)] Before generation $\overline{\sigma}^{(N)}_{N^{\alpha}}$, the number of empty vertices in the graph is smaller than $N^{\alpha}$. The probability that two packs of parasites coming from 2 different vertices $x$ and $y$ attack the same vertex is bound by
    \begin{equation*}
        \frac{|(\mathcal{N}_{x}\cap \mathcal{N}_{y})\backslash\{x,y\}|}{\text{deg}(x)\text{deg}(y)}\leq \frac{1}{\text{deg}(x)\wedge \text{deg}(y)}\leq \frac{1}{\delta_{N,l}}\in \Theta\Big(\frac{1}{N^{\beta}}\Big),
    \end{equation*}
    where $\mathcal{N}_{x}=\{y\in \mathcal{V}^{(N)}: x\in\mathcal{E}^{(N)}\}$ denotes the neighbourhood of $x\in \mathcal{V}^{(N)}$. An application of Markov's inequality yields that the total number of packs of parasites generated before generation  $\overline{\sigma}^{(N)}_{N^{\alpha}}$ is with high probability bounded by $N^{\alpha}\log(N)$, as in Lemma 4.8 in \cite{BrouardEtAl2022}. Hence, each pack of parasites is involved in an event of type $b)$ with probability at most $N^{\alpha}\log(N) \cdot (\delta_{N,l})^{-1}=\Theta \left(\frac{\log(N)}{N^{\beta-\alpha}}\right)$.
\end{itemize}
Set $\vartheta_N:=2 \frac{N^\alpha \log(N)}{\delta_{N,l}}\in \Theta \left(\frac{\log(N)}{N^{\beta-\alpha}}\right)$, then $\vartheta$ is an upper bound on the probability that a pack of parasites is involved in one of the events of type a) or b). For $\alpha < \beta$ we have $\vartheta_N \in o(1)$. By the exact same line of arguments as in Lemma \ref{Lemma: InvasionTildeProcess} one can conclude that
\begin{equation}\label{lab: proba 1 local infection}
    \lim_{N \to \infty} \mP \left( \exists g \in \mathbb{N}_{0}:  \overline{I}^{(N)}_{g} \geq N^{\alpha}\right)=1
\end{equation}
for any $\alpha<\beta$. Using the same approach as in the proof of Lemma \ref{Lem: time control RGG}, where the only difference is a coupling from below with Galton-Watson processes instead of DBPC, one can show that under the event $\left\{\exists g \in \mathbb{N}_{0}:  \overline{I}^{(N)}_{g} \geq N^{\alpha}\right\}$  it exists a box of diameter $r_N$, in which at least  $\frac{N^{\alpha}}{\log^2(N)}$ hosts got infected.  \\
Taking $\frac{\beta}{2}<\alpha<\beta$ and using a similar approach as in the proof of Lemma \ref{Lem: every vertices are killed}, one shows that it exists a box of diameter $r_N$ where all the hosts are killed and the number of infected individuals is of the order $N^{\beta}$. Then arguing as in Subsection  \ref{section; travelling wave} one shows that whp every vertices are killed by the infection process. Combined with \eqref{lab: proba 1 local infection} the result follows. 
\end{proof}

\subsection{Proof of Theorem \ref{MainResult} 2)}
\begin{proof}[Proof of Theorem \ref{MainResult} 2)]

At each generation a parasite can move at most to a distance $r_N$ meaning that the minimal number of generations the infection process $\overline{\mathbf{I}}^{(N)}$ needs for killing every host is at least the number of boxes of diameter $r_N$ that separates the initial vertex to the boundaries of the domain. In dimension $n$, using the max norm, this number is exactly $\frac{1}{2r_N}$, giving the result $\mP \left(\left\lfloor \frac{1}{2r_N}\right\rfloor \leq T^{(N)}\right)=1.$ \\
According to \eqref{lab: control time intial phase}, Lemma \ref{Lem: boxes with sufficient infection} and applying a similar reasoning as in Lemma \ref{Lem: every vertices are killed}, one shows that it exists $C>0$ such that 
\begin{align}
    \mP \left(  \sigma^{(N)}\leq C \log(\log(N))\Big{\vert} \sigma^{(N)}<\infty\right) \underset{N \to \infty}{\longrightarrow} 1,
\end{align}
where $\sigma^{(N)}:= \inf \left\{ g \in \mathbb{N}: \exists k \in \mathbf{K}, \overline{I}^{(N)}_{g}(k)= \vert \mathcal{V}^{(N)}(k)\vert\right\}$. Moreover Equation \eqref{lab: control time invasion} gives that after time $\sigma^{(N)}$, the remaining time up to total infection is upper bounded by $\frac{1}{2r_N}\left(1+\mathcal{O} \left(\frac{\varepsilon_N}{r_N}\right)\right)$. Combining these two facts gives that 
\begin{align}
    \mP \left( T^{(N)} \leq  \left\lceil \frac{1}{2r_N} \right\rceil+\mathcal{O}\left(\max \left(\log(\log(N)),\frac{\varepsilon_N}{r_N^2}\right)\right)\Big{\vert} T^{(N)}<\infty\right)\to 1. 
\end{align}
\end{proof}

\paragraph*{Acknowledgements}
MS and HT were supported by the LOEWE programme of the state of Hessen (CMMS).
CP acknowledges support from the German Research Foundation through grant PO-2590/1-1.
HT acknowledges support from the German Research Foundation through grant ME-2088/5-1.
VB and CP acknowledge support during the JTP 2022 "Stochastic Modelling in the Life Science"  funded by the Deutsche Forschungsgemeinschaft (DFG, German Research Foundation) under Germany's Excellence Strategy – EXC-2047/1 – 390685813.

\section{Supplementary material}

\begin{proof}[Proof of Lemma \ref{Lemma: extinction-explosion principle DBPC}]
If $p_{0,o}=0$ the result follows by applying the extinction-explosion principle to the super-critical Galton-Watson process formed by the offspring generated by the initial individual. Due to assumption $p_{1,o}\neq 1$ this process is super-critical. 

For the remaining cases we will show that all states except 0 are transient states, which yields the result.

Assume that $p_{0,o}\neq 0$ and that $p_{0,c}=0$. Note that this means that we get at least one offspring from every possible cooperation of parents. Thus if we have $Z_g=k$ parents, we get at least $\frac{k(k-1)}{2}$ many offspring due to cooperation. But it holds that $\frac{k(k-1)}{2}>k$ for $k\geq 4$. Thus, if at some generation $n$ we have that $Z_g \geq 4$, then we know that $Z_{g+1}>Z_g$ almost surely. This implies that $Z_g\to \infty$ as $n\to \infty$ almost surely, if $Z_{g_0} \geq 4$ for some $g_0\in \mathbbm{N}$. On the other hand since we exclude that $p_{0,o}=1$ and $p_{1,c}=1$ we have
\begin{align*}
    \mP(0<Z_g\leq 3|0<Z_{g-1}\leq 3\big)\leq c 
\end{align*}
for some $c<1$. Consequently,
the event $\{0<Z_g\leq 3\,\forall n\geq 0\}$ is a null-set and so all states but 0 are transient.

Assume that $p_{0,o}>0$ and $p_{0,c}>0$.  
Let us assume that in some generation $g_0$ we have $Z_{g_0}=k$ for some $k\geq 1$. If the process dies out in the next generation it enters the trap $0$ such that it can never return to $k$. Thus, an obvious lower bound for the probability to never hit $k$ again is
\begin{equation*}
    \mP(Z_g\neq k \forall g \geq g_0| Z_{g_0}=k)\geq p_{0,o}^k p_{0,c}^{\binom{k}{2}}>0 
\end{equation*}
for $k\geq 1$. But this already implies that $\mP(Z_g =k \text{ for some } g> g_0 |Z_{g_0}=k) <1$, i.e.\ the state $k$ is transient.

\end{proof}

A consequence of the extinction-explosion principle is the following lemma, which states that the probability of reaching an arbitrary high level, that tends to $\infty$, at some generation or up to some generation is asymptotically the same as surviving. It is a special case of Proposition \ref{Lemma: convergence survival probability DBPC}. We need it to prove Proposition \ref{Lemma: convergence survival probability DBPC} and other statements.

\begin{Lem}\label{lem:ReachingLevelToSurvival}
Let $\ZZ$ be a DBPC with survival probability $\pi>0$ and satisfying $p_{1,o}\neq 1$ and $(p_{0,o},p_{1,c})\neq (1,1)$. Then for any sequence $(b_N)_{N \in \mathbb{N}}$  satisfying $b_N \to \infty$ we have 
\begin{align}
    \lim_{N \to \infty}\mathbb{P} \left( \exists g \in \mathbb{N}_{0}: \overline{Z}_{g} \geq b_N\right)=\lim_{N \to \infty}\mathbb{P} \left( \exists g \in \mathbb{N}_{0}: Z_{g} \geq b_N\right)=\pi. 
\end{align}
\end{Lem}
The proof follows basically along the same arguments as the corresponding Lemma 3.7 in \cite{BrouardEtAl2022}.

\begin{proof}[Proof of Lemma \ref{lem:ReachingLevelToSurvival}]
First we will show that $\mathbb{P} \left( \exists g \in \mathbb{N}_{0}: Z_{g} \geq b_N\right) \to \pi$. By the extinction-explosion principle for DBPC, proven in Lemma \ref{Lemma: extinction-explosion principle DBPC}, we have that $\pi \leq \mathbb{P}\left( \exists g \in \mathbb{N}: Z_g \geq b_N\right)$. Then 
\begin{align}
    \pi&=\mathbb{P} \left( Z_g>0, \forall g \in \mathbb{N}_0\right) \\
    &=\mathbb{P} \left( \left\{\exists g \in \mathbb{N}_{0}: Z_g \geq b_N \right\} \cap \left\{Z_g>0, \forall g \in \mathbb{N}_0\right\}\right) \\
    & =\mathbb{P} \left(\{\exists g \in \mathbb{N}_{0}: Z_{g} \geq b_N\} \right) \mathbb{P} \left( \{Z_{g}>0, \forall g \in \mathbb{N}_{0}\}\vert \{\exists g \in \mathbb{N}_{0}: Z_{g} \geq b_N\}\right).
\end{align}
Using the strong Markov property, one can show that 
\begin{align}
    \mathbb{P} \left( \{Z_{g}>0, \forall g \in \mathbb{N}_{0}\}\vert \{\exists g \in \mathbb{N}_{0}: Z_{g} \geq b_N\}\right) \geq \mathbb{P} \left(Z_{g}>0, \forall g \in \mathbb{N}_{0}\vert Z_{0}=b_{N}\right). 
\end{align}
Then because the interaction is a cooperative one, a DBPC starting in $b_N$ can be coupled with $b_{N}$ independent DBPCs starting in 1  such that we get 
\begin{equation}\label{Equation: DBPC Superbranching property}
    \mathbb{P} \left(\exists g \in \mathbb{N}: Z_{g}=0\vert Z_{0}=b_{N}\right) \leq \left(\mathbb{P}\left(\exists g \in \mathbb{N}: Z_{g}=0\vert Z_{0}=1\right)\right)^{b_{N}}.
\end{equation}
Introducing $q<1$ as the extinction probability of the DBPC, we finally obtain 
\begin{align}
    \pi \geq \mathbb{P} \left( \left\{\exists g \in \mathbb{N}_{0}: Z_g \geq b_N \right\} \right) \left( 1-q^{b_N}\right).
\end{align}
It follows that
\begin{align}
    \mathbb{P} \left( \left\{\exists g \in \mathbb{N}_{0}: Z_g \geq b_N \right\} \right) \leq \frac{\pi}{1-q^{b_N}}\to \pi. 
\end{align}
Hence we have shown that
\begin{align}
\label{Equation: DBPC survives when reaching high number}
    \mathbb{P} \left( \exists g \in \mathbb{N}_{0}: Z_{g} \geq b_N\right) \underset{N \to \infty}{\longrightarrow} \pi.
\end{align}
For proving the remaining equality it remains to show that  
\begin{equation*}
    \mP \left(\Big{\{}\exists \ g \in \mathbb{N}_0: \overline{Z}_g \geq b_{N} \Big{\}} \cap \Big{\{} \exists \ g \in \mathbb{N}_0: \  Z_g=0\Big{\}}\right)=o(1)
\end{equation*}
due to the extinction-explosion principle for DBPCs.
Let $(c_N)_{N\in \mathbb{N}}$ be a sequence with $c_N \underset{N \to \infty}{\to} \infty$ and $\frac{b_{N}}{c_{N}} \underset{N \to \infty}{\to} \infty$ and consider the subsets
\begin{align}
        &A^{(N)}:=\Big{\{}\exists \ g \in \mathbb{N}_0:  \overline{Z}_g \geq b_{N}, \ \exists \ i\leq g, \ Z_i \geq c_N \Big{\}} \cap \Big{\{} \exists \ g \in \mathbb{N}_0: \  Z_g=0\Big{\}}, \\ 
        &B^{(N)}:=\Big{\{}\exists \ g \in \mathbb{N}_0:\overline{Z}_g \geq b_{N}, \ \forall \  i\leq g, \ Z_i < c_N \Big{\}} \cap \Big{\{} \exists \ g \in \mathbb{N}_0: \  Z_g=0\Big{\}}.
\end{align}
By definition
\begin{equation*}
        \Big{\{}\exists \ g \in \mathbb{N}_0: \overline{Z}_g \geq b_{N} \Big{\}} \cap \Big{\{} \exists \ g \in \mathbb{N}_0: \  Z_g=0\Big{\}}=  A^{(N)} \sqcup B^{(N)}.
\end{equation*}
The extinction-explosion principle together with \eqref{Equation: DBPC survives when reaching high number} yields that
\begin{equation*}
  \mP \left( A^{(N)}  \right) \leq  \mP \left( \Big{\{} \exists \ g \in \N, \ Z_g \geq c_N \Big{\}} \cap \Big{\{} \exists \ g \in \mathbb{N}_0, \  Z_g=0\Big{\}}\right) \underset{N \to \infty}{\to}0.
\end{equation*}

Furthermore 
\begin{equation*}
        B^{(N)} \subset \Big{\{} Z_{\floor{\frac{b_N}{c_N}}}>0 \Big{\}} \cap \Big{\{} \exists \ g \in \mathbb{N}_0, \  Z_g=0 \Big{\}}, 
\end{equation*}
which gives
\begin{equation*}
\begin{split}
    \mP \left({B^{(N)}}^c\right) \geq & \mP\left(\Big{\{} Z_{\floor{\frac{b_N}{c_N}}}= 0 \Big{\}} \sqcup \Big{\{} \forall \ g \in \mathbb{N}, \  Z_g>0 \Big{\}} \right)\\&=\mP \left( \Big{\{} Z_{\floor{\frac{b_N}{c_N}}}=0 \Big{\}}\right)+\mP \left( \Big{\{}\forall \ g \in \mathbb{N}, \  Z_g>0 \Big{\}}\right) \\ 
    &\to 1,
\end{split}
\end{equation*}
because for any sequence $u_N \to \infty$ we have $\mathbb{P}\left(Z_{u_N}=0\right)\to 1-\pi$, which follows by monotonicity of the events since $\{Z_{g+1}=0\}\subset \{Z_{g}=0\}$ for all $g\geq0$. Hence, we have $\mP \left( A^{(N)} \sqcup B^{(N)} \right) \underset{N \to \infty}{\to} 0.$
\end{proof}

\begin{proof}[Proof of Lemma \ref{lem:technicalBound}]
We set
\begin{align*}
    E(k):=\E[Z_g|Z_{g-1}=k]=\Big(k \mu_o + \binom{k}{2} \mu_c\Big)\geq \binom{k+1}{2}\mu,\\
    V(k):=\V(Z_g|Z_{g-1}=k)=\Big(k \nu_o + \binom{k}{2} \nu_c\Big)\leq \binom{k+1}{2}\nu,
\end{align*}
where $\mu:=\min\{\mu_o,\mu_c\}$ and $\nu:=\max\{\nu_o,\nu_c\}$.

Since we assumed that the first and second moments of the offspring and cooperation distributions converge, it exists a $N_0$ such that  
 \begin{align*}
     \frac{1}{2} \bigg(\frac{k(k+1)}{2} \mu\bigg)\leq  E^{(N)}(k) \qquad \text{ and } \qquad  V^{(N)}(k) \leq \frac{3}{2} \bigg( \frac{k(k+1)}{2} \nu\bigg),
 \end{align*}
 for $N \geq N_0$. By Tchebychev's inequality for any $k\geq L$  we have
 \begin{align}\label{Equation: Tchebychev inequality yields doubling}
 \mP\left( Z_g^{(N)} \geq \frac{k^2 \mu}{8}\Big| Z_{g-1}^{(N)} =k\right) &\geq
  \mP\left( Z_g^{(N)} \geq \frac{E^{(N)}(k)}{2}\Big| Z_{g-1}^{(N)} =k\right) \\
  &\geq  \mP\left( \Big|Z_g^{(N)} - E^{(N)}(k)\Big| \leq \frac{E^{(N)}(k)}{2}\Big| Z_{g-1}^{(N)} =k\right)\\
  &\geq 1- \frac{4V^{(N)}(k)}{\left(E^{(N)}(k) \right)^2} 
  \geq 1- \frac{48 ( k(k+1) \nu)}{ (k(k+1) \mu)^2} \geq 1- \frac{48 \nu}{ k^2 \mu^2}.
 \end{align}
We choose $f_i(k)= \frac{k^{2^i}\mu^{2^{i}-1}}{8^{2^{i}-1}}$, where 
$f_0(k)=k$. We see that by this choice we have the relation
\begin{equation}\label{Equation: doubling steps i to i+1 DBPC}
    \frac{f_{i-1}^2(k)\mu}{8}=\bigg(\frac{k^{2^{i-1}}\mu^{2^{i-1}-1}}{8^{2^{i-1}-1}}\bigg)^2\frac{\mu}{8}=\bigg(\frac{k^{2^{i}}\mu^{2^{i}-2}}{8^{2^{i}-2}}\bigg)\frac{\mu}{8}=
     \frac{k^{2^i}\mu^{2^{i}-1}}{8^{2^{i}-1}}=f_{i}(k).
\end{equation}
Furthermore, we assumed that
$k>\mu^{-1}(8+\nu)^2$, and therefore  
\begin{equation}\label{Equation: Lower bound in doubling steps DBPC}
    f_{i}(k)=\frac{k^{2^{i}}\mu^{2^{i}-1}}{8^{2^{i}-1}}> \frac{(8+\nu)^{2^{i+1}}\mu^{2^{i}-1}}{\mu^{2^{i}}8^{2^{i}-1}}=
    \frac{(8+\nu)^{2^{i}}(8+\nu)^{2^{i}}}{\mu 8^{2^{i}-1}}
    > \frac{8(8+\nu)^{2^{i}}}{\mu}.
\end{equation}
Now applying \eqref{Equation: Tchebychev inequality yields doubling} and \eqref{Equation: doubling steps i to i+1 DBPC} recursively implies that
\begin{align}\label{lab:KeyInqualitySurv0}
    \mP\left( \bigcap_{i=1}^{M} \left\{ Z_{g+i}^{(N)} > f_{i}(k) \right\} \Big| Z_{g}^{(N)} =k\right)  \geq \prod_{i=1}^{M} \Big(1- \frac{48\nu}{f^2_{i-1}(k) \mu^2}\Big) = \prod_{i=1}^{M} \Big(1- \frac{6\nu}{f_{i}(k) \mu}\Big) 
\end{align}
and by \eqref{Equation: Lower bound in doubling steps DBPC} it follows $f_i(k)\mu>8(8+\nu)^{2^i}$, which yields that
\begin{align}\label{lab:KeyInqualitySurv}
    \mP\left( \bigcap_{i=1}^{M} \left\{ Z_{g+i}^{(N)} > f_{i}(k) \right\} \Big| Z_{g}^{(N)} =k\right)  \geq \prod_{i=1}^{M} \bigg(1-\frac{3}{4(8+\nu)^{2^{i}-1} }\bigg), 
\end{align}
where we used that $\nu (8+\nu)^{-1}<1$.
\end{proof}

\begin{proof}[Proof of Lemma \ref{lem:UniformlyConvergence}]
    By Lemma~\ref{lem:technicalBound} follows that if $N_0$ is large enough such that 
    $b_N>L=\lceil \mu^{-1}(8+\nu)^2 \rceil$ for all $N\geq N_0$,
     \begin{equation*}
        \mP\left( \bigcap_{i=1}^{\infty}\left\{ Z_{g+i}^{(N)} > f_i(b_N) \right\} \Big| Z_{g}^{(N)} =b_N\right)  \geq \prod_{i=1}^{\infty} \Big(1- \frac{6\nu}{f_{i}(b_N) \mu}\Big).
    \end{equation*}
    where $f_i(b_N)= \frac{b_N^{2^i}\mu^{2^{i}-1}}{8^{2^{i}-1}}$ and 
    $f_0(b_N)=b_N$. Without loss of generality we can assume that $b_N$ is monotonically increasing in $N$, which implies that $\log \big(1- \frac{6\nu}{f_{i}(b_N) \mu}\big)$ is monotonically increasing to $0$ as $N\to \infty$. Furthermore, note that for $N$ large enough $\inf_{i\geq 0}\frac{6\nu}{f_{i}(b_N)\mu }\leq \frac{1}{2}$, and thus
    \begin{equation*}
         0\geq \sum_{i=1}^{\infty}\log \Big( 1-\frac{6\nu}{f_{i}(b_N) \mu}\Big)\geq -\sum_{i=1}^{\infty}\frac{12\nu}{f_{i}(b_N) \mu}>-\infty,
    \end{equation*}
    for all $N$ large enough, where we used that $1-x\geq e^{-2x}$ for $x\in[0,\tfrac{1}{2}]$. 
    Now by using continuity of $\exp(\cdot)$ and $\log(\cdot)$ and applying the monotone convergence theorem we obtain that
    \begin{align*}
       \lim_{N\to \infty} \exp\bigg(\log \bigg(\prod_{i=1}^{\infty} \Big(1- \frac{12\nu}{f_{i}(b_N) \mu}\Big)\bigg)\bigg)=&
\exp\bigg(\sum_{i=1}^{\infty} \log \Big( \lim_{N\to \infty} \Big( 1- \frac{12\nu}{f_{i}(b_N) \mu}\Big)\Big)\bigg)=1,
    \end{align*}
    since $1- \frac{12\nu}{f_{i}(b_N) \mu}\to 1$ as $N\to \infty$ for all $i\geq 1$.
\end{proof}

\begin{proof}[Proof of Proposition \ref{Lemma: convergence survival probability DBPC}]
Due to Assumption~\ref{Ass:SequenceDBPC} we have that neither $p_{1,o}^{(N)} \neq 1$ nor $(p^{(N)}_{0,o},p^{(N)}_{1,c})\neq (1,1)$ for $N$ large enough. Due to Lemma 3.1, for an arbitrary $(b_N)_{N\in \N}$ such that $b_N \to \infty$ we have that it exists $A_N$ such that $\mathbb{P}(A_N)=0$ and 
\begin{align}
    \{\forall g \in \mathbb{N}_0: Z_{g}^{(N)}>0\} \backslash A_N \subset \{\exists g \in \mathbb{N}_0: Z_{g}^{(N)} \geq b_N\}.
\end{align}
Then using Lemma \ref{lem:UniformlyConvergence} we obtain that
\begin{align}
    \mathbb{P} \left(\{\exists g \in \mathbb{N}_0: Z_{g}^{(N)} \geq b_N\}\cap \{\exists g \in \mathbb{N}_0: Z_{g}^{(N)}=0\}\right)\to 0.
\end{align}
Consequently using that $\mathbb{P}(\bigcup_{N \in \mathbb{N}}A_N)=0$ it follows that if the limit exists it satisfies  
\begin{align}
\label{Equation: reaching threshold means survival}
    \lim_{N \to \infty} \mathbb{P} \left( \forall g \in \mathbb{N}_0: Z_{g}^{(N)}>0\right)=\lim_{N \to \infty} \mathbb{P} \left( \exists g \in \mathbb{N}_0: Z_{g}^{(N)}\geq b_N\right).
\end{align}

Let $(c_N)_{N\in \mathbb{N}}$ be a sequence with $c_N \to \infty$ and $\frac{b_{N}}{c_{N}} \to\infty$. In order to show, if the limit exists, that
\begin{align}
    \lim_{N \to \infty} \mathbb{P} \left( \exists g \in \mathbb{N}_0: Z_{g}^{(N)}\geq b_N\right)=\lim_{N \to \infty} \mathbb{P} \left( \exists g \in \mathbb{N}_0: \overline{Z}_{g}^{(N)}\geq b_N\right),
\end{align}
it remains to show that
\begin{align}\label{Eq: current size drives total size}
    \mathbb{P} \left( \exists g \in \mathbb{N}_{0}: \overline{Z}_{g}^{(N)} \geq b_N \text{ and } \forall i \leq g, Z_{i}^{(N)} \leq c_N\right) \to 0.
\end{align}
In particular if $\overline{\tau}^{(N)}_{b_N}:= \inf\{g \in \mathbb{N}_{0}: \overline{Z}^{(N)}_g \geq b_N\}$ we have 
\begin{align}
    \mathbb{P} \left( \exists g \in \mathbb{N}_{0}: \overline{Z}_{g}^{(N)} \geq b_N \text{ and } \forall i \leq g, Z_{i}^{(N)} \leq c_N\right) \leq \mathbb{P} \left( \left\lceil\frac{b_N}{c_N}\right\rceil\leq \overline{\tau}^{(N)}_{b_N}<\infty \right).
\end{align}
But according to Proposition \ref{GrowFast} it follows that $\mathbb{P} \left( \overline{\tau}^{(N)}_{b_N} \leq C\log\left(\log\left(b_N\right)\right)\vert \overline{\tau}^{(N)}_{b_N} <\infty\right) \to 1$. So in particular taking $c_N \to \infty$ such that $\log(\log(b_N))=o(\frac{b_N}{c_N})$ implies that  
$\mathbb{P} \left( \lceil\frac{b_N}{c_N}\rceil\leq \overline{\tau}^{(N)}_{b_N}<\infty \right) \to 0,$ which gives \eqref{Eq: current size drives total size}.\\
To conclude the proof it only remains to show that 
\begin{align}
    \lim_{N \to \infty} \mathbb{P} \left( \exists g \in \mathbb{N}_{0}: \overline{Z}_{g}^{(N)} \geq b_N\right)=\pi.
\end{align}
Note that $\sum_{k=0}^{\infty} \vert p_{k,o}^{(N)}-p_{k,o}\vert\leq 2$, and thus by dominated convergences follows that
\begin{equation*}
    \lim_{N\to \infty}\sum_{k=0}^{\infty} \vert p_{k,o}^{(N)}-p_{k,o}\vert=0.
\end{equation*}
Analogously follows that $\lim_{N\to \infty}\sum_{k=0}^{\infty} \vert p_{k,c}^{(N)}-p_{k,c}\vert=0$.  From this follows that 
for a given sequence $(K_N)_{N\in\N}\subset \N$ such that $K_N\to \infty$ as $N\to \infty$ we find a sequence $\left(\varepsilon_N\right)_{N \in \mathbb{N}}$ with $\varepsilon_N\to 0$ as $N\to \infty$ and with 
\begin{align}\label{eq:Bound on offspring distribution 1}
    &\max \left( \sum_{k=0}^{K_N} \vert p_{k,o}^{(N)}-p_{k,o}\vert, \sum_{k=0}^{K_N}\vert p_{k,c}^{(N)}-p_{k,c}\vert \right) \leq \varepsilon_N, \\
   & \max \left( \sum_{k=K_N+1}^{\infty} p_{k,o}, \sum_{k=K_N+1}^{\infty} p_{k,c}\right) \leq \varepsilon_N.
\end{align}
Note that this implies
\begin{align*}
    \sum_{k=K_N+1}^{\infty} p^{(N)}_{k,o}&\leq \bigg|\sum_{k=K_N+1}^{\infty}  (p_{k,o}^{(N)}-p_{k,o})\bigg|+\sum_{k=K_N+1}^{\infty} p_{k,o}\leq \bigg|\sum_{k=0}^{K_N} (p_{k,o}^{(N)}-p_{k,o})\bigg|+\varepsilon_N\leq 2\varepsilon_N.
\end{align*}
By the exact same calculation we get the same bound for the sum of $p_{k,c}^{(N)}$ from $K_N+1$ to $\infty$ and this yields that
\begin{equation}\label{eq:Bound on offspring distribution 2}
    \max \bigg( \sum_{k=K_N+1}^{\infty} p^{(N)}_{k,o}, \sum_{k=K_N+1}^{\infty} p^{(N)}_{k,c}\bigg) \leq 2\varepsilon_N.
\end{equation}
We know by assumption that $\varepsilon_N \to 0$.
Consider now a sequence $\left(e_N\right)_{N \in \mathbb{N}}$ such that $e_N \to \infty$ and $\varepsilon_N e_N^{2} \to 0$. The first step is to prove that
\begin{align}
\label{Equation: total size DBPC}
     \lim_{N \to \infty} \mathbb{P} \left( \exists g \in \mathbb{N}_{0}: \overline{Z}_{g}^{(N)} \geq e_N\right)=\pi.
 \end{align}

We start by showing that whp. the sequence of DBPC $\left(\ZZ^{(N)}\right)_{N \in \mathbb{N}}$ and the limiting DBPC $\ZZ$ can be exactly coupled until their total size reaches the level $b_N$ or they both die out. Introduce the stopping time of the first generation that the total size of $\ZZ$ reaches the level $e_N$ or that it dies out as
\begin{equation}
    \overline{\tau}_{e_N,0}:=\inf \left\{g \in \mathbb{N}_0: \overline{Z}_{g} \geq e_N \text{ or } Z_{g}=0\right\}.
\end{equation}
By definition we have that almost surely $\overline{Z}_{\overline{\tau}_{e_N,0}-1} < e_N$ which means that in order to make an exact coupling between $\ZZ^{(N)}$ and $\ZZ$ until generation $\overline{\tau}_{e_N,0}$, there are at most $e_N^2$ offspring and cooperation independent random variables to couple. Till this point we have not specified the joint distribution of the offspring and cooperation random variables $(X_i,Y_{j,k})_{i\in\N,j<k}$ and $(X^{(N)}_i,Y^{(N)}_{j,k})_{i\in\N,j<k}$. We couple them in such a way that $\mP\big(X_i\neq X^{(N)}_i\big)$ and $\mP\big(Y_{j,k}\neq Y^{(N)}_{j,k}\big)$ are minimized. For a single random variable this can be done for each pair recursively via the maximal coupling, see Theorem~2.9 in \cite{van2016random}, such that
\begin{equation*}
    \mP\big(X_i\neq X^{(N)}_i\big)=\frac{1}{2}\sum_{k=0}^{\infty} \vert p_{k,o}^{(N)}-p_{k,o}\vert\leq 2\varepsilon_N\quad \text{ and } \quad \mP\big(Y_{j,k}\neq Y^{(N)}_{j,k}\big)=\frac{1}{2}\sum_{k=0}^{K_N} \vert p_{k,o}^{(N)}-p_{k,o}\vert\leq2\varepsilon_N,
\end{equation*}
where we used the bounds from \eqref{eq:Bound on offspring distribution 1} and \eqref{eq:Bound on offspring distribution 2}. Since these are families of independent random variables and also the offspring and cooperation random variables across different generations are independent, the probability that the $e_N^2$ relevant offspring and cooperation independent random variables are equal is lower bounded by $(1- 2 \varepsilon_N)^{ e_N^2}\rightarrow 1$ by the
choice of $(e_N)$. In summary we have 
\begin{align}\label{eq:FiniteLevelCoupling}
    \mathbb{P} \big( \overline{Z}^{(N)}_g =\overline{Z}_g, \forall g \leq \overline{\tau}_{e_N,0}\big)\geq (1- 2 \varepsilon_N)^{ e_N^2} \to 1
\end{align}
as $N\to \infty$. Let us now define the event $C_N:=\{\overline{Z}^{(N)}_g =\overline{Z}_g, \forall g \leq \overline{\tau}_{e_N,0}\}$. We see that by monotonicity and Lemma~\ref{lem:ReachingLevelToSurvival} that
\begin{align*}
    \mP(\{\exists g \in \mathbb{N}_{0}: \overline{Z}_{g} \geq e_N\}\cap C_N)\leq  \mP(\exists g \in \mathbb{N}_{0}: \overline{Z}_{g} \geq e_N)\to \pi
\end{align*}
as $N\to \infty$. On the other hand, by monotonicity and Equation \eqref{eq:FiniteLevelCoupling} we see that
\begin{align*}
    \mP(\{\exists g \in \mathbb{N}_{0}: \overline{Z}_{g} \geq e_N\}\cap C_N^c)\leq \mP(\exists g\leq \overline{\tau}_{e_N,0} : \overline{Z}^{(N)}_g \neq \overline{Z}_g)\to 0
\end{align*}
as $N\to \infty$. This yields that
\begin{align*}
    \lim_{N\to \infty}\mP(\{\exists g \in \mathbb{N}_{0}: \overline{Z}_{g} \geq e_N\}\cap C_N)= \lim_{N\to \infty} \mP(\exists g \in \mathbb{N}_{0}: \overline{Z}_{g} \geq e_N)=\pi.
\end{align*}
But $C_N$ states that the $Z$ and $Z^{(N)}$ are coupled until $\overline{\tau}_{e_N,0}$, and therefore we also know that
\begin{align*}
     \mP(\{\exists g \in \mathbb{N}_{0}: \overline{Z}_{g}\geq e_N\}\cap C_N)=\mP(\{\exists g \in \mathbb{N}_{0}: \overline{Z}_{g}^{(N)} \geq e_N\}\cap C_N)
\end{align*}
for all $N>0$. 
But this equality already implies Equation \eqref{Equation: total size DBPC}, i.e.\
\begin{align*}
     \lim_{N\to \infty}\mP(\{\exists g \in \mathbb{N}_{0}: \overline{Z}_{g}^{(N)} \geq e_N\})=\pi.
\end{align*}

This concludes the proof since we have shown it for $b_N=e_N$ and since we have shown it for one specific choice it follows also for an arbitrary sequence $(b_N)_{N\in \N}$ because of the extinction-explosion principle shown in Lemma~\ref{Lemma: extinction-explosion principle DBPC}. 
\end{proof}

\begin{Lem}\label{Lemma: Controll exp-Term}
Consider sequence $(D_N), (V_N), (m_N), (f_N), (g_N), (h_N)$ such that $V_N \sim a\sqrt{D_N}$, $0\leq h_N\leq f_N$, $g_N, m_N\geq 0$ and for $k_N:=\max\{m_N, f_N, g_N \}$ it holds $\frac{k_N^4 V_N^3}{D_N^2} \in o(1)$. Then  
\begin{align*}
\frac{(D_N - f_N)!}{(D_N-f_N - (m_N V_N- g_N))! (D_N- h_N)^{m_N V_N- g_N} } & \geq \exp\left(-\frac{(m_N V_N- g_N)^2}{2 D_N}\right) \exp\left(- \frac{k_N^4 V_N^3}{D_N^2} \right) \\ &\geq \exp\left(-\frac{(m_N V_N)^2}{2 D_N}\right) \exp\left(- \frac{k_N^4 V_N^3}{D_N^2} \right) .
\end{align*}
On the other hand we have 
\begin{align*}
\frac{(D_N - f_N)!}{(D_N-f_N - (m_N V_N- g_N))! (D_N- h_N)^{m_N V_N- g_N} } & \leq \exp\left(-\frac{(m_N V_N- g_N)^2}{2 D_N}\right) \exp\left(\frac{m_NV_N}{D_N}\right).
\end{align*}
\end{Lem}

\begin{proof}[Proof of Lemma \ref{Lemma: Controll exp-Term}]
For completeness we show the inequality 
\[ 1-x \geq \exp(-x) \exp(-x^2)\] for $x\in [0, \tfrac{1}{2}]$ first.
We have
 \[1- x + x^2/2 \geq \exp(-x),\] so
 \[
 (1-x)\left(1 + \frac{x^2}{2(1-x)}\right) \geq \exp(-x)\]
 which is equivalent to
\[1-x  \geq \exp(-x) \frac{1}{\left(1 + \frac{x^2}{2(1-x)}\right)} .\]
We have 
\[1 + \frac{x^2}{2(1-x)} \leq 1+ x^2 < \exp(x^2)\]
 which yields
\[1-x  \geq \exp(-x)\exp(-x^2).\]

We use this inequality to estimate 

\begin{align*}  &\frac{(D_N - f_N)!}{(D_N-f_N - (m_N V_N- g_N))! (D_N- h_N)^{m_N V_N- g_N} }  \\
& = \left(1 + \frac{h_N -f_N}{D_N - h_N}\right)\cdot   ... \cdot 
\left(1 + \frac{h_N -f_N}{D_N - h_N} - \frac{m_N V_N - g_N -1}{D_N - h_N}\right)\\
&\geq  
\exp\left(\frac{-(f_N- h_N)(m_N V_N - g_N)}{D_N-h_N} \right) 
\exp\left(-\sum_{i=1}^{m_N V_N- g_N-1} \frac{i}{D_N - h_N} \right) \\ & \quad
 \cdot \exp\left( -\left(\frac{f_N - h_N - m_N V_N }{D_N - h_N}\right)^2 m_N V_N  \right)\\
 &\geq  \exp\left(-\frac{\left(m_N V_N -g_N  \right)^2}{2 D_N} \right) \exp\left(- \frac{k_N^4 V_N^3}{D_N^2} \right),
\end{align*}
using $\frac{1}{1-x} \leq 1+2x$ for $x \in [0,1/2]$ and $\frac{V_N}{D_N}\leq \frac{k_N V_N^3}{D_N^2}$.
Then we also have that 

\begin{align*}  &\frac{(D_N - f_N)!}{(D_N-f_N - (m_N V_N- g_N))! (D_N- h_N)^{m_N V_N- g_N} }  \\
&\leq  
\exp\left(\frac{-(f_N- h_N)(m_N V_N - g_N)}{D_N-h_N} \right) 
\exp\left(-\sum_{i=1}^{m_N V_N- g_N-1} \frac{i}{D_N - h_N} \right) \\
 &\leq  \exp\left(-\frac{\left(m_N V_N -g_N-1  \right)^2}{2 D_N} \right) \\
 & \leq \exp\left(-\frac{\left(m_N V_N -g_N  \right)^2}{2 D_N}\right) \exp \left( \frac{m_NV_N}{D_N}\right). 
\end{align*}

\end{proof}

\begin{proof}[Proof of Lemma \ref{lab: balls and boxes experiment}]
We have that 
\begin{align}
    \mathbb{P} \left(C_k^{(h_N')}\right)&= \frac{\prod_{i=1}^{k}\binom{m'_N V'_N-2 (i-1)}{2}}{k! (D'_N-m'_N)^{k}}
    \frac{(D'_N-h'_N)!}{(D'_N-h'_N-k)!(D'_N-m'_N)^k}\\
    &\cdot\frac{(D'_N-m'_N-k)!}{[(D'_N-m'_N-k)-(m'_N V'_N-2k)]!(D'_N-m'_N)^{m'_N V'_N-2k}}.
\end{align}
In particular applying Lemma \ref{Lemma: Controll exp-Term} with $f_N=m'_N+k$, $g_N=2k$, and $h_N=m'_N$, one obtains that 
\begin{align}
    \mathbb{P} \left(C_k^{(h_N')}\right) &\leq \left( \frac{m'_NV'_N}{2D'_N}\right)^k \frac{1}{k!} \exp\left( \ell'_N \frac{m'_N}{D'_N}\right) \exp \left( -\frac{(m'_NV'_N-2\ell'_N)^2}{2D'_N}\right) \exp \left( \frac{m'_NV'_N}{D'_N}\right) \\
    &\leq\left( \frac{m'_NV'_N}{2D'_N}\right)^k \frac{1}{k!}\exp \left( -\frac{(m'_NV'_N-2\ell'_N)^2}{2D'_N}\right) \exp \left( \frac{{\ell'_N}^2 V'_N}{D'_N}\right).
\end{align}
 \\
Applying again Lemma \ref{Lemma: Controll exp-Term} with $f_N=m'_N+k$, $g_N=2k$ and $h_N=m'_N$ we obtain 
\begin{align}
   \mathbb{P} \left(C_k^{(h_N')}\right)& \geq \left( \frac{\left(m'_NV'_N-2\ell'_N\right)^2}{2 D'_N}\right)^k \frac{1}{k!} \exp \left( -4\frac{{\ell'}_N^2}{D'_N-m'_N}\right) \exp \left( -\frac{(m'_N V'_N)^2}{2D'_N}\right) \exp \left( -\frac{{\ell'}_N^4 {V'}_N^3}{{D'}_N^2}\right) \\
    & \geq \left( \frac{\left(m'_NV'_N-2\ell'_N\right)^2}{2 D'_N}\right)^k \frac{1}{k!} \exp \left( -\frac{(m'_N V'_N-2\ell'_N)^2}{2D'_N}\right)\exp \left( -\frac{{\ell'}_N^5{V'}_N^3}{{D'}_N^2}\right).
\end{align}
\end{proof}

\begin{proof}[Proof of Lemma \ref{Lem: experiment boxes and balls}]
Denote by $S_N=D_N-\varphi_1(N)-\varphi_2(N)$, $b_N=D_N-\varphi_1(N)$, and $h(N):=H(N)(v_N-\varphi_3(N))$. Introduce for all $i \leq b_N$ the random variable $G_i^{(N)}$ which counts the number of balls in box $i$. We have 
\begin{align}
    &\mathbb{P}\left(G_1^{(N)} \leq 1\right)=\mathbb{P}\left(G_1^{(N)}=0\right)+\mathbb{P}\left(G_1^{(N)}=1\right) \\
    &=\left(1-\frac{1}{b_N}\right)^{h(N)}+\frac{h(N)}{b_N}\left(1-\frac{1}{b_N}\right)^{h(N)-1} \\
    &=\exp \left(h(N) \log\left(1-\frac{1}{b_N}\right)\right) \left[1+\frac{h(N)}{b_N}\frac{1}{1-\frac{1}{b_N}}\right]\\
    &= \left[1-\frac{h(N)}{b_N}+\frac{1}{2}\left(\frac{h(N)}{b_N}\right)^2+\mathcal{O}\left(\left(\frac{h(N)}{b_N}\right)^3\right)\right]\left(1-\mathcal{O}\left(\frac{h(N)}{b_N^2}\right)\right) \\
    &\hspace{2cm}\cdot\left[1+\frac{h(N)}{b_N}+\mathcal{O}\left(\frac{h(N)}{b_N^2}\right)\right] \\
    &=1-\frac{1}{2}\left(\frac{h(N)}{b_N}\right)^2+\mathcal{O}\left(\frac{h(N)}{b_N^{2}}\right).
\end{align}

Combining this previous computation and the fact that $S_N\frac{h^{2}(N)}{b_{N}^{2}}=\Theta(H^{2}(N))$ gives that  
\begin{align}
    \mathbb{E}\left[G^{(N)}\right]=\mathbb{E}\left[\sum_{i=1}^{S_N}\1_{\left\{G_i^{(N)}\geq 2\right\}}\right]=S_N \left(1-\mathbb{P}\left(G_{1}^{(N)}\leq 1\right)\right)=\Theta\left(H^{2}(N)\right). 
\end{align}
By using Markov Inequality we obtain that 
\begin{align}
    \mathbb{P} \left( G^{(N)} \geq H^{2}(N)f_2(N)\right) \leq \frac{\mathbb{E}[G^{(N)}]}{H^{2}(N)f_2(N)}=\Theta \left(\frac{1}{f_2(N)}\right),
\end{align}
which gives exactly \eqref{Eq: max number pair of balls}. \\
Due to the scaling of $H(N)$ it could happen that some boxes have more than 3 balls in it. In order to deal with such a fact introduce 
\begin{align}
    \ell:=\inf \left\{i \geq 2: H(N)^{i+1}=o\left(\sqrt{D_N}^{i-1}\right)\right\}.
\end{align} 
Now we have 
\begin{align}
    \mathbb{P}\left( \exists i \leq b_N, G^{(N)}_{i} \geq \ell+1\right) \leq b_N \frac{h(N)^{\ell+1}}{b_N^{\ell+1}}=\Theta \left( \frac{H(N)^{\ell+1}V_N^{\ell+1}}{D_N^{\ell}}\right)= \Theta \left( \frac{H(N)^{\ell+1}}{\sqrt{D_N}^{\ell-1}}\right) \to 0,
\end{align}
by definition of $\ell$. Then under the event $\left\{\forall i \leq b_N, G_{i}^{(N)}\leq \ell\right\}$, getting $G^{(N)}=k$ is obtained via different scenarios: any $k$ packs of balls of a number between $2$ and $\ell$ and all the other remaining balls are alone. This number of scenarios is easily upper bounded by $k^{\ell}$, because for all $2 \leq j \leq \ell$ the number of boxes getting exactly $j$ balls is upper bounded by $k$. Choosing the packs of balls and where they are going is upper bounded by $\frac{H^{k\ell}(N)V_{N}^{k\ell}}{b_N^{k}}$. The number of remaining balls to be placed on different boxes is lower bounded by $h(N)-k\ell$. At the end we get that  
\begin{align} \label{lab: upper control final phase}
    \mathbb{P}\left( G^{(N)} \leq \frac{H^{2}(N)}{f_1(N)} \Big\vert \{\forall i \leq b_N, G_i \leq \ell\}\right)\leq \sum_{k=0}^{H^{2}(N)f_1^{-1}(N)}k^{\ell} H^{k \ell}(N) V_{N}^{k\ell}\prod_{i=0}^{h(N)-k\ell}\left(1-\frac{i}{b_N}\right).
\end{align}
Also using that $\prod_{i=0}^{h(N)-k\ell}\left(1-\frac{i}{b_N}\right)=\exp \left(\sum_{i=0}^{h(N)-k\ell} \log\left(1-\frac{i}{b_N}\right)\right)$ and that $\log(1-x)\leq -x$ for all $x \in [0,1)$, one can show, using that $h(N)=o(b_N)$, that for all $k \leq H^{2}(N)f_{1}^{-1}(N)$ that 
\begin{align}
\label{Equation: bound factorial term}
    \prod_{i=0}^{h(N)-k\ell} \left(1-\frac{i}{b_N}\right)\leq \exp \left(-\sum_{i=0}^{h(N)-k\ell}\frac{i}{b_N}\right) \leq \exp \left(-\sum_{i=0}^{h(N)-H^{2}(N)\ell}\frac{i}{b_N}\right).
\end{align}
Because $\sum_{i=0}^{h(N)-H^{2}(N) \ell}\frac{i}{b_N}\sim \frac{H^{2}(N)v^{2}_{N}}{2b_{N}}\sim H^{2}(N)\frac{c^{2}}{2}$, it follows that for $N$ large enough we have 
\begin{align}
    \exp \left(-\sum_{i=0}^{h(N)-H^{2}(N) \ell}\frac{i}{b_N}\right) \leq \exp \left( -H^{2}(N)\frac{a^{2}}{4}\right). 
\end{align}
Then using the natural bound $k^{\ell} \left(H(N)V_N\right)^{k \ell} \leq H^{2\ell}(N) \left(H(N)V_N\right)^{H^{2}(N)f_{1}^{-1}(N)\ell}$ for all $k \leq H^{2}(N)f_{1}^{-1}(N)$, we get that for $N$ large enough
\begin{align}
    \eqref{lab: upper control final phase} &\leq \exp \left( -H^{2}(N)\frac{a^{2}}{4}\right) H^{2\left(\ell+1\right)}(N) \left(H(N)V_N\right)^{H^{2}(N)f_{1}^{-1}(N)\ell} \\
    &= \exp \left(-\frac{H^{2}(N)}{4}\left[a^{2}-4 \ell f_{1}^{-1}(N)\log\left(H(N)V_{N}\right)\right]\right)H^{2\left(\ell+1\right)}(N) \\
    & \leq \exp \left(-\frac{H^{2}(N)a^{2}}{8}\right),
\end{align}
where for the last inequality we use that $\log(H(N)V_{N})=o(f_{1}(N))$. \\
Finally to conclude the proof we have 
\begin{align}
    \mathbb{P}\left( G^{(N)} \leq \frac{H^{2}(N)}{f_1(N)} \right) &\leq \mathbb{P}\left( \exists i \leq b_N, G^{(N)}_{i} \geq \ell+1\right) +\mathbb{P}\left( G^{(N)} \leq \frac{H^{2}(N)}{f_1(N)} \Big\vert \{\forall i \leq b_N, G_i \leq \ell\}\right) \\
    & \leq \Theta \left( \frac{H(N)^{\ell+1}}{\sqrt{D_N}^{\ell-1}}\right).
\end{align}

\end{proof}

\begin{Lem}\label{Lemma: InvasionTildeProcess}
Let $0<\alpha< \beta$, $a>0$ and consider the sequence of processes $\big(\widetilde{\mathcal{S}}^{(N,a)},\widetilde{\mathcal{J}}^{(N,a)},\widetilde{\mathcal{R}}^{(N,a)}\big)_N$ defined in Section \ref{Section: Proof of Theorem 2.7 (iii)}. It holds
\begin{align*}
\lim_{N\rightarrow \infty}\mathbf{P}(\{
\exists g\geq 0:|
\widetilde{\mathcal{R}}_g^{(N,a)}|
\geq N^\alpha\})
\geq \varphi_a,
\end{align*}
 where $\varphi_a$ denotes the survival probability of a Galton-Watson process with Poi$\left(\tfrac{a^2}{2}\right)$ offspring distribution.
\end{Lem}

In the proof of Lemma \ref{Lemma: InvasionTildeProcess}
we will couple  $\big(\widetilde{\mathcal{S}}^{(N,a)},\widetilde{\mathcal{J}}^{(N,a)},\widetilde{\mathcal{R}}^{(N,a)}\big)_N$
 with the Galton-Watson process defined next.

\begin{Def}(Lower Galton-Watson Process)
\label{lower branching process} \\
Let $0<\delta<\frac{1}{2}$ and $(b_N)_{N \in \N}$ be a sequence with $b_{N} \rightarrow \infty$ and $a_N \in o\left( \sqrt{D_N}\right)$. Furthermore assume $(\theta_N)_{N \in \N}$ is a $[0,1]$-valued sequence with
$\theta_{N} \rightarrow 0$. 
Let $\mathbf{X}_{l}^{(N)}=\big( X^{(N)}_{g,l}\big)_{g \in \N_0}$ be a Galton-Watson process with mixed binomial offspring distribution $\text{Bin}\big(\widetilde{X}^{(N)},1-\theta_{N}\big)$, where the probability weights $\big(\widetilde{p}_k^{(N)}\big)_{k \in \N_0}$ of $\widetilde{X}^{(N)}$ are
for all   $ 1 \leq j \leq b_{N}$
\begin{linenomath*}
\begin{align}
\label{first offspring distribution LGWP without removal 1}
    &\widetilde{p}_j^{(N)}:=\left(\dfrac{(V^{(a)}_{N}-2b_{N})^{2}}{2D_{N}} \right)^{j}\frac{1}{j!}\exp \left(-\frac{\left(V^{(a)}_{N}\right)^{2}}{2D_{N}}\right) \left(1-\frac{1}{D_N^{\delta}} \right),
    \end{align}
    \end{linenomath*}
    and 
    \begin{linenomath*}
    \begin{equation}\label{first offspring distribution LGWP without removal 2}
    \widetilde{p}_0^{(N)}:=1-\sum_{j=1}^{b_N}\widetilde{p}_j^{(N)}.
    \end{equation}
\end{linenomath*}
Denote by $\Phi_{l}^{(N)}$ the generating function of the offspring distribution $\big(p_{k,l}^{(N)}\big)_{k \in \N_{0}}$
of $\mathbf{X}_{l}^{(N)}$, and by $\overline{X}^{(N)}_{g,l} := \sum_{i=0}^{g} X_{i,l}^{(N)}$ the total size of the Galton-Watson process until generation $g$ and $\overline{\mathbf{X}_{l}} = \big(\overline{X}_{g,l}^{(N)}\big)_{g \in \mathbb{N}_0} $ the corresponding process.  
\end{Def}

\begin{proof}[Proof of Lemma \ref{Lemma: InvasionTildeProcess}]

Introduce $\overline{\sigma}^{(N)}_{N^{\alpha}}=\inf \left\{ g \in \mathbb{N}_{0}: \vert \widetilde{\mathcal{J}}^{(N,a)}_{g} \cup \widetilde{\mathcal{R}}^{(N,a)}_{g}\vert \geq N^{\alpha}\right\}$. We are going to show that 
\begin{align}\label{lab: coupling below GWP}
    \mP \left( \vert \widetilde{\mathcal{J}}^{(N,a)}_{g} \cup \widetilde{\mathcal{R}}^{(N,a)}_{g} \vert \geq \overline{X}^{(N)}_{g,l} \, \forall \, g \leq \overline{\sigma}^{(N)}_{N^{\alpha}}\right) \to 1
\end{align}
for the process $\overline{\bf{X}}^{(N)}_{l}$ defined in Definition \ref{lower branching process}.

Consider $D_N-1$ boxes, assume $V_N$ many balls are put uniformly at random into the boxes. Denote by $C^{(N)}_{j}$ the event that exactly $j$ boxes contain at least 2 balls.
One can show using similar calculations as in the proof of
Proposition 3.5 of \cite{BrouardEtAl2022} that 
\begin{align}
    \widetilde{p}^{(N)}_j \leq \mP\left(C^{(N)}_{j}\right), \forall 1 \leq j \leq a_{N},
\end{align}
  This means that whenever a vertex $x$ gets infected one can estimate from below how many of its neighbors are visited by at least $2$ of its $V_N$ parasites, which we call a \textit{pack of parasites}, by the corresponding number of offspring in the Galton-Watson process $\mathbf{X}_{l}^{(N)}$, since $\widetilde{p}^{(N)}_{0}=1-\sum_{i=1}^{b_N}\Tilde{p}^{(N)}_j$.

However, in the process $\big(\widetilde{\mathcal{S}}^{(N,a)},\widetilde{\mathcal{J}}^{(N,a)},\widetilde{\mathcal{R}}^{(N,a)}\big)$ ``ghost'' infections may occur, when a) an already empty vertex is attacked by at least $2$ parasites coming from the same infected vertex, or when b) a vertex is attacked by at least two packs of parasites coming from different vertices.

We show next that each pack of parasites of size at least $2$ generated by an infected vertex before generation $\overline{\sigma}^{(N)}_{N^{\alpha}}$ is involved in one of the events a) or b) (independently of the other packs of parasites) with probability at most $\theta_N$. Consequently, by removing packs of parasites with probability $\theta_N$ the number of new infection generated by an infected hosts can with high probability be bounded from below by the number of offspring drawn according to the distribution with weights $(p_{k, l}^{(N)})_{k\in \mathbb{N}_0}$ from Definition \ref{lower branching process} for any generation $n \leq\overline{\sigma}^{(N)}_{N^{\alpha}}$.

Now we upper bound the probabilities of the events a) and b).
\begin{itemize}
    \item[a)] Before generation $\overline{\sigma}^{(N)}_{N^{\alpha}}$ the probability that a pack of parasites goes to an empty vertex is bounded from above by $\frac{N^{\alpha}}{D_N-1} =\Theta \left(\frac{1}{N^{\beta-\alpha}}\right)$.
    \item[b)] Before generation $\overline{\sigma}^{(N)}_{N^{\alpha}}$, the number of empty vertices in the graph is smaller than $N^{\alpha}$. The probability that two packs of parasites coming from 2 different vertices attack the same vertex is $\frac{1}{D_N-2}=\Theta\left( \frac{1}{N^{\beta}}\right)$. Using Markov inequality one can show that the total number of packs of parasites generated before generation  $\overline{\sigma}^{(N)}_{N^{\alpha}}$ is with high probability bounded by $N^{\alpha}\log(N)$, as in Lemma 4.8 in \cite{BrouardEtAl2022}. Hence, each pack of parasites is involved in an event of type $b)$ with probability at most $N^{\alpha}\log(N) \cdot \frac{1}{D_N-2}=\Theta \left(\frac{\log(N)}{N^{\beta-\alpha}}\right)$.
\end{itemize}
In summary, $\theta_N= 2 \frac{N^\alpha \log(N)}{D_N-2} =\Theta \left(\frac{\log(N)}{N^{\beta-\alpha}}\right)$ yields an upper bound on the probability that a pack of parasites is involved in one of the events of type a) or b). Since $\alpha < \beta$ we have $\theta_N \in o(1)$.\\
Then taking the Galton Watson process $\mathbf{X}_{l}^{(N)}$ of Definition \ref{lower branching process} with $\theta=2 \frac{N^\alpha \log(N)}{D_N-2}$ we have exactly proven \eqref{lab: coupling below GWP}.\\ \\
As in the proof of Proposition 3.3 of \cite{BrouardEtAl2022} one can show the uniform convergence of $\Phi_{l}^{(N)}$ to the generating function of a $\text{Poi}\big(\frac{a^2}{2}\big)$ distribution, such that applying Lemma 3.7 of \cite{BrouardEtAl2022} gives that 
\begin{equation} \label{lab: GWP reaching size}
    \lim_{N \to \infty} \mP \left( \exists g \in \mathbb{N}_{0}: \overline{X}^{(N)}_{g,l} \geq N^{\alpha}\right)=\varphi_a.
\end{equation}

Combining \eqref{lab: coupling below GWP} and \eqref{lab: GWP reaching size} we get that 
\begin{align}
    \liminf_{N \to \infty} \mP \left( \exists g \in \mathbb{N}_{0}: |
\widetilde{\mathcal{R}}_g^{(N,a)}|
\geq N^\alpha\right) \geq \varphi_a.
\end{align}
\end{proof}

For the proof of Lemma \ref{Lemma: moderate deviation number of vertices}  we have to control moderate deviations of Poisson distributed random variables. 
 
\begin{Lem} \label{Lemma: moderate deviation number of vertices}
Let $0<\gamma <1$, $N\in \N$ and $X \sim \text{Poi}(N^\gamma)$ then for any $0<\varepsilon < \gamma/2$ it holds that
    \begin{align} \label{Eq:ModerateDeviation}
        \lim_{N\to \infty}\frac{-2}{N^{2\varepsilon}}\log(\mP(|X -N^\gamma| > N^{\frac{\gamma}{2} +\varepsilon})) =1
    \end{align}
\end{Lem}

The proof of Lemma \ref{Remark: uniform bound number vertices} we base on the following lemma.

\begin{Lem}\label{Proposition: moderate deviation for Poisson process}
    Let $X=(X_t)$ be a Poisson process with intensity $1$ on $[0,\infty)$ and $a(t)$ be a function such that $\tfrac{t}{a(t)}\to \infty$ and $\tfrac{a(t)}{\sqrt{t}}\to \infty$. Then for every Borel-set $B\subset\R$ it holds that
    \begin{align*}
    \limsup_{t\to \infty} \frac{t}{a(t)^2}\log\Big(\mP\Big(\frac{1}{a(t)}(X_t-t)\in B\Big)\Big)&\leq -\inf_{x\in \overline{B}} \frac{x^2}{2} \quad \text{ and }\\
    \liminf_{t\to \infty} \frac{t}{a(t)^2}\log\Big(\mP\Big(\frac{1}{a(t)}(X_t-t)\in B\Big)\Big)&\geq -\inf_{x\in \mathring{B}}\frac{x^2}{2},
\end{align*}
where $\mathring{B}$ denote the interior and $\overline{B}$ the closure of $B$.
\end{Lem}

\begin{proof}[Proof of Lemma \ref{Proposition: moderate deviation for Poisson process}]
    This is a direct consequence of Theorem 1.1' found in \cite{gao1996moderate}. We will now check the conditions of this theorem. First denote by $Y_t=X_t-t$ the compensated Poisson process, note that $(Y_t)_{t\geq 0}$ is a martingale. Let $\delta>0$ not it holds that 
    \begin{equation*}
        \E\big[\exp\big(\delta\sup_{s\leq t\leq s+1}|Y_t-Y_s|\big)\big|\sigma(Y_u:u\leq s)\big] = \E\big[\exp\big(\delta\sup_{0\leq t\leq 1}|Y_t|\big)\big],
    \end{equation*}
    where we used that the Poisson process has independent and stationary increments and that $Y_0=0$. Furthermore it holds that
    \begin{equation*}
        \E\big[\exp\big(\delta\sup_{0\leq t\leq 1}|Y_t|\big)\big]
        \leq \E\big[\exp\big(\delta X_1\big)\big]e^{1}=\exp(\exp(\delta)), 
    \end{equation*}
    where we used that $X_t$ is a monotone process and that the moment generating function of a Poisson distribution with parameter $1$ is given through $t\mapsto \exp(e^t-1)$. This provides $(A1)'$. The second condition $(A2)'$ follows from the fact that 
    \begin{equation*}
        \E\Big[\frac{1}{t}(Y_{s+t}-Y_s)^2-1\Big|\sigma(Y_u:u\leq s)\Big]= \E\Big[\frac{1}{t}Y_t^2\Big]-1=0,
    \end{equation*}
    where we again used again that the process has stationary and independent increments. Now the claimed moderate deviation principle follows from Theorem 1.1' in \cite{gao1996moderate}.
\end{proof}

\begin{proof}[Proof of Lemma \ref{Lemma: moderate deviation number of vertices}]
    This follows directly from Proposition~\ref{Proposition: moderate deviation for Poisson process} choose $a(t)=t^{\frac{1}{2}+\frac{\varepsilon}{\gamma}}$, $B=[0,1]^{c}$ and consider the subsequent $(X_{N^{\gamma}}-N^{\gamma})_{N\geq 0}$ instead of $(X_t-t)_{t\geq 0}$. Then plugging the choices in we get that
    \begin{equation*}
        \lim_{N\to \infty}\frac{1}{N^{2\varepsilon}}\log(\mP(|X -N^\gamma| > N^{\frac{\gamma}{2} +\varepsilon})) =-\frac{1}{2}.
    \end{equation*}
    which provides the claim.
\end{proof}

\end{document}